\documentclass[a4paper,12pt]{article}
\usepackage[left=3cm, right=3cm, top=2cm, bottom=2cm]{geometry}

\usepackage{amsmath, amssymb, amsthm, bm}

\usepackage{anyfontsize}
\usepackage[T1]{fontenc}
\usepackage{newtxtext,newtxmath}

\usepackage{graphicx}
\usepackage{booktabs}
\usepackage{subcaption}
\usepackage{multirow}
\usepackage[table]{xcolor}
\usepackage[export]{adjustbox}

\usepackage{tikz}
\usetikzlibrary{chains, scopes, positioning, backgrounds, shapes, fit, shadows, calc, 3d, arrows.meta, arrows}

\usepackage{cite}
\usepackage{hyperref}
\usepackage{orcidlink}
\usepackage{cleveref}
\usepackage{csquotes} % choose the right quotes for the language/font
\usepackage{microtype}

\allowdisplaybreaks[4]
\numberwithin{equation}{section}
\graphicspath{{figure/}}

\newtheoremstyle{myproblemstyle}
{\topsep}
{\topsep}
{\itshape}
{}
{\bfseries}
{}
{ }
{\thmname{#1}\thmnumber{ #2}\thmnote{ (#3)}}
\theoremstyle{myproblemstyle}
\newtheorem{problem}{Problem}[section]

\theoremstyle{plain}

\theoremstyle{definition}

\crefname{hypothesis}{Hypothesis}{Hypotheses}

\def\bs{\bigskip}

\begin{document}

\title{A Deep Learning-Enhanced Fourier Method for the Multi-Frequency Inverse Source Problem with Sparse Far-Field Data}

\author{
	Hao Chen\footnote{Research Center for Mathematics, Beijing Normal University, Zhuhai 519087, P. R. China. {\it 202331130059@mail.bnu.edu.cn}}, \
	Yan Chang\footnote{School of Mathematics, Harbin Institute of Technology, Harbin 150001, P. R. China. {\it changyan@hit.edu.cn} }, \
	Yukun Guo\footnote{School of Mathematics, Harbin Institute of Technology, Harbin 150001, P. R. China. {\it ykguo@hit.edu.cn}}, \
	Yuliang Wang\orcidlink{0000-0001-6264-650X}\footnote{Corresponding author. International Frontier Interdisciplinary Research Institute, Wenzhou-Kean University, Wenzhou 325060, P. R. China. {\it yulwang@wku.edu.cn} }
}
\date{}
\maketitle

\begin{abstract}
	This paper introduces a hybrid computational framework for the multi-frequency inverse source problem governed by the Helmholtz equation. By integrating a classical Fourier method with a deep convolutional neural network, we address the challenges inherent in sparse and noisy far-field data. The Fourier method provides a physics-informed, low-frequency approximation of the source, which serves as the input to a U-Net. The network is trained to map this coarse approximation to a high-fidelity source reconstruction, effectively suppressing truncation artifacts and recovering fine-scale geometric details. To enhance computational efficiency and robustness, we propose a high-to-low noise transfer learning strategy: a model pre-trained on high-noise regimes captures global topological features, offering a robust initialization for fine-tuning on lower-noise data. Numerical experiments demonstrate that the framework achieves accurate reconstructions with noise levels up to 100\%, significantly outperforms traditional spectral methods under sparse measurement constraints, and generalizes well to unseen source geometries.
\end{abstract}


\section{Introduction}
The inverse source scattering problem seeks to recover an unknown source from far-field measurements, and it underpins applications such as medical imaging \cite{Felicio2019, bao2002, Arridge1999}, pollution source identification \cite{BARATI2021, MALYSHEV1989}, and antenna synthesis \cite{Devaney2007, ramm1999multidimensional}. Despite its practical relevance, the problem is intrinsically ill-posed: measurements are taken only in the far field, while the source is localized, so small perturbations can lead to large reconstruction errors. In practical settings, reconstructions must also contend with sparse observation angles and measurement noise, which further degrade stability and resolution. These challenges motivate the development of efficient algorithms that can extract reliable source information from limited and imperfect data.

Numerical algorithms for this problem generally fall into two categories: iterative and non-iterative. Iterative approaches, such as Newton-type \cite{alves2009iterative, martins2012iterative} and boundary integral methods \cite{kress2015nonlinear}, generally offer high accuracy but suffer from high computational costs and sensitivity to initialization. Recursive algorithms using multi-frequency data \cite{bao2015inverse, bao2011inverse, bao2015recursive} help mitigate local minima. Non-iterative methods, including filtered backprojection \cite{griesmaier2012inverse, griesmaier2013inverse} and eigenfunction expansions \cite{eller2009acoustic}, offer efficiency but often rely on dense data. Recently, Direct Sampling Methods (DSM) \cite{LIU2025114251, CiCP-25-5} and Fourier series expansion methods \cite{zhang2015fourier, wang2017fourier} have been developed to reconstruct targets via indicator functions or spectral approximations, though their performance degrades with sparse data.

Deep learning has emerged as a powerful alternative, addressing the efficiency and robustness limitations of classical solvers. Deep learning-based approaches can be broadly categorized into three groups: supervised end-to-end learning, physics-informed learning, and hybrid methods. Supervised end-to-end methods, such as those employing U-Nets \cite{wei2018deep}, SwitchNet \cite{khoo2019switchnet}, or operator-learning frameworks like Fourier-DeepONet \cite{zhu2023fourier} and Invertible Fourier Neural Operators (iFNO) \cite{long2024invertible}, approximate the inverse operator directly but often lack physical interpretability. Physics-informed learning embeds PDE constraints into the training process, as seen in PINNs for parameter recovery \cite{Chen20, rasht2022physics, chen2024surface}. While rigorous, these methods can be computationally intensive to train.

Hybrid algorithms, combining neural networks with classical inversion techniques, have gained traction for effectively leveraging domain knowledge. For the inverse medium problem, methods integrating the Direct Sampling Method (DSM) with deep learning \cite{ning2023direct}, or employing two-step enhanced strategies \cite{yao2019two}, have shown promise in reconstructing scatterers. Similarly, Xu et al. \cite{xu2020deep} developed hybrid schemes combining contrast source inversion with CNNs for phaseless data, while Zhou et al. \cite{zhou2020improved} utilized a modified contrast scheme to enhance 2D and 3D reconstructions. Regarding the inverse source problem, hybrid strategies have also been explored. Li et al. \cite{li2023data} proposed a data-assisted hybrid approach for reconstructing the mean and variance of random sources. Du et al. \cite{du2023divide} developed a hybrid framework integrating deep neural networks with Bayesian inversion for point source identification.

Contributing to this line of hybrid research, we propose a deep-learning-enhanced Fourier method for the deterministic inverse source problem. Distinct from the aforementioned strategies targeting statistical or discrete sources, our approach addresses continuous source functions, aiming to recover fine-scale geometric details from sparse far-field data. Unlike operator-learning frameworks that map function spaces directly, our method decouples physical modeling from data-driven refinement. We use the classical Fourier method to generate a coarse, physics-informed initialization from sparse multi-frequency data. This approximation acts as a \enquote{warm start} for a U-Net, which is trained to recover diverse source geometries and suppress truncation artifacts. This formulation treats the learning process as one of artifact correction and super-resolution, significantly enhancing performance under sparse data constraints. We further introduce a high-to-low noise transfer learning strategy, where pre-training on high-noise regimes enables the model to capture robust structural priors, accelerating convergence on lower-noise data.

The remainder of this paper is organized as follows. Section 2 reviews the classical Fourier method. Section 3 details the proposed deep neural network model and transfer learning strategy. Section 4 presents numerical experiments validating the method's effectiveness. Finally, Section 5 concludes the paper.

\section{Problem setup}\label{section:fourier}
Let $S\in L^2(\mathbb{R}^2)$ be a source function such that
\(
	\text{supp }S\subset\subset V_0,
\)
where $V_0$ denotes a rectangular domain centered at the origin. The propagation of the acoustic wave $u$ generated by $S$ is modeled by the Helmholtz equation
\begin{align}\label{eq:Helmholtz}
	\Delta u+k^2u=S,\quad \text{in}\,\mathbb{R}^2,
\end{align}
where $k>0$ is the wavenumber. We assume $S$ to be independent of $k.$
Assuming the wave field $u$ satisfies the Sommerfeld radiation condition
\begin{align}\label{eq:Sommerfeld}
	\lim\limits_{r=|x|\to\infty}\sqrt{r}\left(\frac{\partial u}{\partial r}-\mathrm{i}ku\right)=0,
\end{align}
the solution to \eqref{eq:Helmholtz}--\eqref{eq:Sommerfeld} is given by
\[
	u(x;k)=-\int_{V_0}\Phi(x,y;k)S(y)\mathrm{d}y,
\]
where
\[
	\Phi(x,y;k)=\frac{\mathrm{i}}{4}H_0^{(1)}(k|x-y|)
\]
is the fundamental solution to the Helmholtz equation and $H_0^{(1)}$ is the Hankel function of the first kind of zero order.
From the asymptotic behavior of $H_0^{(1)}$ \cite{colton2019inverse}, $u(x;k)$ admits the asymptotic expansion
\begin{align}\label{eq:asymptotic}
	u(x;k)=\frac{\mathrm{e}^{\mathrm{i}k|x|}}{\sqrt{|x|}}\left\{
	u^\infty(\hat{x};k)+\mathcal{O}\left(\frac{1}{|x|}\right)
	\right\},\quad|x|\to\infty,
\end{align}
uniformly in all directions $\hat{x}=x/|x|.$
In \eqref{eq:asymptotic}, $u^\infty$ is the far-field pattern defined by
\begin{equation*}
	u^\infty(\hat{x};k)=-\gamma\int_{V_0}S(y)\mathrm{e}^{-\mathrm{i}k\hat{x}\cdot y}\mathrm{d}y,\
	\hat{x}\in\mathbb{S},
\end{equation*}
where $\mathbb{S}=\left\{x\in\mathbb{R}^2:|x|=1\right\}$, and
\(
	\gamma=\frac{\mathrm{e}^{\mathrm{i}\pi/4}}{\sqrt{8\pi k}}.
\)

In this paper, we are concerned with the following inverse problem:
\begin{problem}[Multi-frequency ISP with far-field pattern]
Given a finite number of frequencies $\{k\}$, determine the source function $S(x)$ from the far-field data $\{u^\infty(\hat{x}_k;k)\},$ where the observation direction $\hat{x}_k$ depends on $k$.
\end{problem}

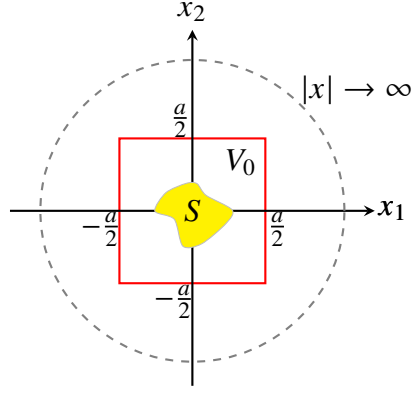
\begin{figure}
	\begin{center}
		\begin{tikzpicture}[scale=0.8]
			\draw [red, thick](-1.2,-1.2) rectangle (1.2,1.2);
			\draw [thick,-stealth] (-3,0) -- (3,0);
			\draw [thick,-stealth] (0,-2.9) -- (0,3) node at (0,3.3) {$x_2$} node at (2.7,2) {$|x|\to\infty$};
			\draw  node at (3.3,0) {$x_1$}  node at (1.4,-0.3) {$\frac{a}{2}$}  node at (-1.5,-0.3) {$-\frac{a}{2}$};
			\draw  node at (3.3,0) {$x_1$}  node at (-.2,1.5) {$\frac{a}{2}$}  node at (-.3,-1.5) {$-\frac{a}{2}$};
			\draw [thick, dashed,gray] (0,0) circle (2.5);\draw node at (0.8, 0.8) {$V_0$};
			\pgfmathsetseed{1235}
			\draw (2,2) plot[gray, smooth cycle, samples=10, domain={1:12}] (2+\x*360/12+3*rnd:.2cm+.5cm*rnd) node at (0, 0) {$S$}[lightgray, fill=yellow];
		\end{tikzpicture}
	\end{center}
	\caption{Illustration of the inverse source scattering problem.}\label{fig:setup}
\end{figure}

\Cref{fig:setup} schematically illustrates the geometric setting, where the yellow domain represents the support of $S$ and the red square designates the integration domain
\(
	V_0=\left(-\frac{a}{2},\frac{a}{2}\right)^2
\)
, chosen such that $S\subset\subset V_0$. Define the Fourier basis functions $\phi_{\bm l}(x)=\mathrm{exp}\left(\mathrm{i}\frac{2\pi}{a}{\bm l}\cdot x\right)$ for ${\bm l}\in\mathbb{Z}^2, x \in V_0$. The source function $S(x)$ is approximated by the truncated Fourier expansion
\begin{align}\label{eq:SN}
	S_N(x)=\sum_{|\bm l|_\infty\le N}\hat{s}_{\bm l}\phi_{\bm l}(x),
\end{align}
where $\hat{s}_{\bm l}$ are the Fourier coefficients and $N$ is the truncation frequency.

For ${\bm l} \in \mathbb{Z}^2$, define the wavenumbers $k_{\bm l}$ and the observation directions $\hat{x}_{\bm l}$ by
\begin{align}\label{eq:fre-dir}
	k_{\bm l}=\left\{
	\begin{aligned}
		 & \frac{2\pi}{a}|\bm l|, &  & {\bm l}\in\mathbb{Z}^2\backslash\{\bm 0\}, \\
		 & \frac{2\pi}{a}\lambda, &  & {\bm l}={\bm0},
	\end{aligned}
	\right.\quad
	\hat{x}_{\bm l}=\left\{
	\begin{aligned}
		 & \frac{\bm l}{|\bm l|}, &  & {\bm l}\in\mathbb{Z}^2\backslash\{\bm 0\}, \\
		 & (1,0),                 &  & {\bm l}={\bm0},
	\end{aligned}
	\right.
\end{align}
where $\lambda>0$ is a constant such that $\frac{2\pi}{a}\lambda<\frac{1}{2}$.
The computational formula for $\hat{s}_{\bm l}$ from the corresponding far-field data $u^\infty(\hat{x}_{\bm l};k_{\bm l})$ follows from \cite{wang2017fourier}:
\begin{align}
	\label{eq:sl} & \hat{s}_{\bm l}=-\frac{u^\infty(\hat{x}_{\bm l};k_{\bm l})}{a^2\gamma},\quad{\bm l}\in\mathbb{Z}^2\backslash\{\bm 0\}, \\
	\label{eq:s0} & \hat{s}_{\bm 0}=-\frac{\lambda\pi}{a^2\sin\lambda\pi}\left(
	\frac{u^\infty(\hat{x}_{\bm0};k_{\bm0})}{\gamma}+\sum_{1\le|\bm l|_\infty\le N}\hat{s}_{\bm l}
	\int_{V_0}\phi_{\bm l}(y)\overline{\phi_{{\bm l}_0}(y)}\mathrm{d}y
	\right), \quad {\bm l}_0=(\lambda,0).
\end{align}

The explicit formulas \eqref{eq:sl}--\eqref{eq:s0} provide a direct method to reconstruct the Fourier coefficients (see \cite{wang2017fourier} for derivation). In practice, if the far-field data is contaminated by noise, denoted as $u^\infty_\delta$, we substitute $u^\infty_\delta$ into \eqref{eq:sl}--\eqref{eq:s0} to obtain the perturbed Fourier coefficients $\hat{s}_{\bm l}^\delta$. The noisy Fourier approximation is then defined as $S_N^\delta(x) = \sum_{|\bm l|_\infty\le N}\hat{s}_{\bm l}^\delta\phi_{\bm l}(x)$.
	
\section{The Deep-learning-enhanced Fourier method}
Although the classical Fourier method allows for explicit reconstruction of source functions, its accuracy is strictly limited by the truncation frequency $N$. As implied by \eqref{eq:SN} and \eqref{eq:fre-dir}, capturing fine-scale details requires a large $N$, which in turn demands dense multi-frequency and multi-angle measurements. In many practical scenarios, such data is either unavailable or contaminated by noise. When restricted to sparse data (small $N$), the truncated Fourier series $S_N$ yields a smooth but low-resolution approximation, effectively filtering out high-frequency details. While increasing $N$ allows for the recovery of finer structures, it typically introduces ringing artifacts (Gibbs phenomenon) around discontinuities and amplifies susceptibility to noise.

To overcome the resolution limits of the classical method, we formulate the problem as learning the inverse of the truncation operator $\mathcal{A}$, which maps the true source $S$ to its band-limited approximation $S_N = \mathcal{A}(S)$. We approximate the inverse mapping $\mathcal{A}^{-1}$ using a deep neural network $\mathcal{G}_\Theta$, parameterized by $\Theta$. This is cast as an image-to-image translation task where the network recovers the source from its Fourier approximation:
\[
	S \approx \mathcal{G}_\Theta(S_N).
\]
To learn this mapping, we minimize a loss metric $\mathcal{L}$ quantifying the discrepancy between the network output and the ground truth. The parameters $\Theta$ are updated via stochastic gradient descent:
\begin{equation*}
	\Theta \leftarrow \Theta - \frac{1}{|\mathcal{I}|}\sum_{k\in\mathcal{I}}\tau\nabla_\Theta\mathcal{L}\left(\mathcal{G}_\Theta(S_{N, k}), S_k\right),
\end{equation*}
where $\tau$ is the learning rate and $\mathcal{I}$ denotes the index set for a mini-batch of training data. Through exposure to diverse examples, the network learns to suppress noise, eliminate artifacts, and restore high-frequency details lost in the band-limited Fourier reconstruction. In contrast to end-to-end approaches that directly map far-field data to the source function, this hybrid methodology combines the physical interpretability of the Fourier method with the pattern-recognition capabilities of deep learning. This yields more accurate and stable reconstructions under challenging conditions characterized by sparse or noisy measurements. Importantly, the neural network does not learn the forward or inverse scattering operator from scratch; instead, it numerically approximates the inverse of a fixed, low-frequency Fourier truncation operator, allowing the learning task to focus on artifact suppression and resolution enhancement.

\subsection{The training set}
We use deep learning to enhance the initial Fourier reconstruction. The success of this approach hinges on the construction of an effective dataset, where the choice of the truncation frequency $N$ plays a pivotal role. This choice represents a critical stability--resolution trade-off:

\begin{itemize}
	\item Data Acquisition Feasibility: As established in \eqref{eq:fre-dir}, the method requires far-field data at specific wavenumbers $k_{\bm l}$ and observation directions $\hat{x}_{\bm l}$. The number of required measurements scales with $N$, making large $N$ prohibitively expensive in practice.
	\item High-Frequency Instability: High-frequency Fourier coefficients are notoriously susceptible to measurement noise. A smaller $N$ acts as a regularization parameter, ensuring a stable, albeit coarse, initial reconstruction.
\end{itemize}

\subsection{The U-Net architecture}
\begin{figure}[htbp]
	\centering
	\includegraphics[width=0.85\linewidth]{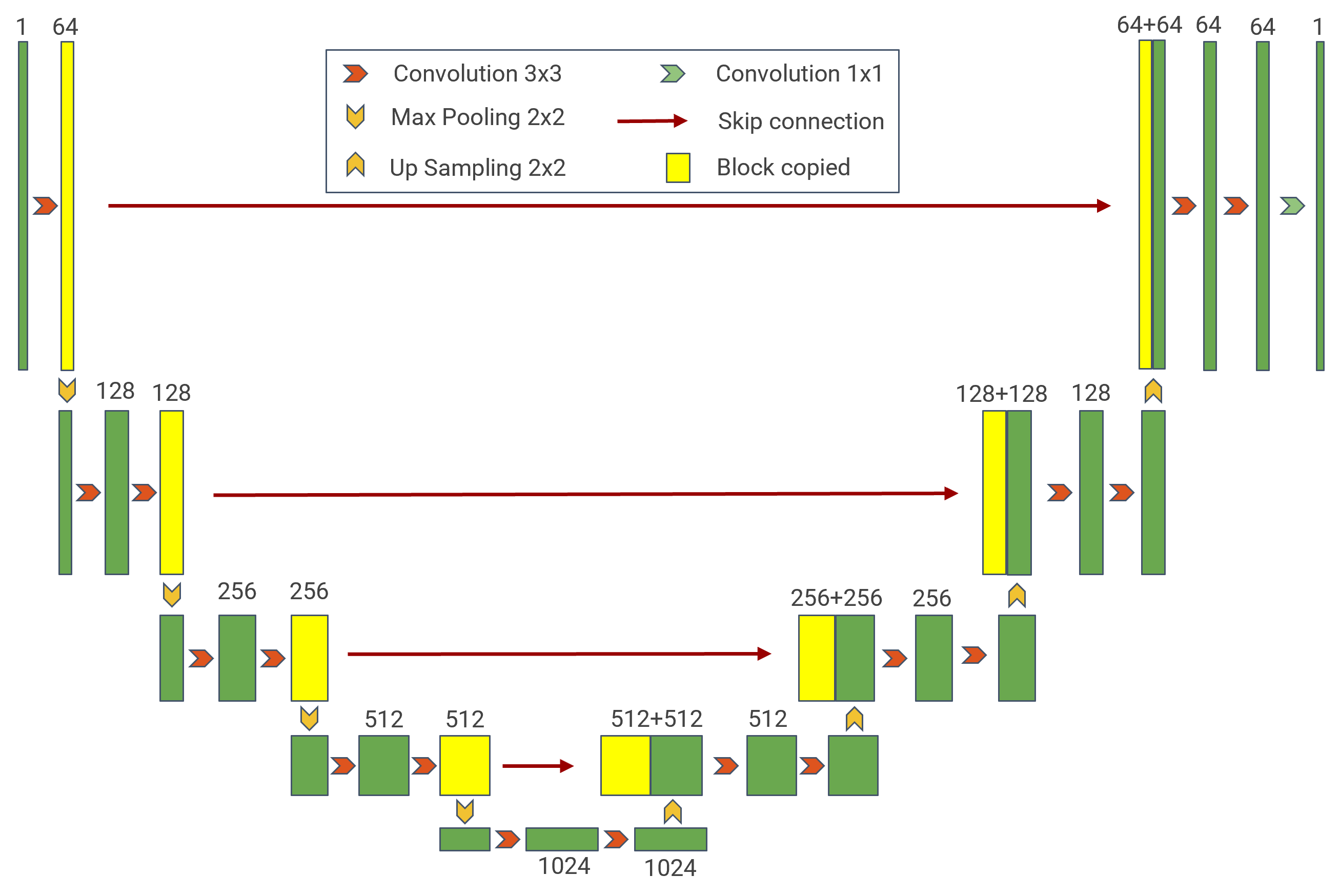}
	\caption{Schematic of the U-Net architecture used in this study.}
	\label{fig:U-Net}
\end{figure}
We employ the U-Net architecture \cite{ronneberger2015u} as the backbone of our deep learning model. As a fully convolutional network, U-Net has demonstrated superior performance compared to conventional Convolutional Autoencoders, particularly when training data is limited. The architecture, illustrated in \Cref{fig:U-Net}, features a symmetric U-shaped design comprising two primary paths. The contracting path (left) consists of repeated blocks, each containing two sequential convolution operations (convolution, batch normalization, and ReLU activation) followed by a max-pooling layer. After each convolution step, the spatial dimensions of the feature map are halved while the number of channels is doubled. Conversely, the extensive path (right) mirrors the contracting path but employs $3\times3$ transposed convolutions for upsampling, restoring dimensions to match the corresponding levels of the contracting path. Crucially, skip connections link the two paths, transferring high-resolution feature maps to the expansive path to preserve fine-grained spatial information lost during downsampling. 

Let $\{(S_{N,i}^\delta, S_i)\}_{i=1}^M$ be the training dataset. Here, $S_{N,i}^\delta \in\mathbb{R}^{n\times n}$ is the matrix computed from measurement data via the Fourier method (with grid size $n=64$), $S_i\in\mathbb{R}^{n\times n}$ is the discretized ground truth, and $M$ is the total number of samples. To train the U-Net $\mathcal{G}_\Theta$, we minimize the Mean Squared Error (MSE) loss:
\begin{equation*}
	\mathcal{L} = \mathbb{E}(\|S-\mathcal{G}_\Theta(S_N^\delta)\|_{L^2}^2),
\end{equation*}
where $\mathbb{E}$ represents the expectation over the distribution induced by the training dataset. In practice, this expectation is approximated by averaging the loss over a mini-batch. While other metrics such as $L^1$ loss or determining regularizers (e.g., total variation) were considered, our experiments indicate that the $L^2$ loss yields superior reconstruction quality.

\subsection{Transfer learning strategy}
The proposed framework allows for training on specific noise levels for a fixed truncation frequency $N$. However, a standard U-Net model is typically specialized to the noise level encountered during training. Training a distinct model from scratch for every new noise condition is computationally inefficient and fails to exploit the shared mathematical structure of the problem. Transfer learning addresses this by reusing a model trained on one dataset to initialize training on a related task, thereby improving efficiency and performance.

We employ transfer learning primarily to accelerate the training process, proposing a ``high-to-low'' noise transfer strategy. First, we train a U-Net model on a high-noise dataset (e.g., $100\%$ noise), yielding a model $\mathcal{G}_\Theta^\delta$ that learns to identify source signals amidst significant interference. Physically, training on high-noise data forces the network to disregard localized fluctuations and focus on the underlying global structure and salient topological features of the source. This pre-trained model $\mathcal{G}_\Theta^\delta$ then serves as the initialization for a new task involving a lower noise level $\delta' < \delta$. The parameters $\Theta$ are subsequently fine-tuned on the new dataset to produce the final model $\mathcal{G}_\Theta^{\delta'}$. This strategy acts as a regularizer, leveraging structural priors learned from the more difficult task to guide optimization in the lower-noise task, helping to avoid poor local minima and accelerating convergence.

\bs

\section{Numerical experiments}

This section presents a series of numerical experiments to evaluate the effectiveness of the proposed method. To assess stability, random noise is added to the asymptotic data $u^{\infty}$. The noisy data are modeled as:
\begin{equation*}
	u^{\infty}_{\delta} := u^{\infty} + \delta |u^{\infty}| r_1 e^{i\pi r_2},
\end{equation*}
where $r_1$ and $r_2$ are independent random variables uniformly distributed in $[-1, 1]$, and $\delta > 0$ represents the noise level.

To obtain multi-frequency far-field data, we fix the truncation frequency $N$ and define the set of wavenumbers:
\begin{equation*}
	\mathbb{K}_N := \{2\pi|\bm{l}|: \bm{l}\in\mathbb{Z}^2, 1\leq|\bm{l}|_{\infty}\leq N\}\cup\{2\pi\lambda\}, \quad \lambda = 10^{-3}.
\end{equation*}
This results in a total of $(2N+1)^2$ scattering data points. The training pairs consist of images with pixel values given by $S$ and $S_N^\delta$, discretized on a $64 \times 64$ uniform grid over the sampling domain $V_0= (-0.5, 0.5)^2$.

Our U-Net implementation employs $3\times3$ convolutional kernels and batch normalization. Training and testing are performed in PyTorch using the Adam optimizer \cite{kingma2014adam} with a batch size of $32$. All computations are accelerated on an NVIDIA A100 80GB GPU. The reconstruction performance is quantified using the normalized mean squared error (NMSE), defined as:
\begin{equation*}
	\text{NMSE} = \frac{\|\mathcal{G}_\Theta(S_N^\delta) - S\|_2^2}{\|S\|_2^2}.
\end{equation*}

\subsection{Baseline Validation: Reconstruction of Disk Sources}
\label{example:disks}
In the first example, we simulate sources consisting of one to three disks. The number of disks, their centers (within $V_0$), and radii (drawn from $[0.1, 0.2]$) are sampled uniformly, identifying only those configurations entirely contained within $V_0$. The intensity of $S(x)$ on each disk is drawn uniformly from $[-1, 1]$, and overlapping regions are assigned the value of the most recently added disk. We set $N = 3$ and compute the initial source functions $S_N^\delta(x)$ via the Fourier method. The dataset comprises $2000$ samples, split into $1600$ for training and $400$ for testing. Training proceeds for $200$ epochs, with an initial learning rate of $0.001$ that decays by a factor of $0.9$ every $5$ epochs.
\begin{figure}[htbp]\small
	\begin{center}
		\begin{tabular}{ccccc}
			Ground truth & Fourier ($N=3$) & U-Net & Fourier ($N=10$) & Noise Level \\
			\includegraphics[valign=m,width=0.15\textwidth]{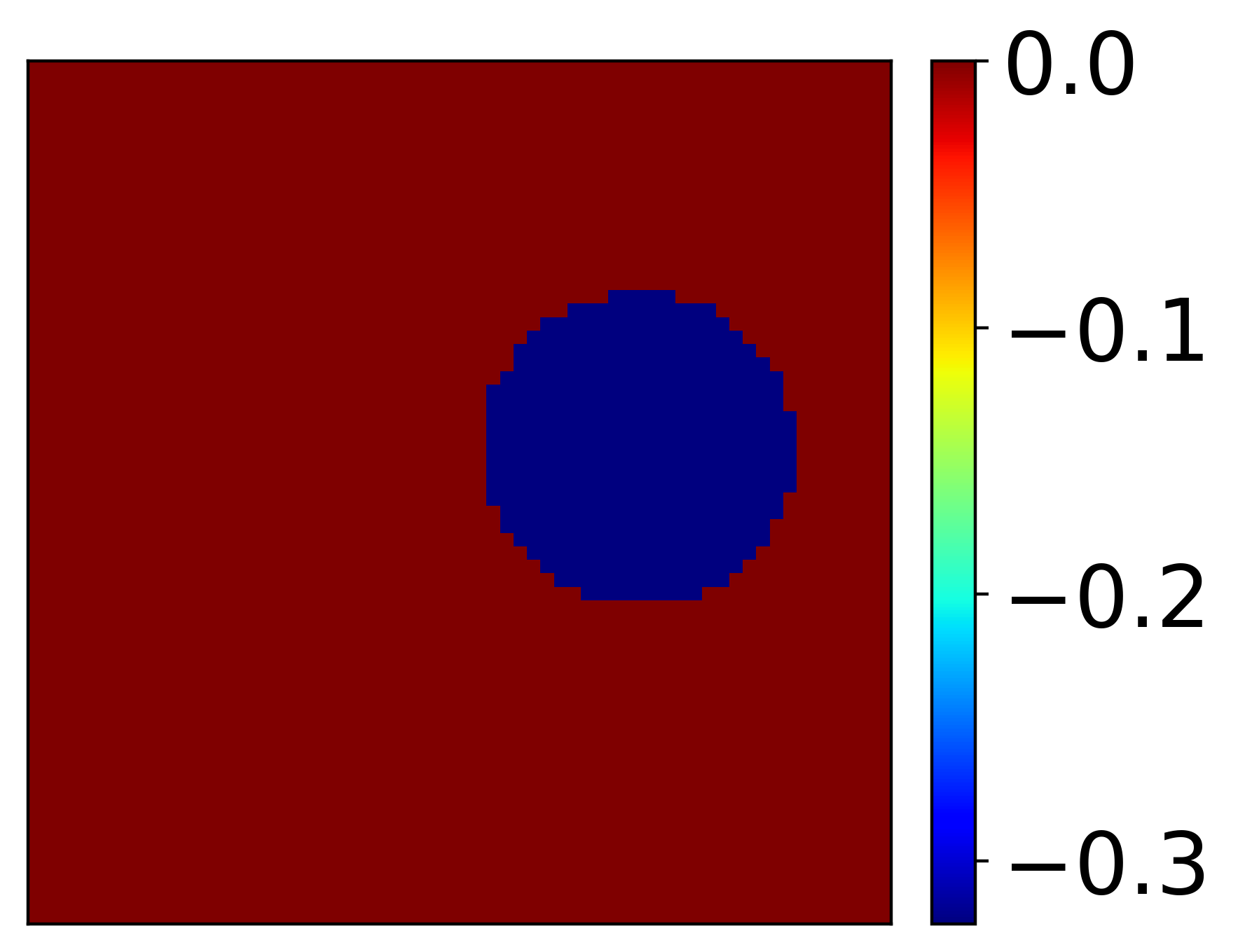} &
			\includegraphics[valign=m,width=0.15\textwidth]{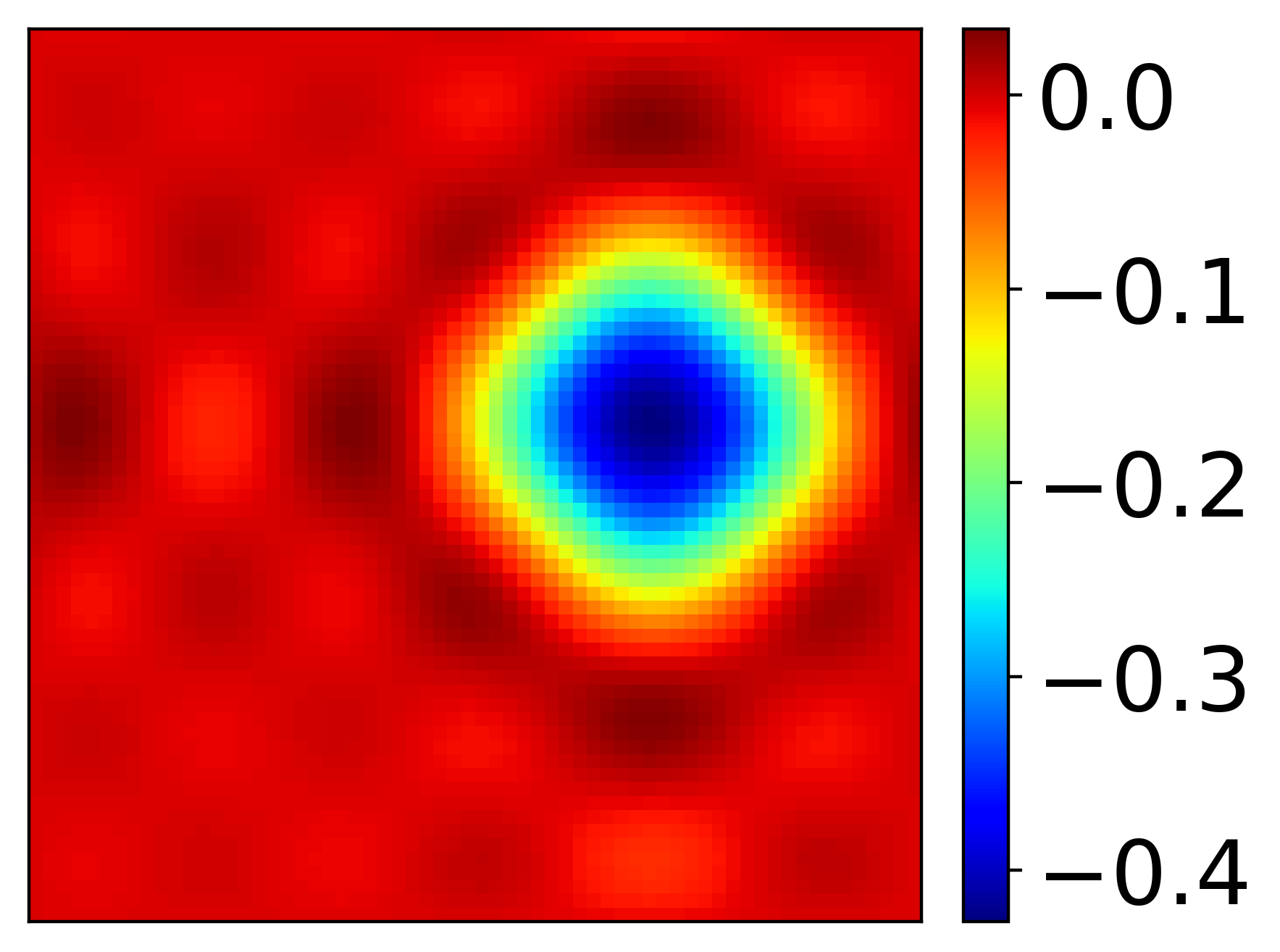} &
			\includegraphics[valign=m,width=0.15\textwidth]{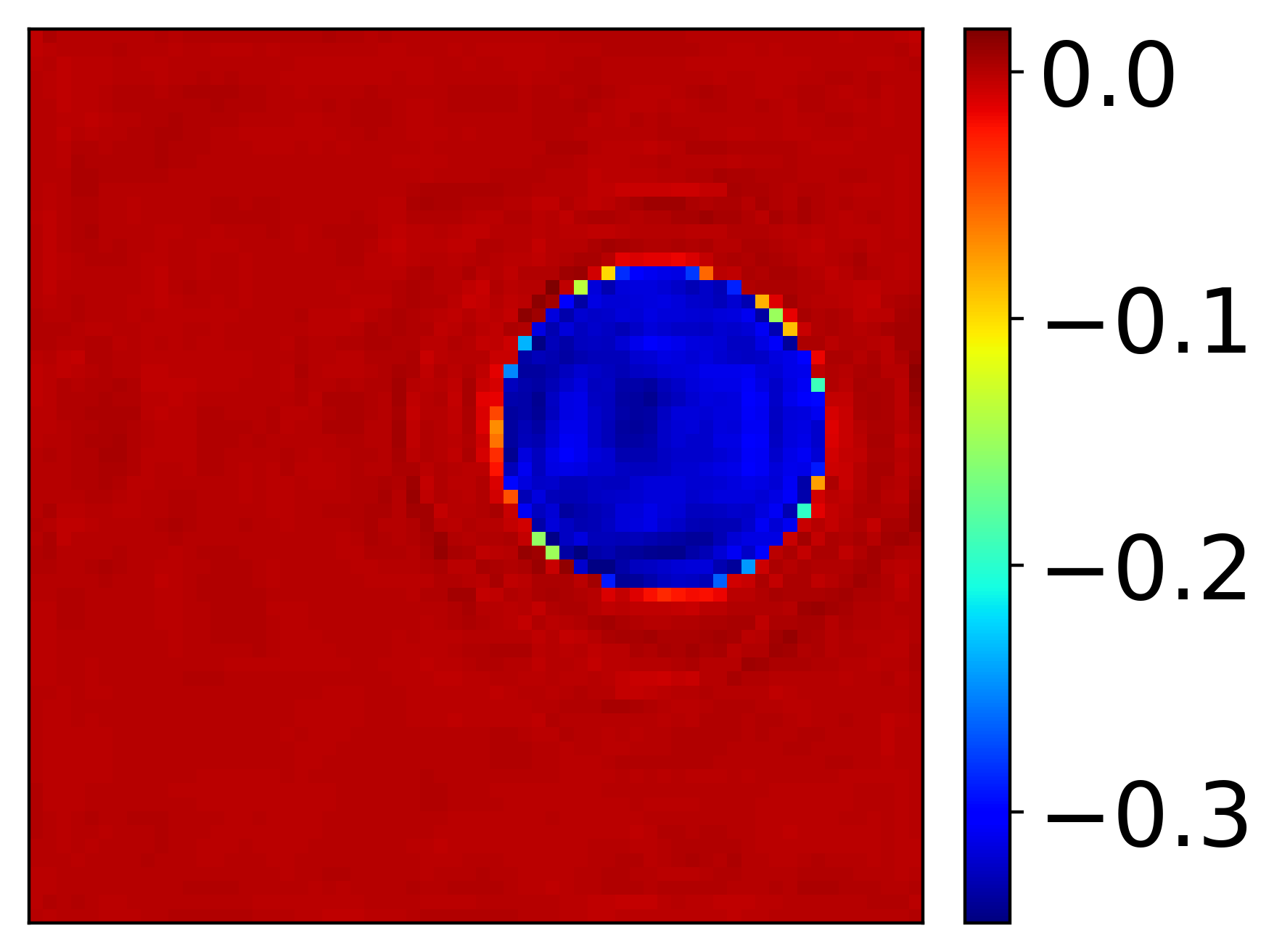} &
			\includegraphics[valign=m,width=0.15\textwidth]{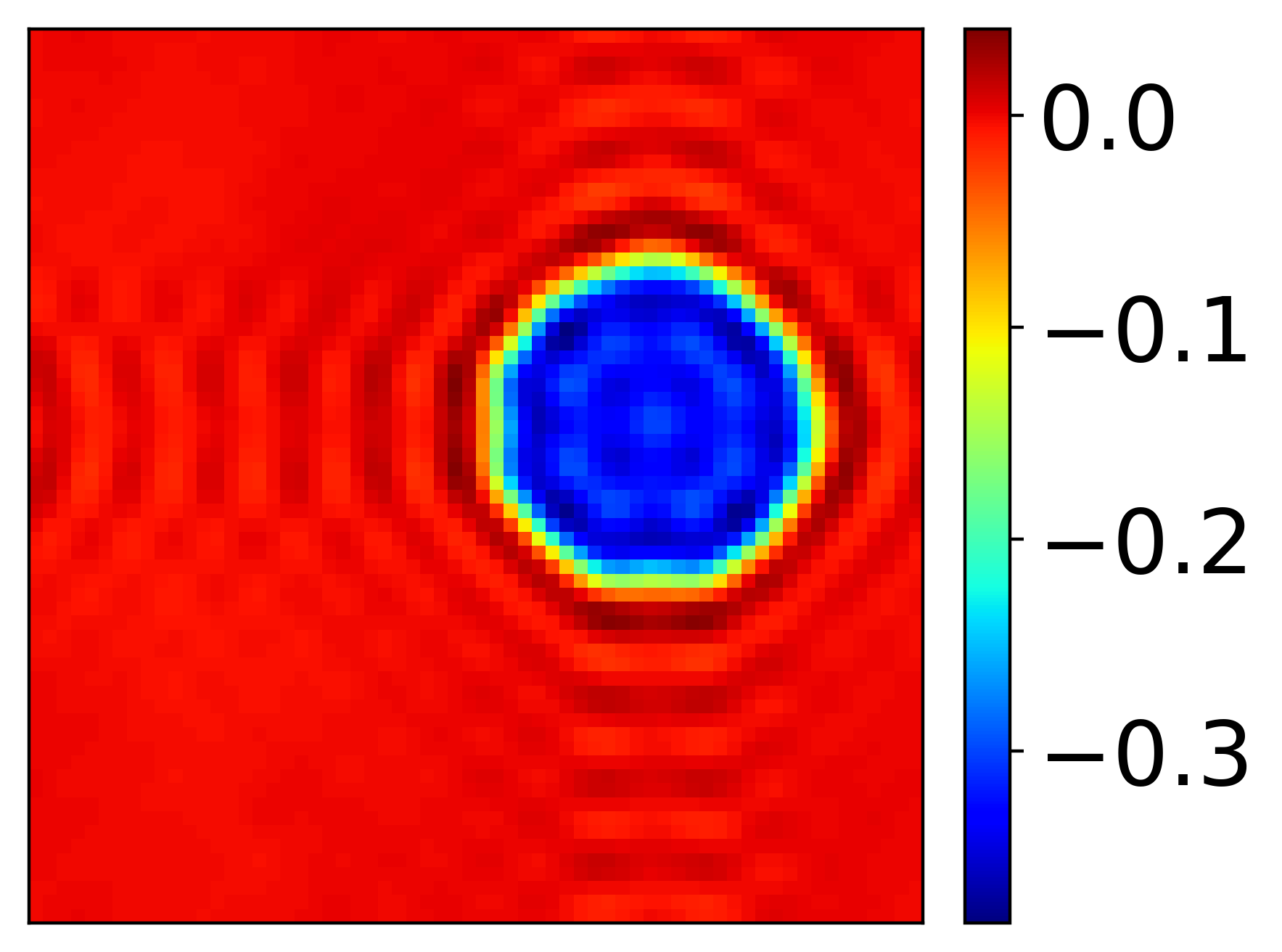} & 5\% \\
			& \includegraphics[valign=m,width=0.15\textwidth]{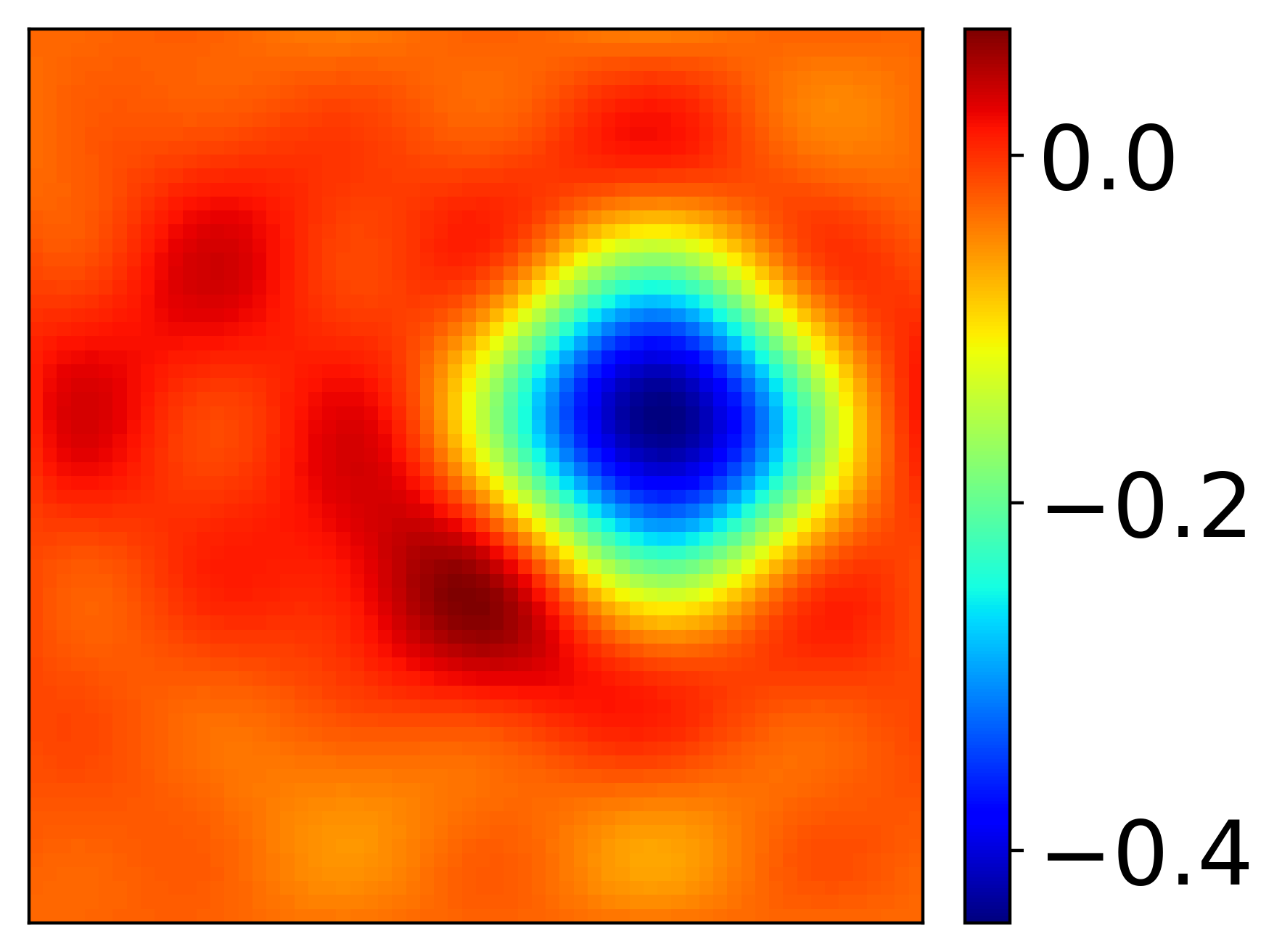} &
			\includegraphics[valign=m,width=0.15\textwidth]{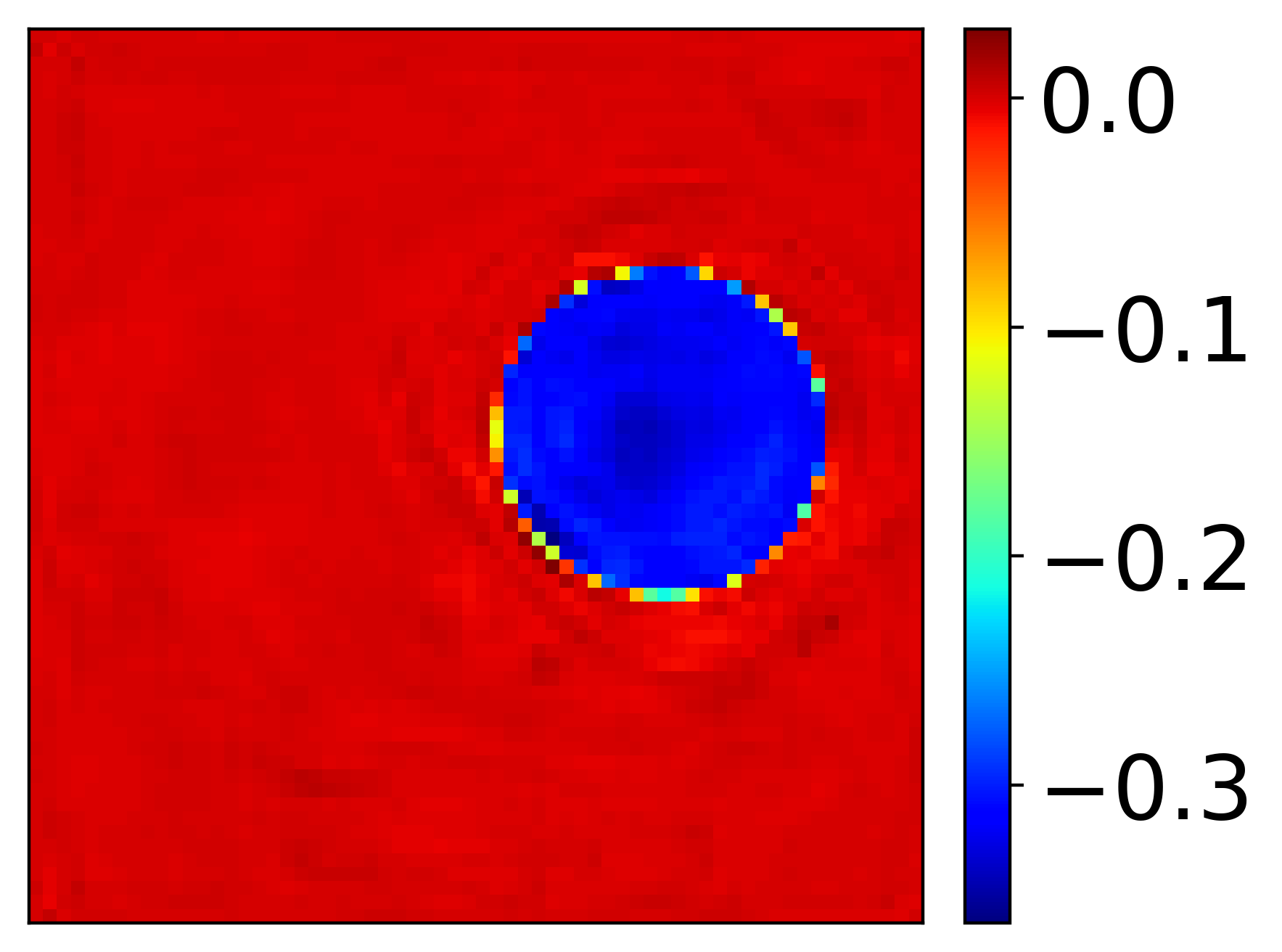} &
			\includegraphics[valign=m,width=0.15\textwidth]{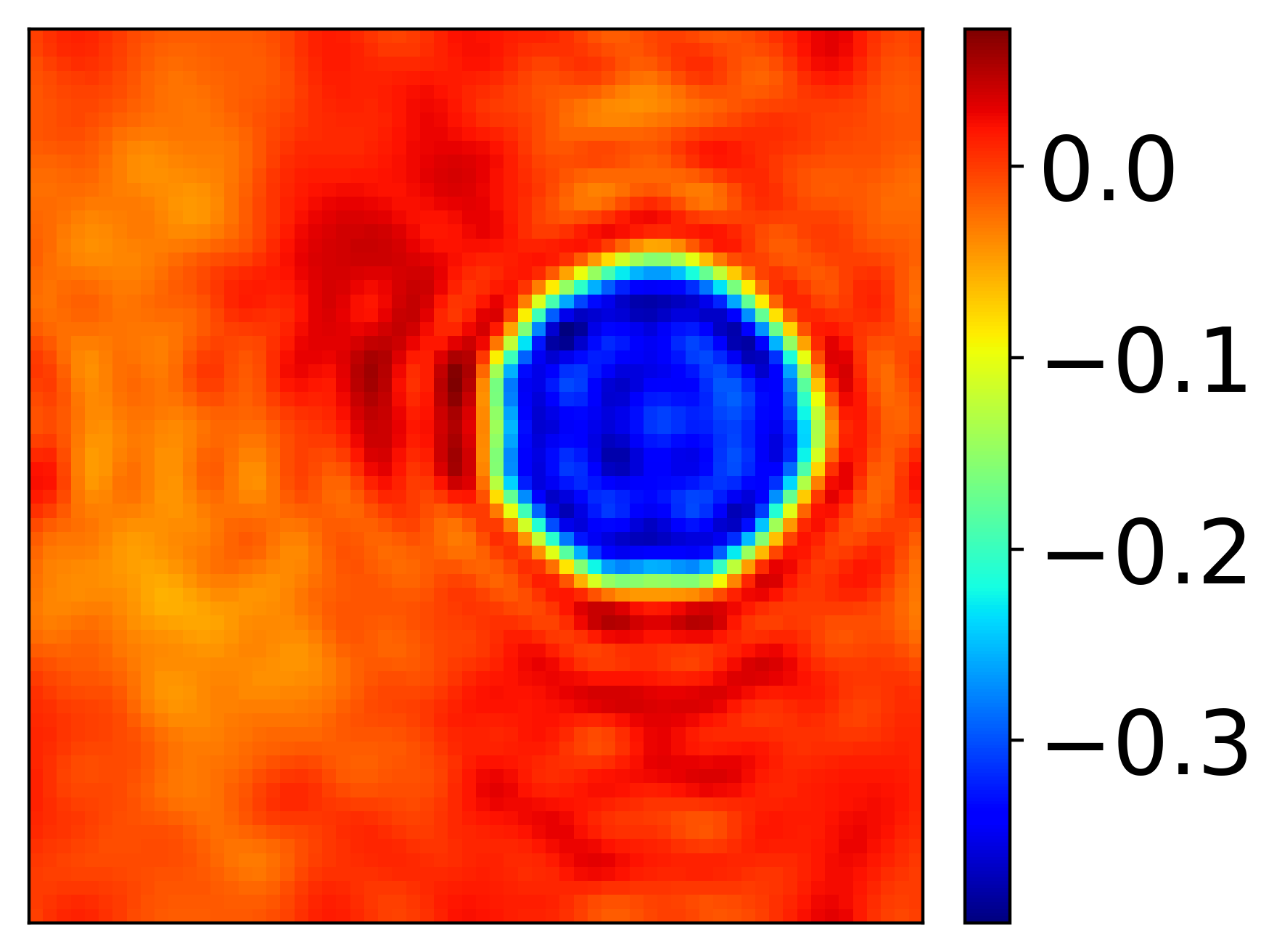} & 50\% \\
			& \includegraphics[valign=m,width=0.15\textwidth]{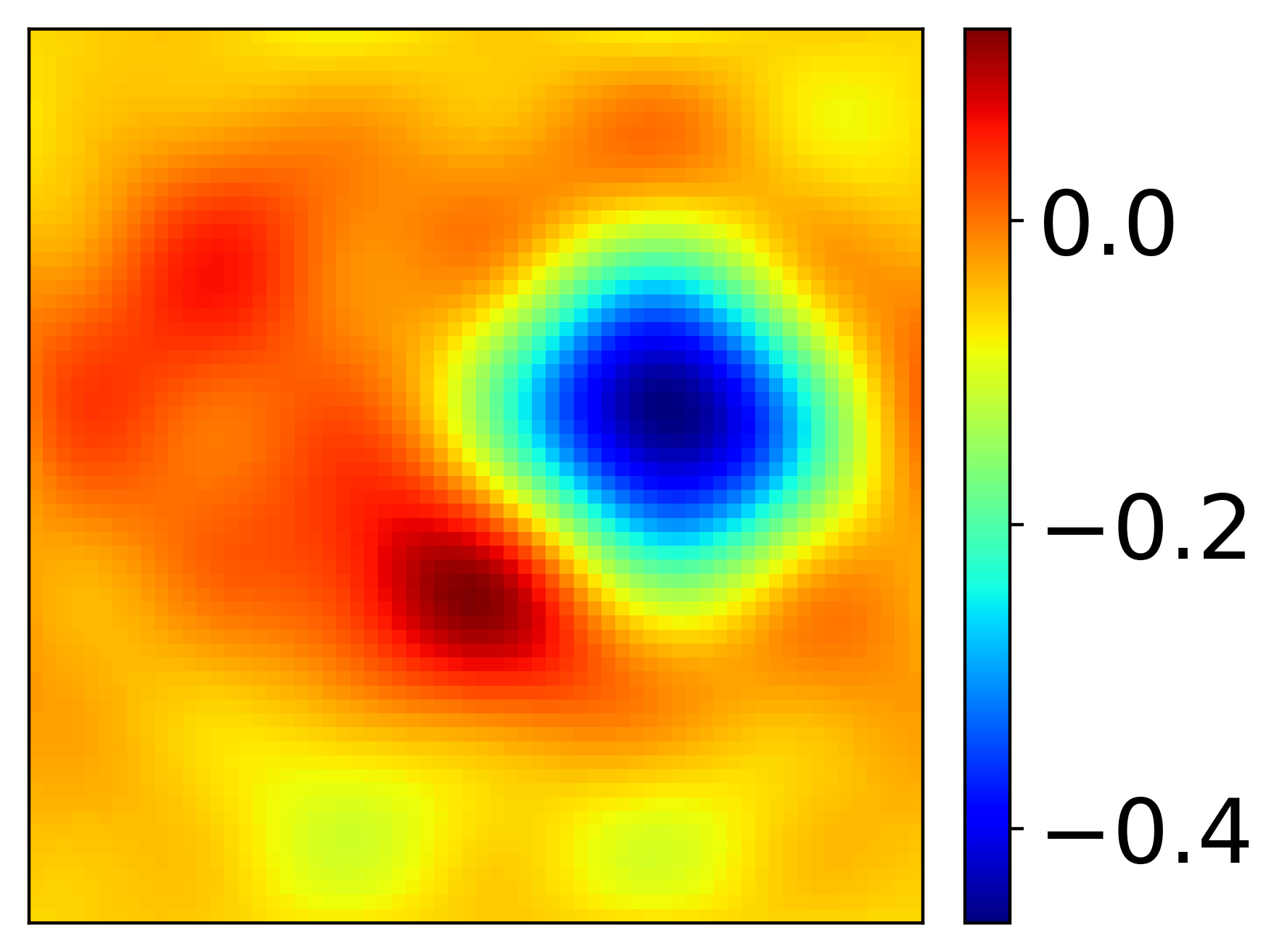} &
			\includegraphics[valign=m,width=0.15\textwidth]{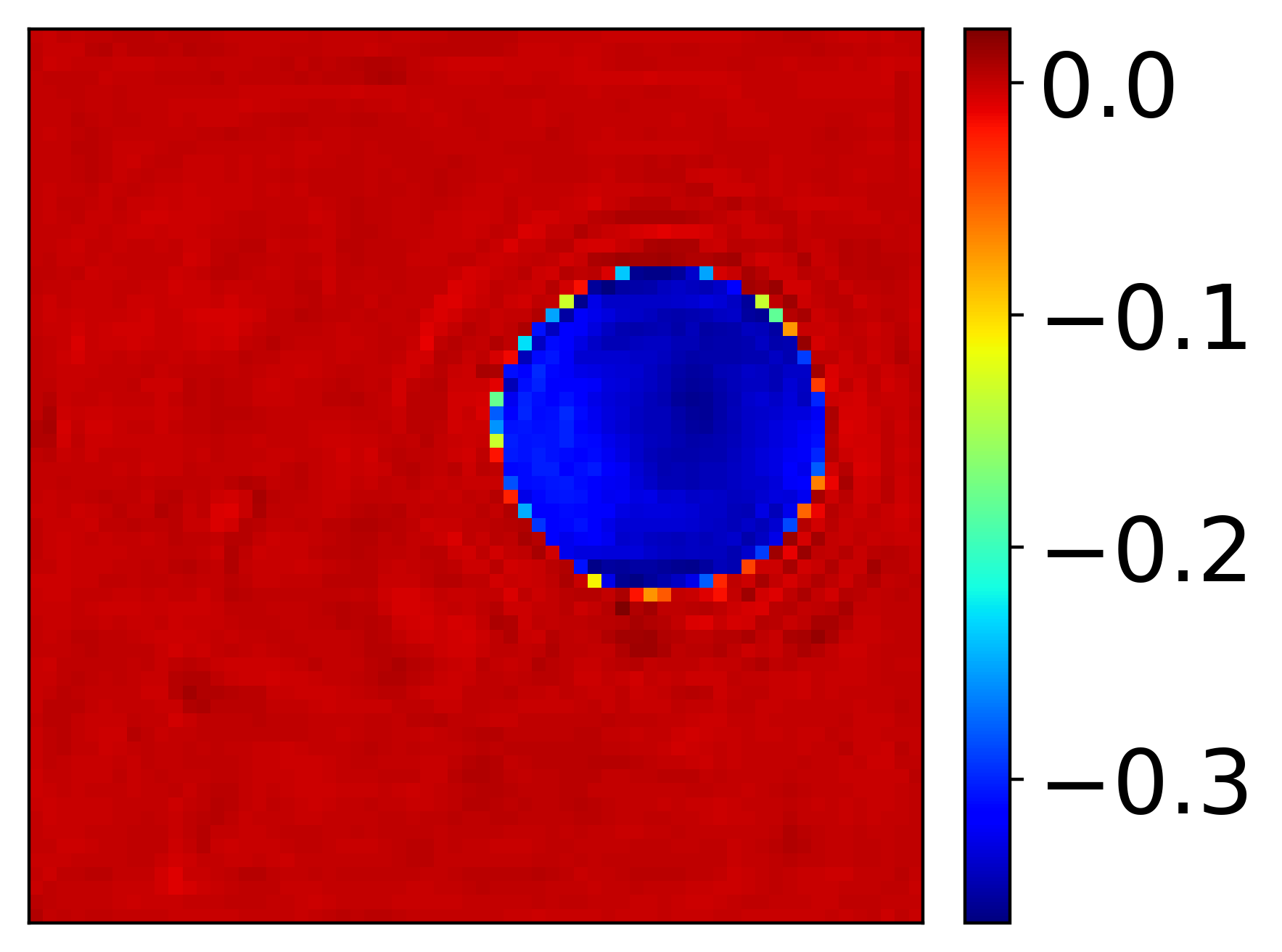} &
			\includegraphics[valign=m,width=0.15\textwidth]{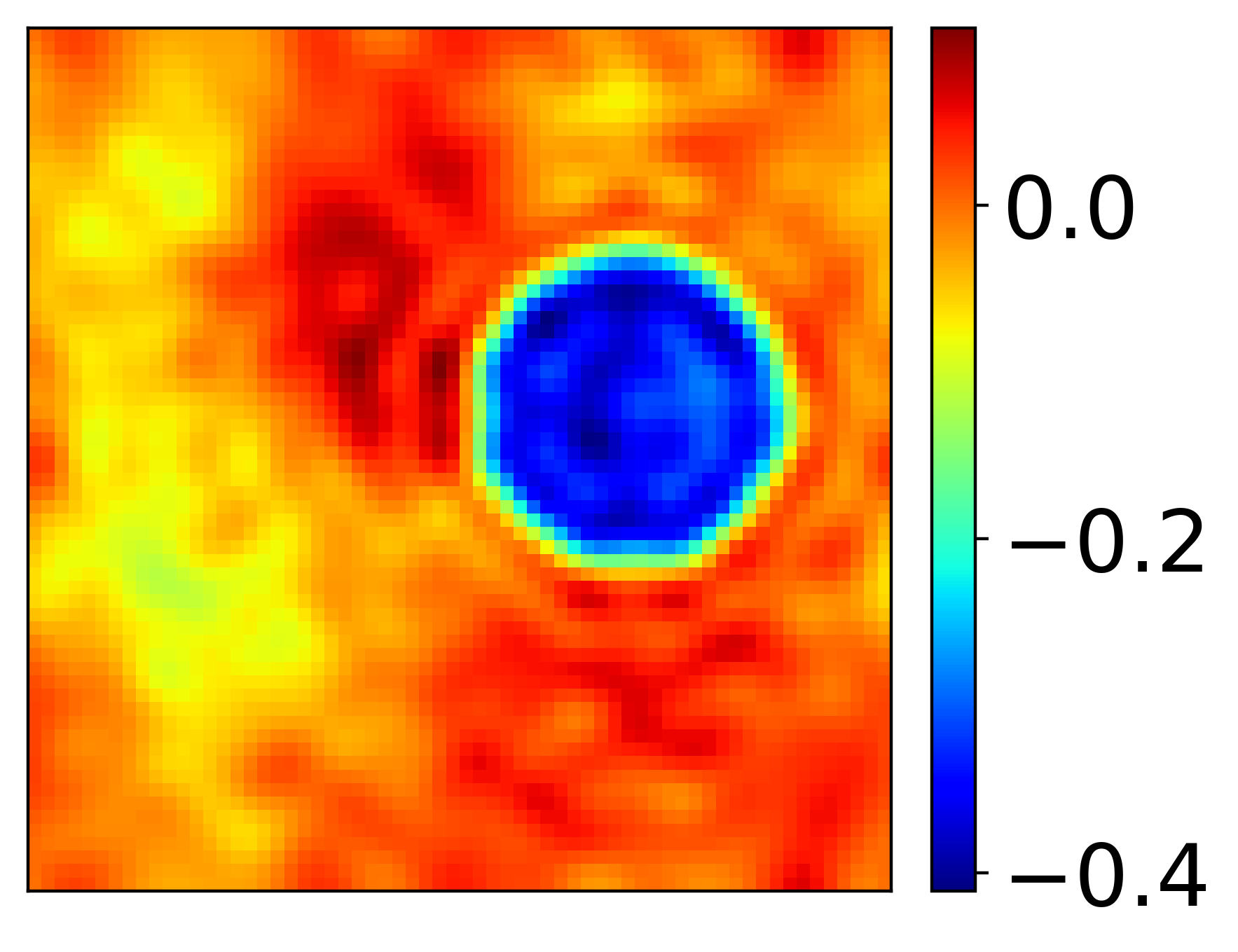} & 100\% \\
			\includegraphics[valign=m,width=0.15\textwidth]{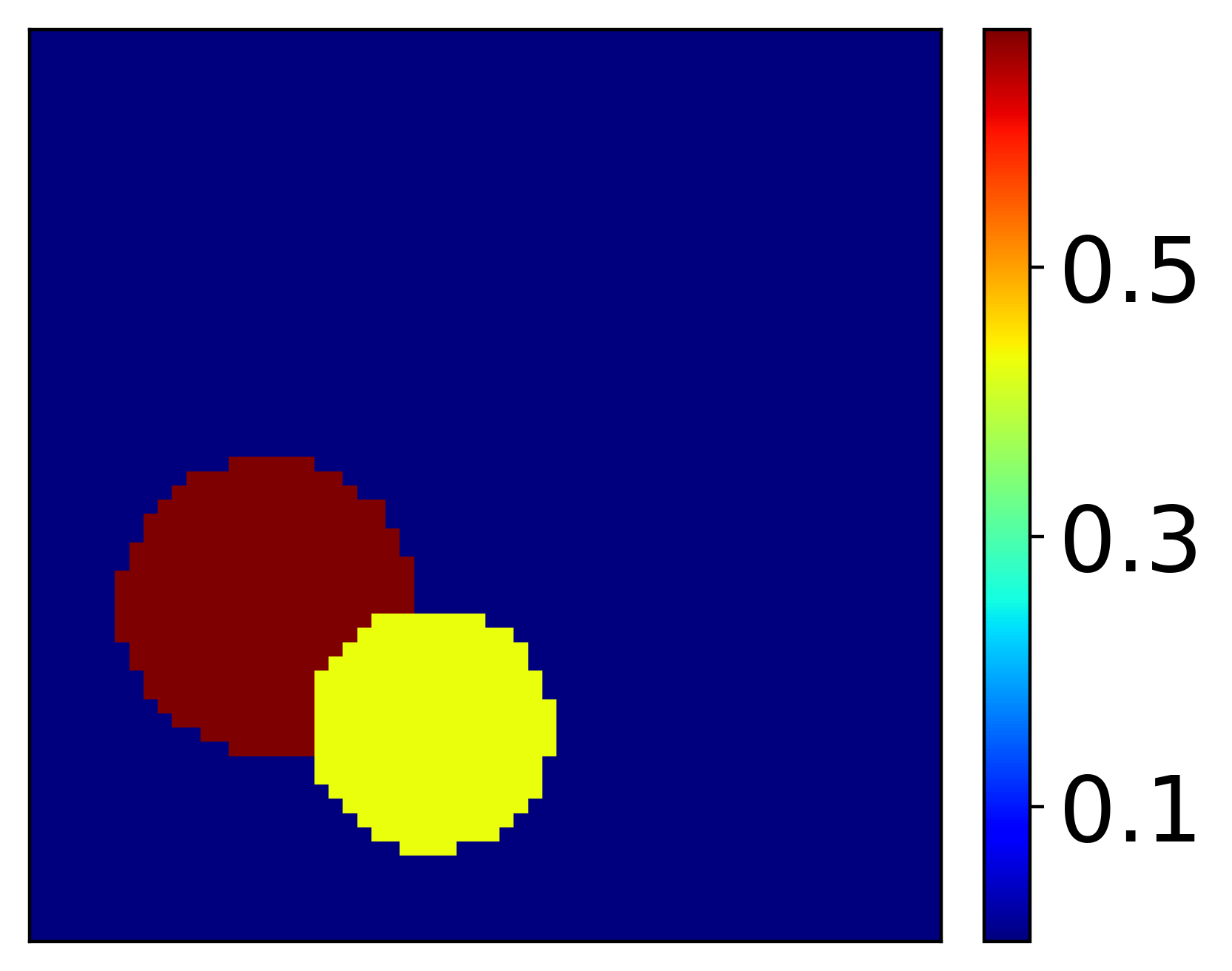} &
			\includegraphics[valign=m,width=0.15\textwidth]{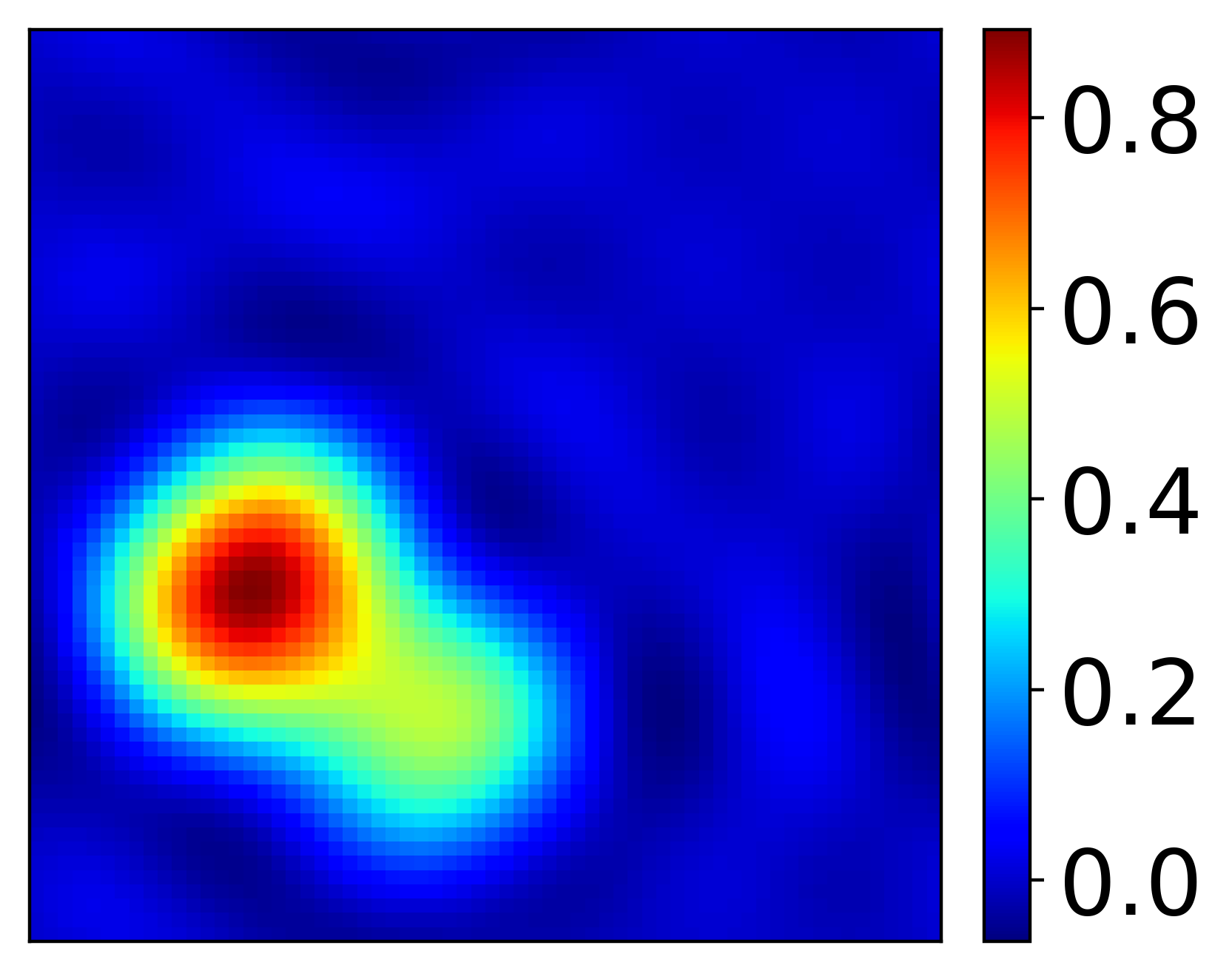} &
			\includegraphics[valign=m,width=0.15\textwidth]{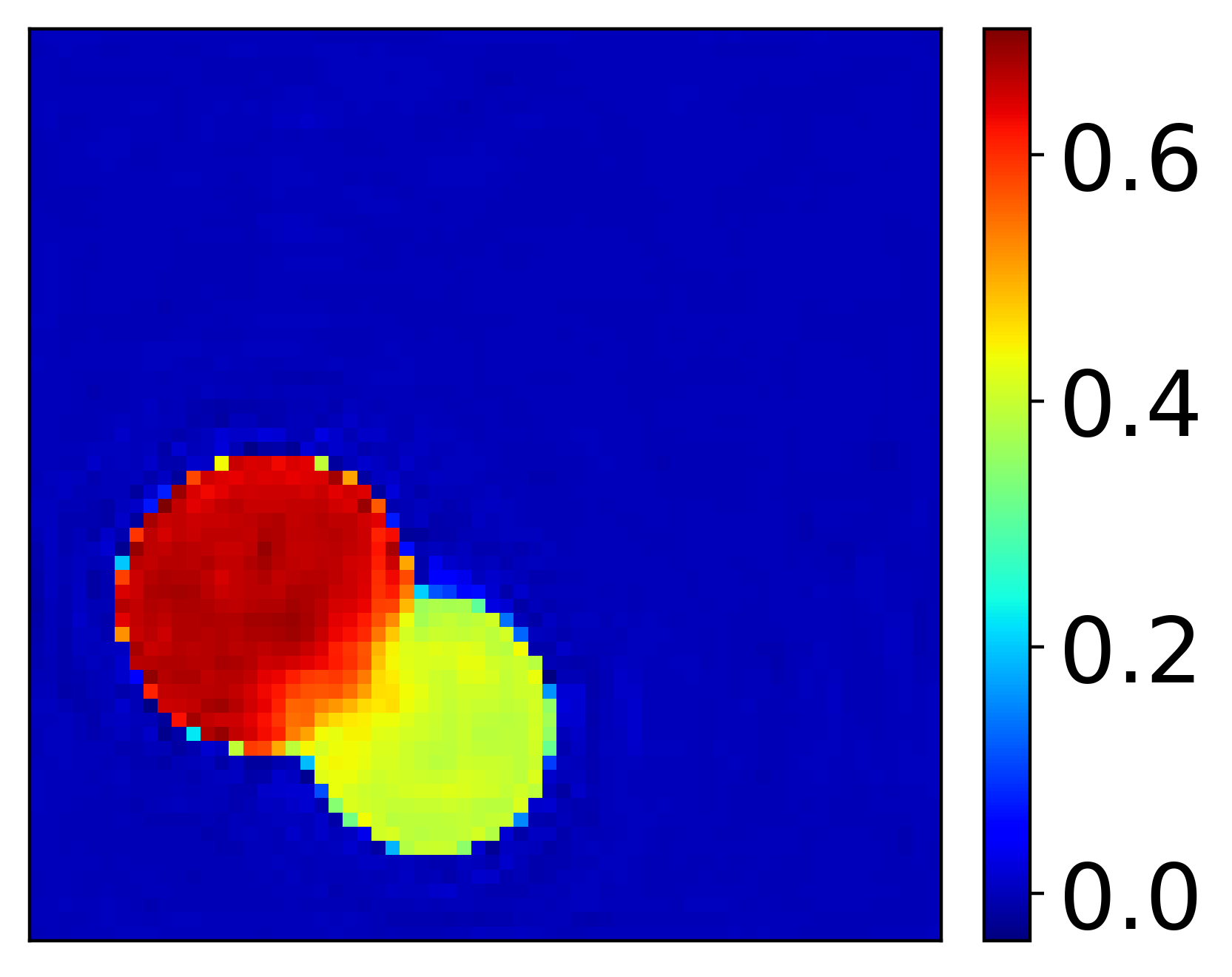} &
			\includegraphics[valign=m,width=0.15\textwidth]{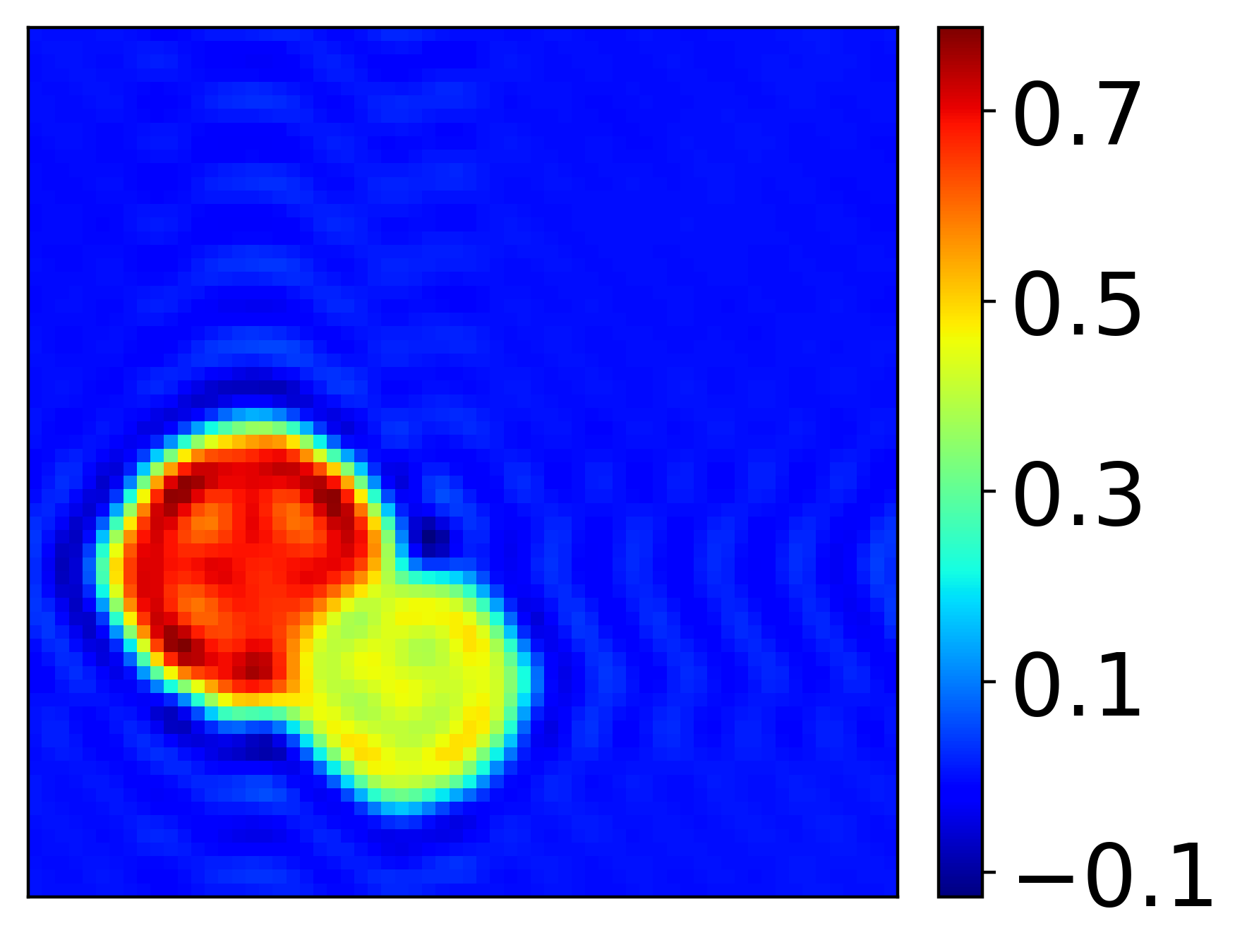} & 5\% \\
			& \includegraphics[valign=m,width=0.15\textwidth]{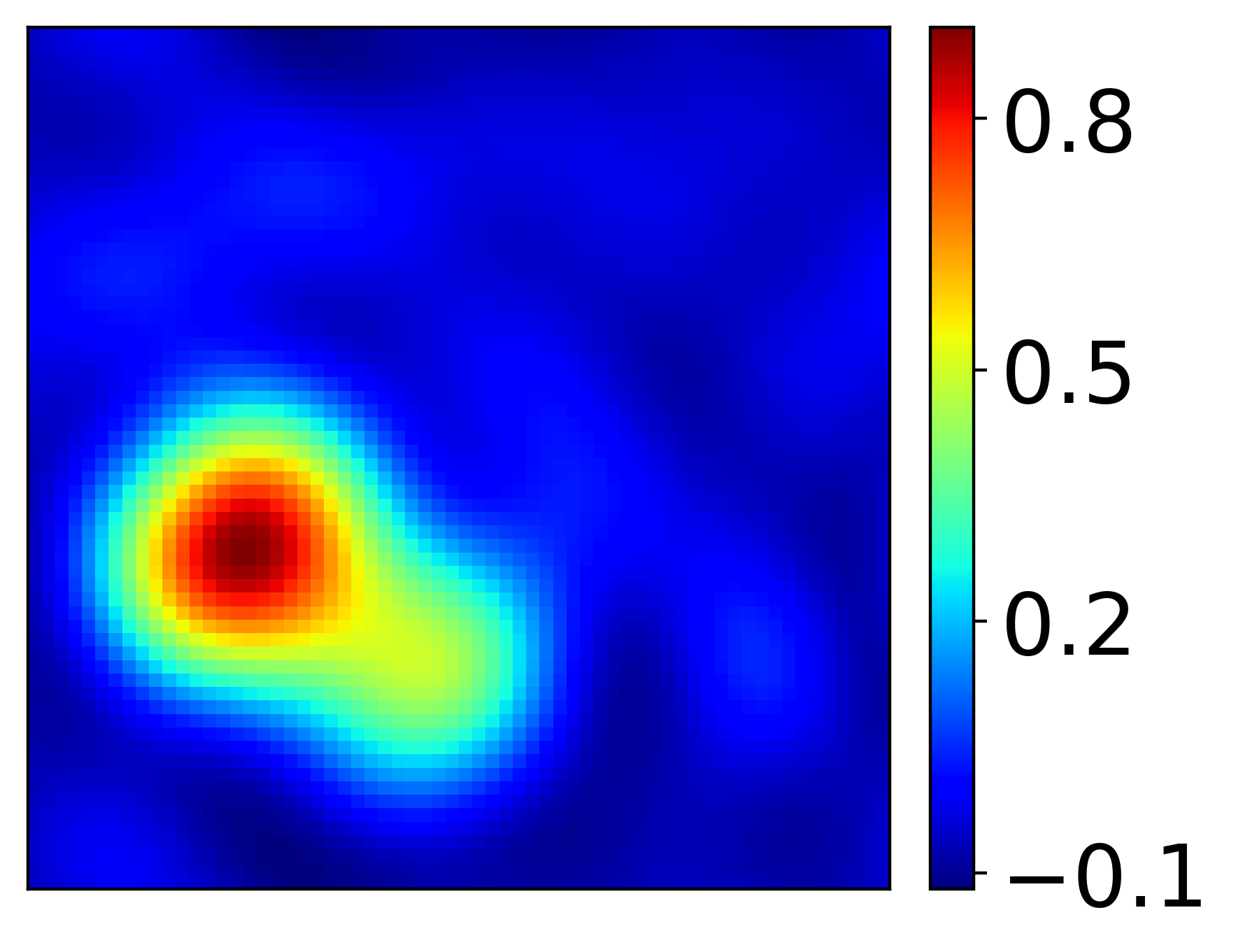} &
			\includegraphics[valign=m,width=0.15\textwidth]{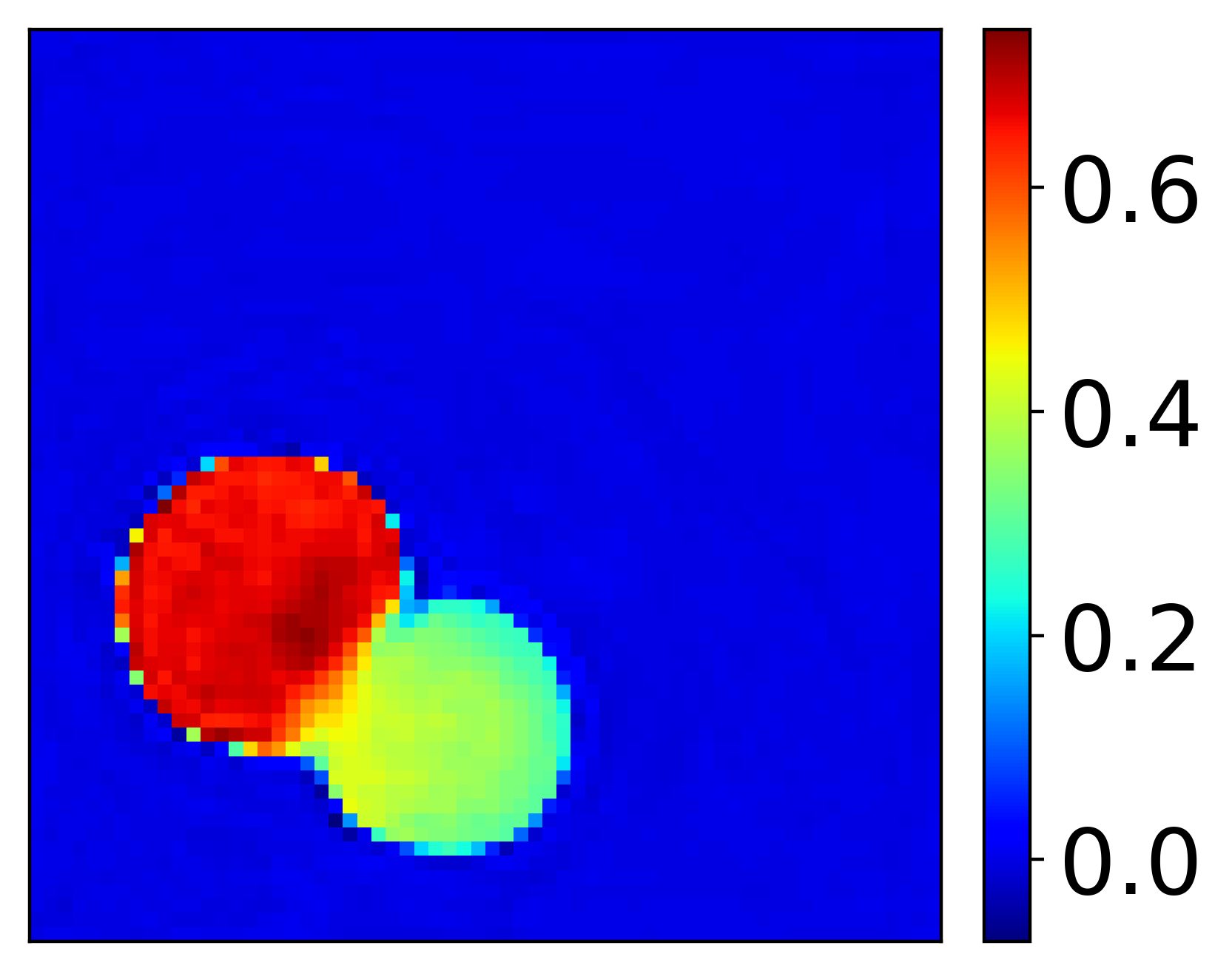} &
			\includegraphics[valign=m,width=0.15\textwidth]{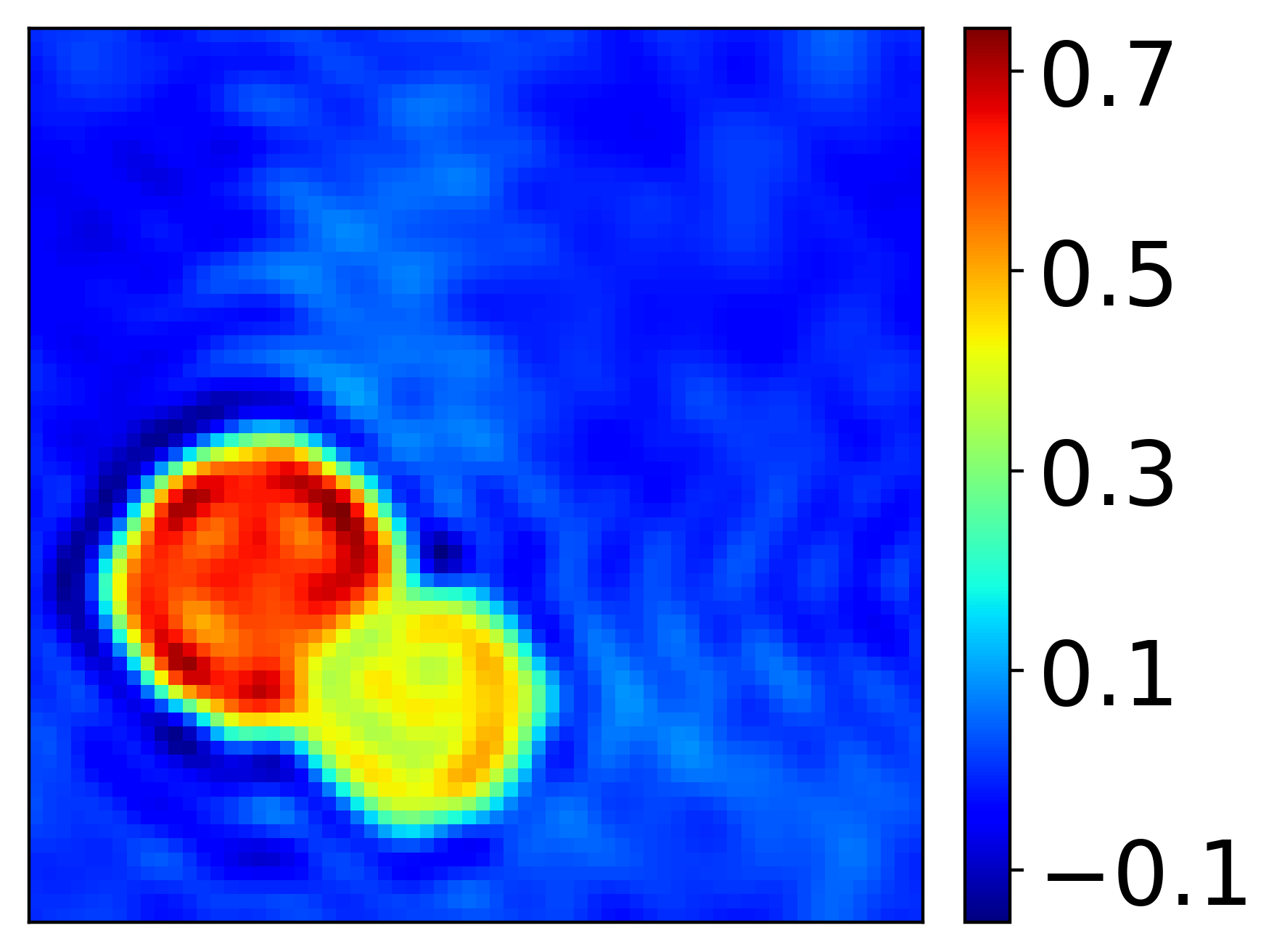} & 50\% \\
			& \includegraphics[valign=m,width=0.15\textwidth]{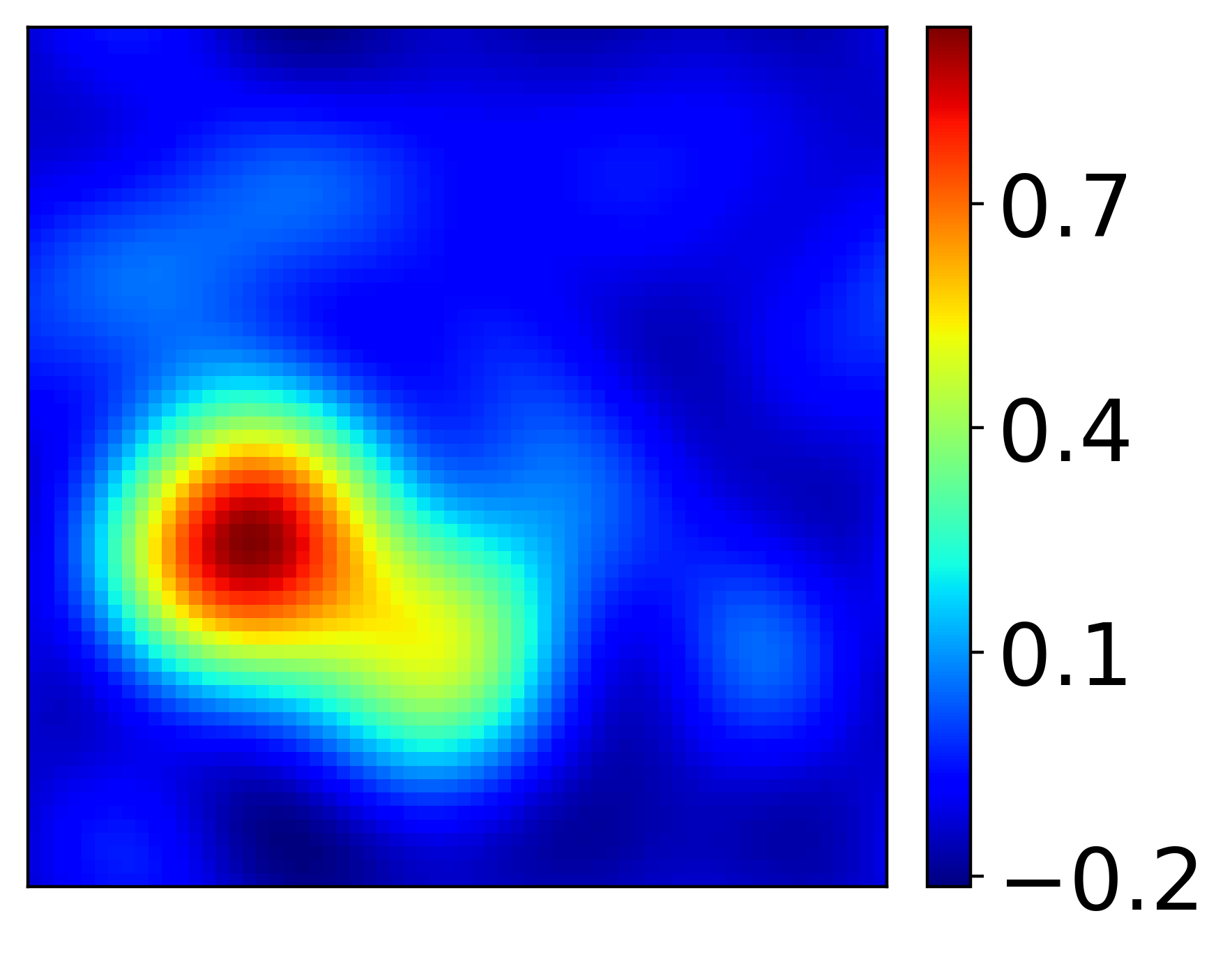} &
			\includegraphics[valign=m,width=0.15\textwidth]{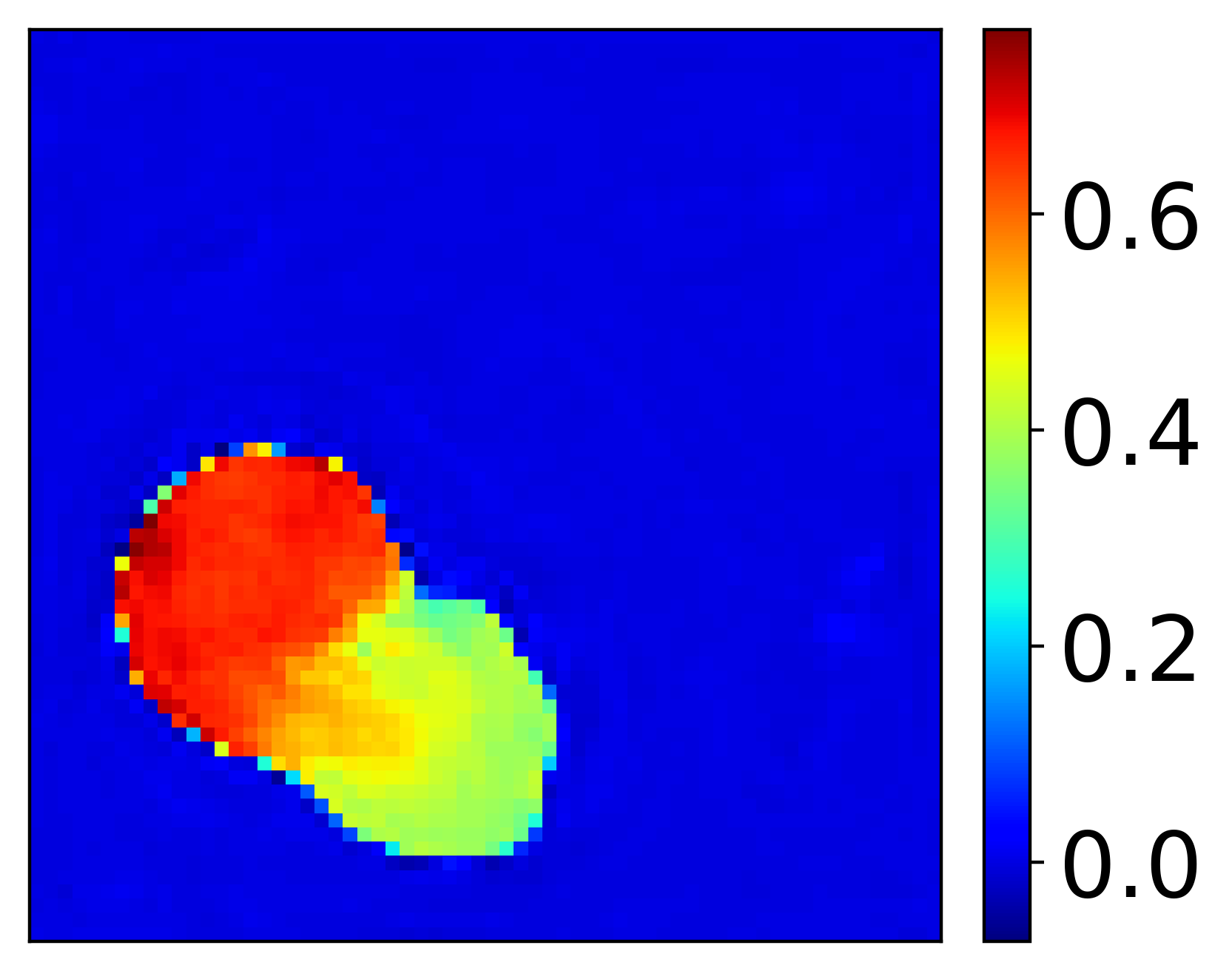} &
			\includegraphics[valign=m,width=0.15\textwidth]{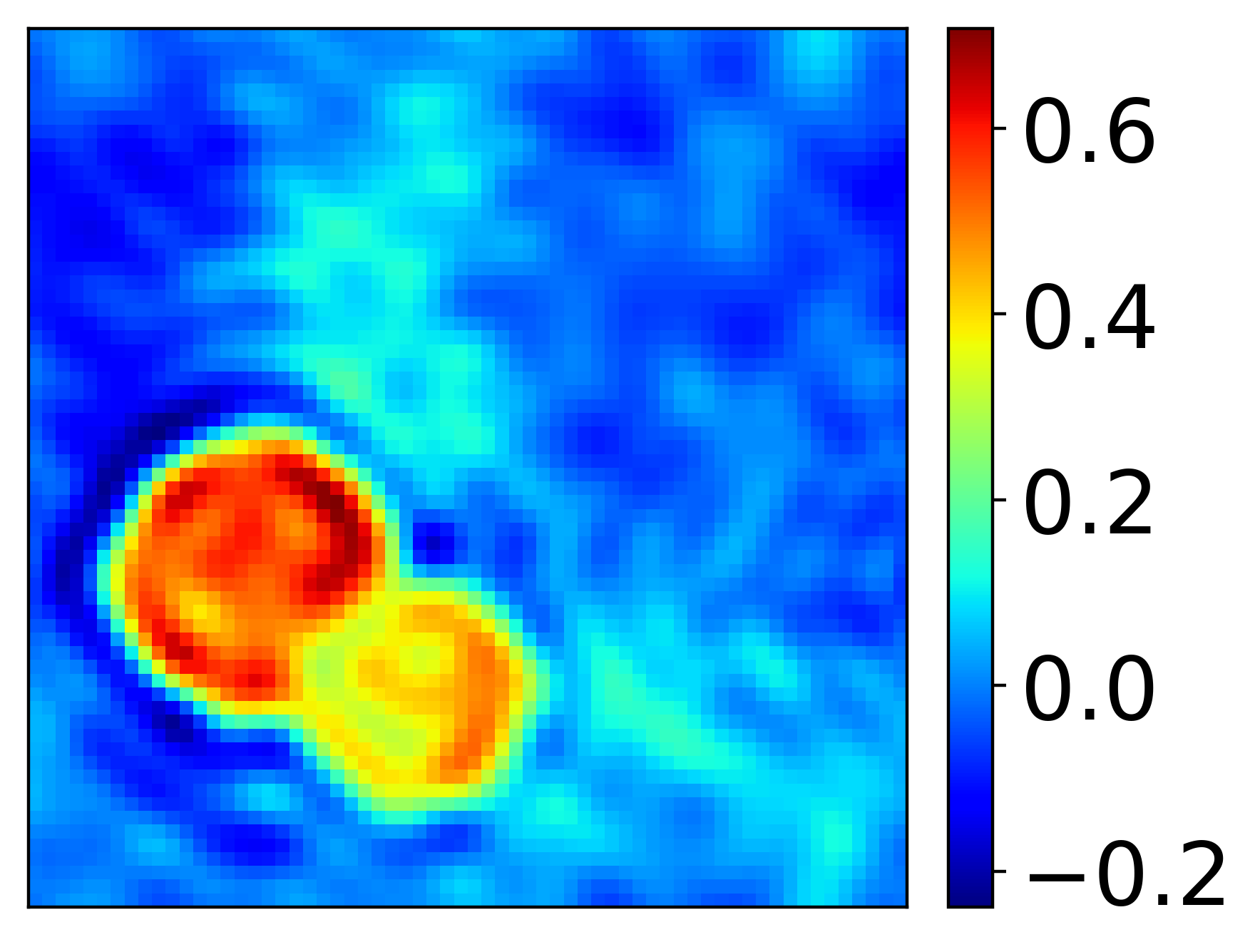} & 100\% \\
			\includegraphics[valign=m,width=0.15\textwidth]{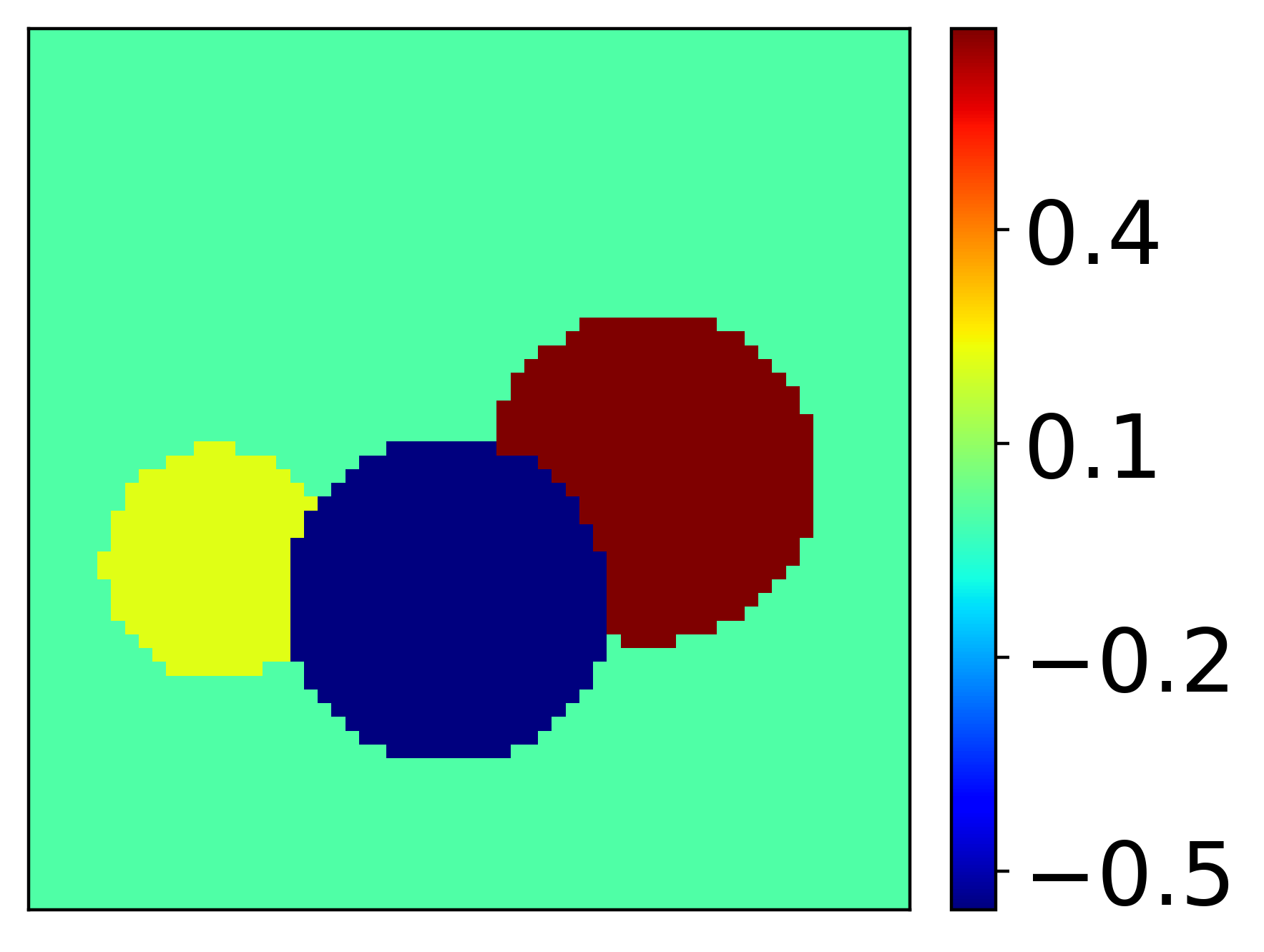} &
			\includegraphics[valign=m,width=0.15\textwidth]{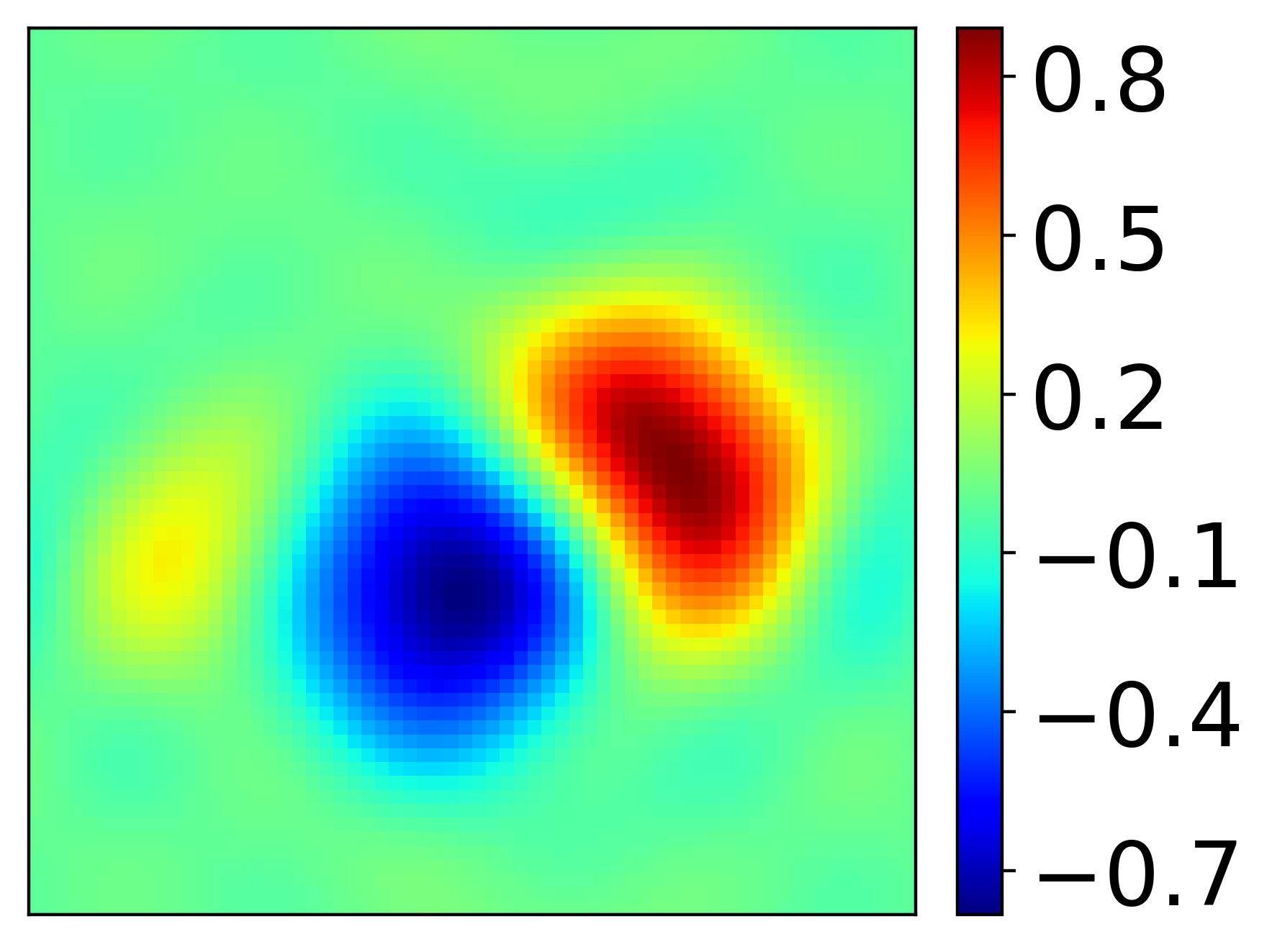} &
			\includegraphics[valign=m,width=0.15\textwidth]{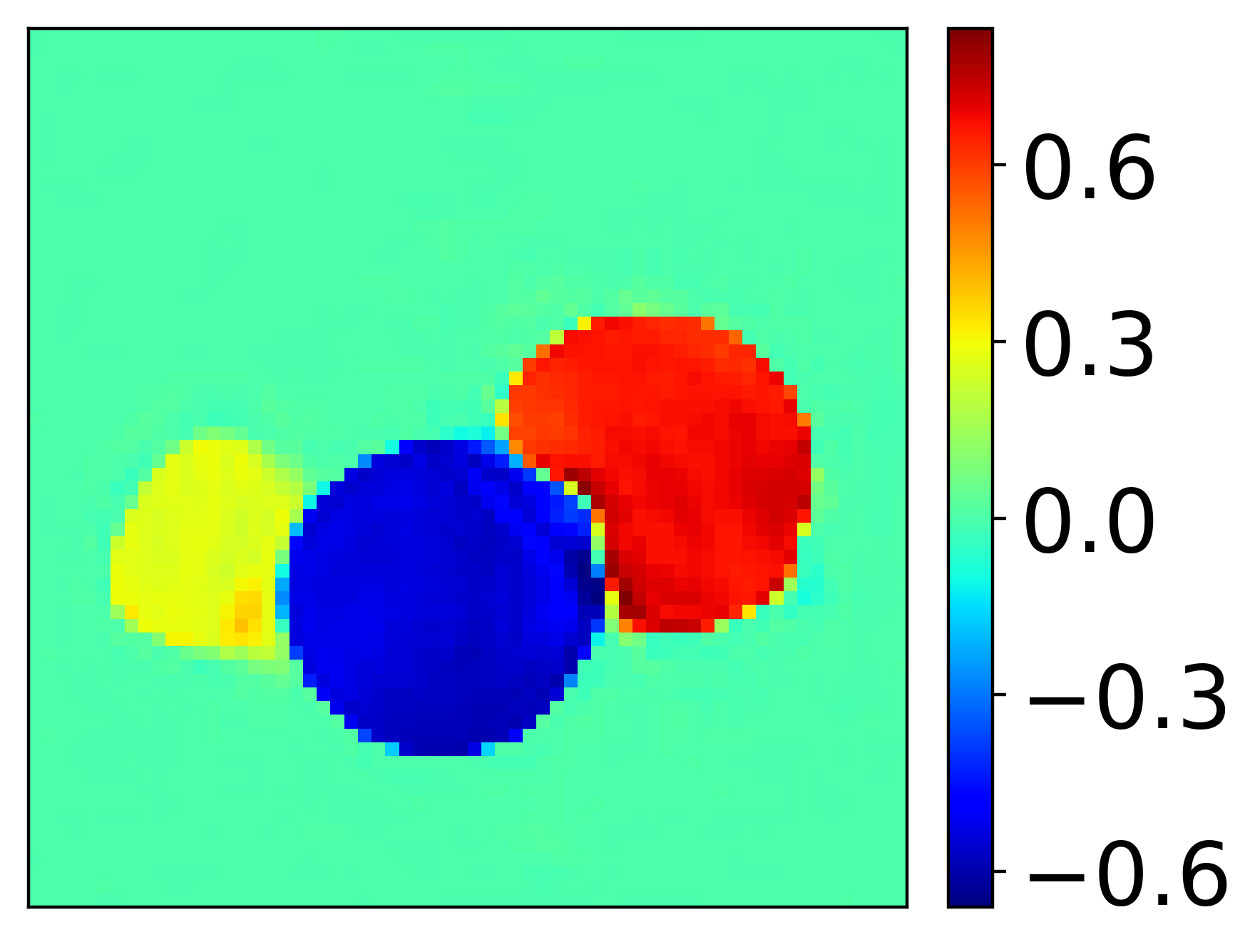} &
			\includegraphics[valign=m,width=0.15\textwidth]{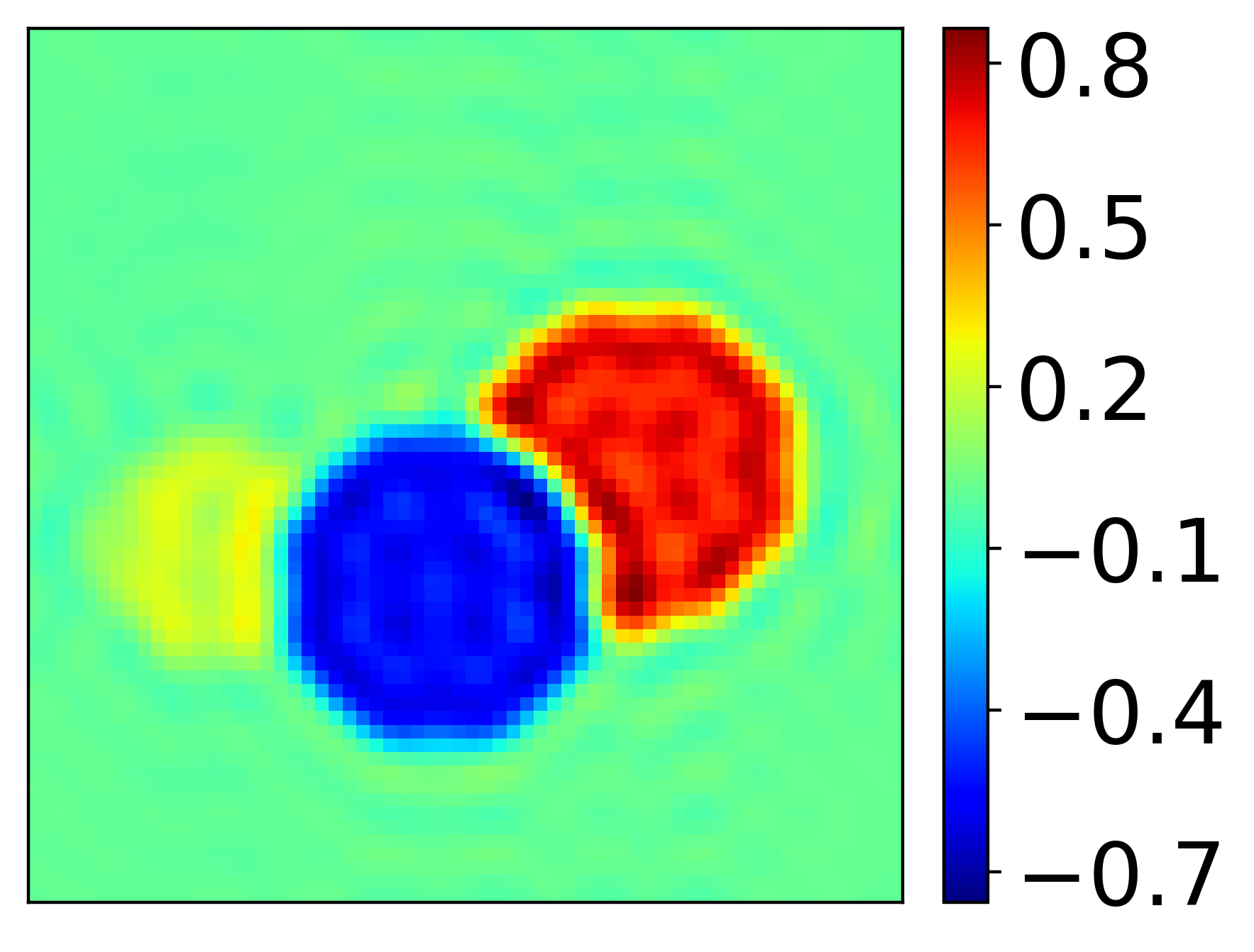} & 5\% \\
			& \includegraphics[valign=m,width=0.15\textwidth]{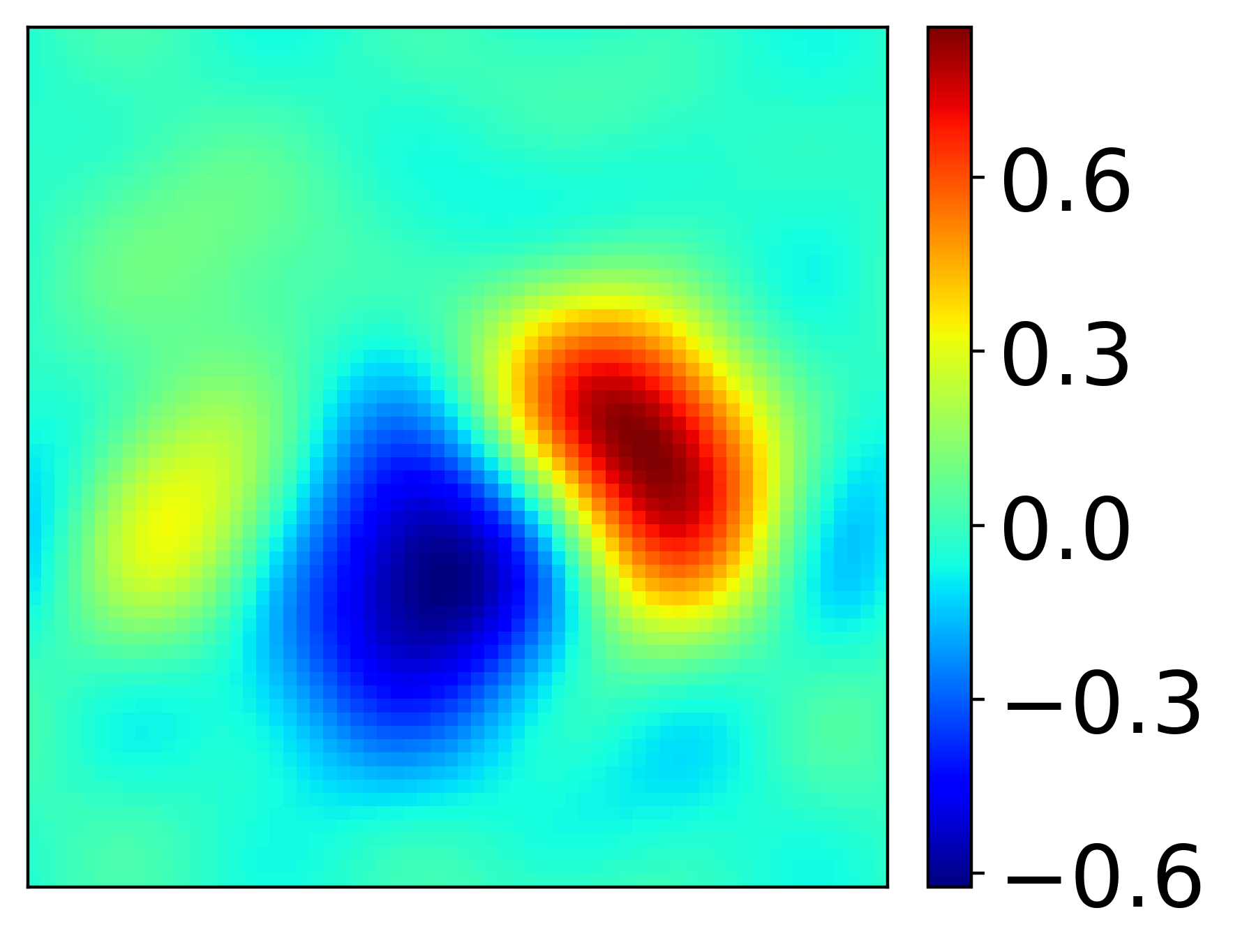} &
			\includegraphics[valign=m,width=0.15\textwidth]{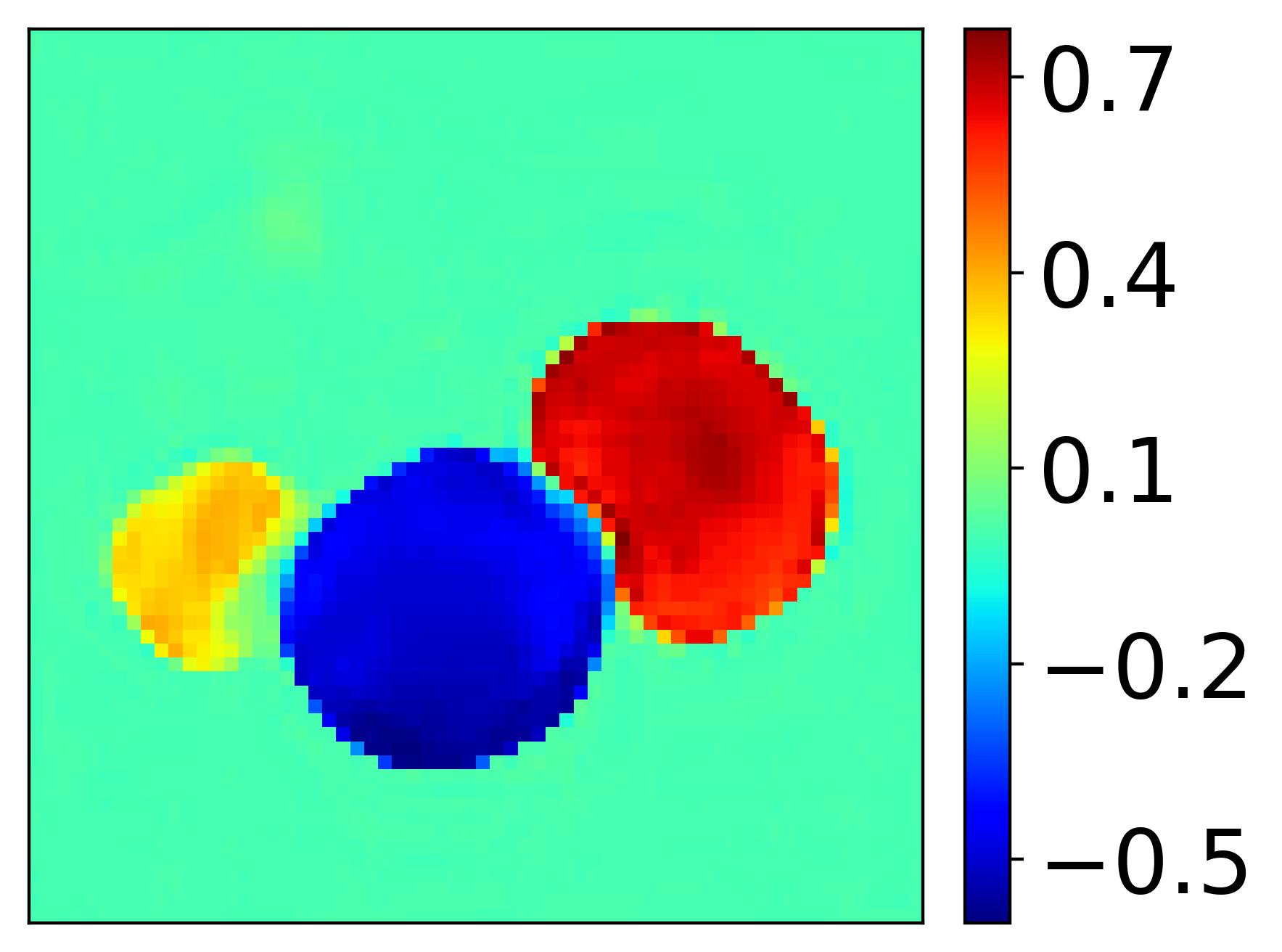} &
			\includegraphics[valign=m,width=0.15\textwidth]{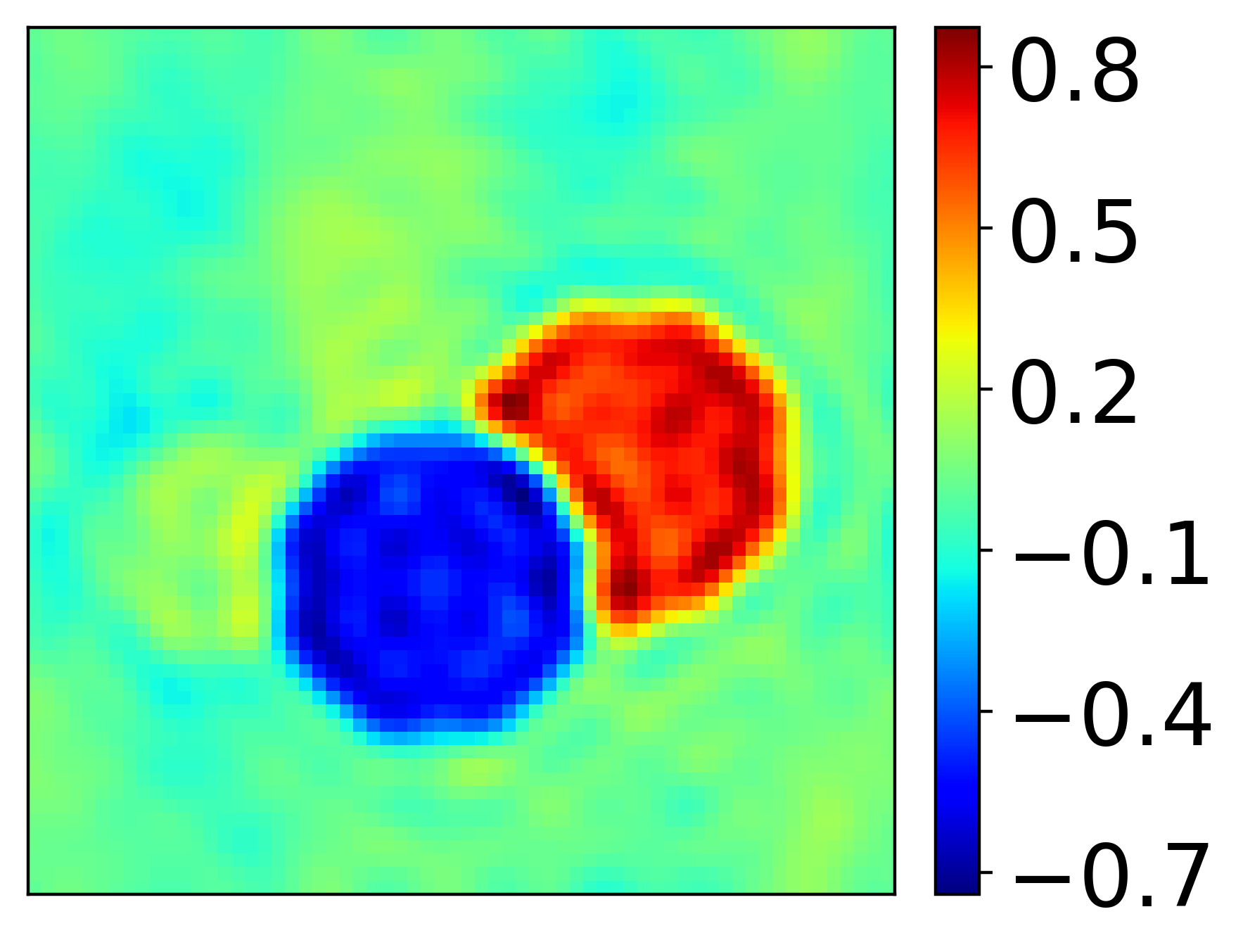} & 50\% \\
			& \includegraphics[valign=m,width=0.15\textwidth]{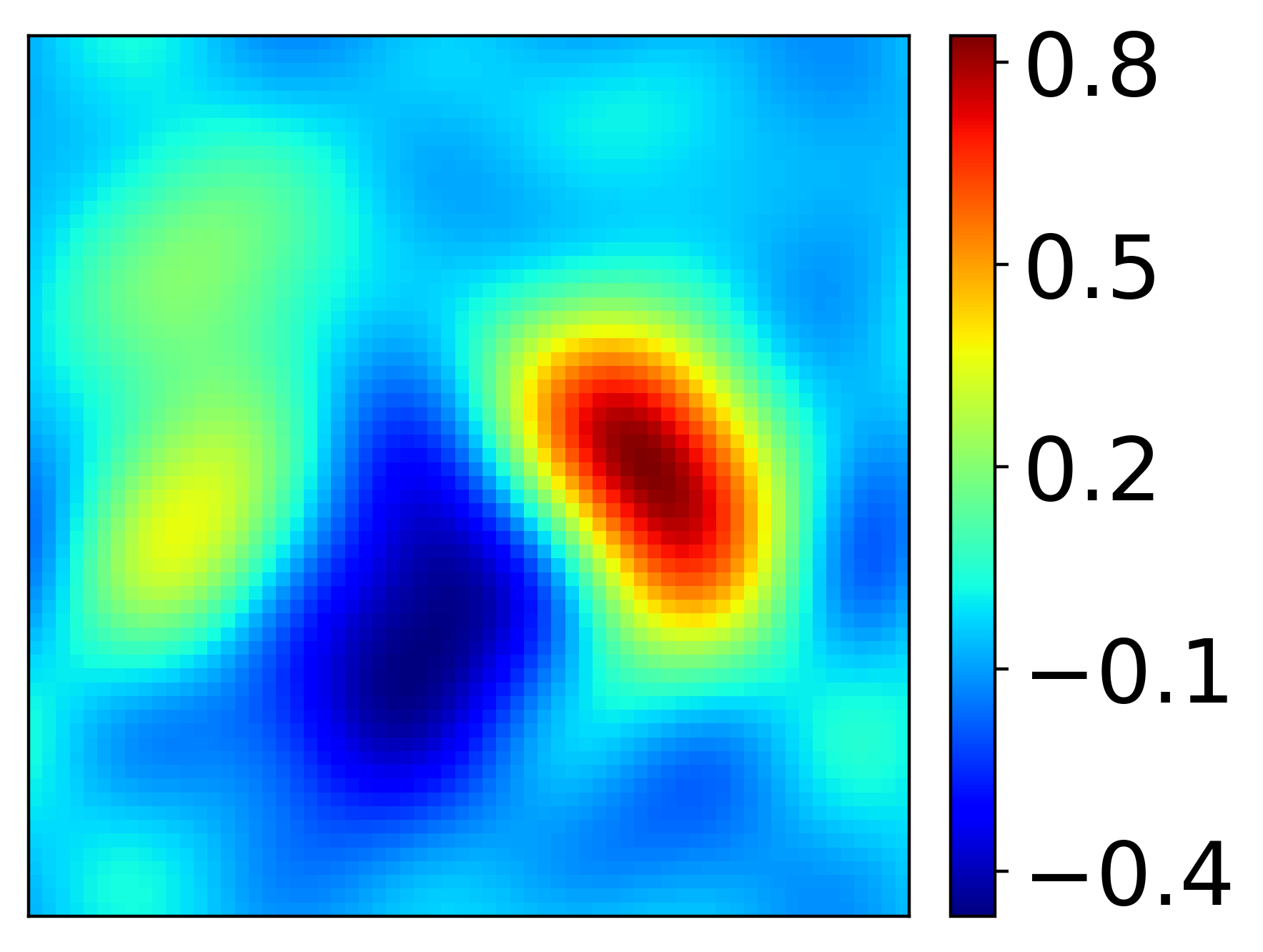} &
			\includegraphics[valign=m,width=0.15\textwidth]{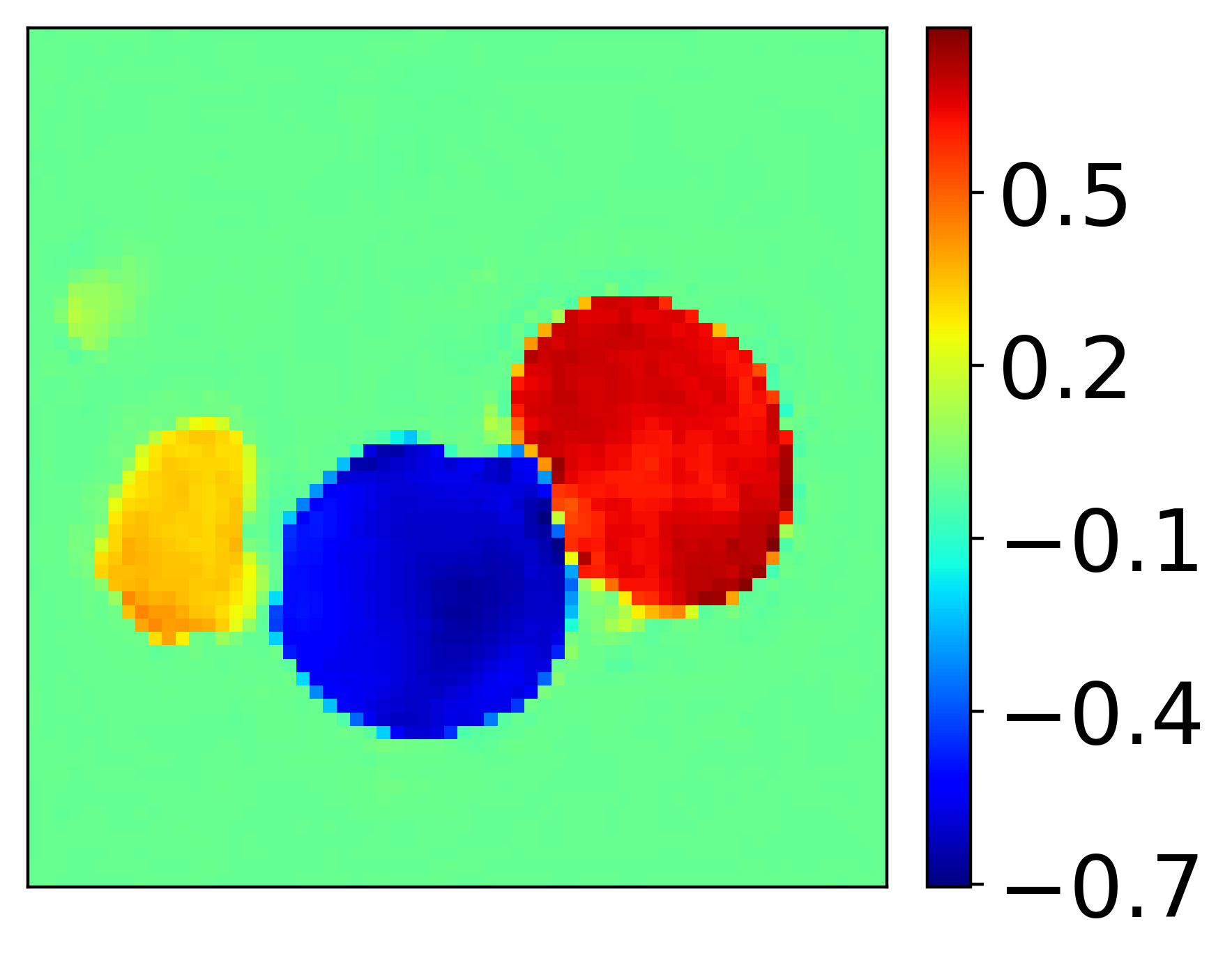} &
			\includegraphics[valign=m,width=0.15\textwidth]{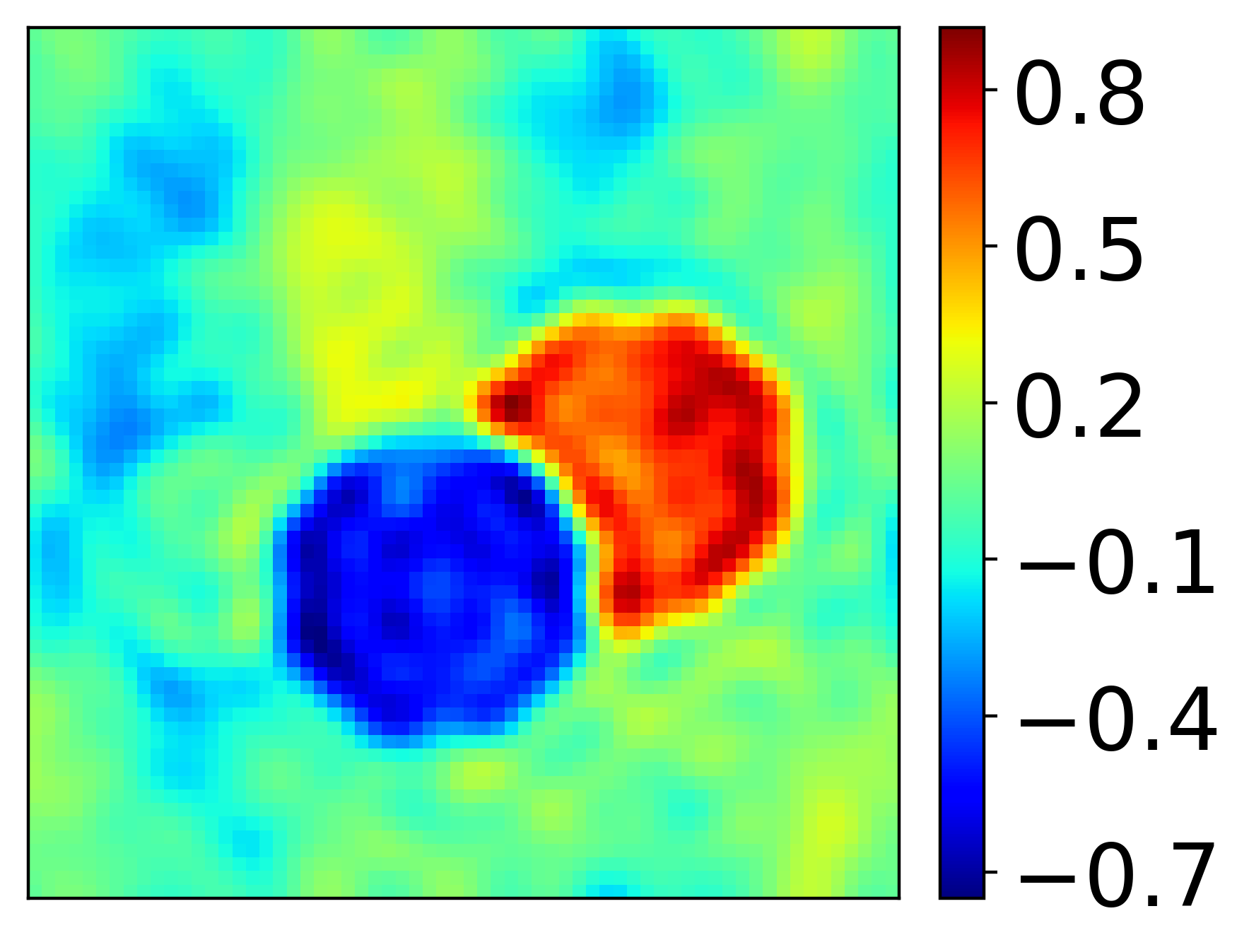} & 100\% \\
		\end{tabular}
	\end{center}
	\setlength{\abovecaptionskip}{-5pt}
	\caption{Reconstruction results for disk sources under $5\%$, $50\%$, and $100\%$ noise levels. Columns from left to right: ground-truth source, classical Fourier reconstruction ($N=3$), U-Net-enhanced reconstruction ($N=3$), and classical Fourier reconstruction ($N=10$).}
	\label{fig:Circle}
\end{figure}
Reconstruction results for noise levels of $5\%$, $50\%$, and $100\%$ are illustrated in \Cref{fig:Circle}, alongside a comparison with the classical Fourier method ($N=10$). As noted in \cite{zhang2015fourier}, the Fourier method typically suffers from the Gibbs phenomenon when reconstructing piecewise constant sources. Our hybrid approach effectively suppresses these artifacts, resolving sharp boundaries and overlapping regions clearly. 
\begin{table}[htbp]
	\centering
	\begin{tabular}{l|ccc}
		\hline
		Noise level      & $5\%$    & $50\%$    & $100\%$   \\
		\hline
		Fourier ($N=10$) & $9.24\%$ & $12.80\%$ & $23.58\%$ \\
		U-Net-enhanced           & $5.35\%$ & $8.48\%$  & $9.62\%$  \\
		\hline
	\end{tabular}
	\caption{NMSE comparison for disk source reconstruction (Example \ref{example:disks}) across various noise levels.}
	\label{tab:disks}
\end{table}
Quantitative comparisons in \Cref{tab:disks} indicate that the U-Net achieves lower errors and superior noise tolerance compared to the baseline, even when the latter utilizes significantly more measurements ($N=10$). Furthermore, the deep-learning-enhanced reconstructions exhibit a cleaner background and significantly reduced oscillatory artifacts relative to purely spectral methods.

\subsection{Transfer Learning Verification: High-to-Low Noise Strategy}
\begin{figure}[htbp]
	\centering
	\subcaptionbox{}{
		\includegraphics[width=0.25\linewidth]{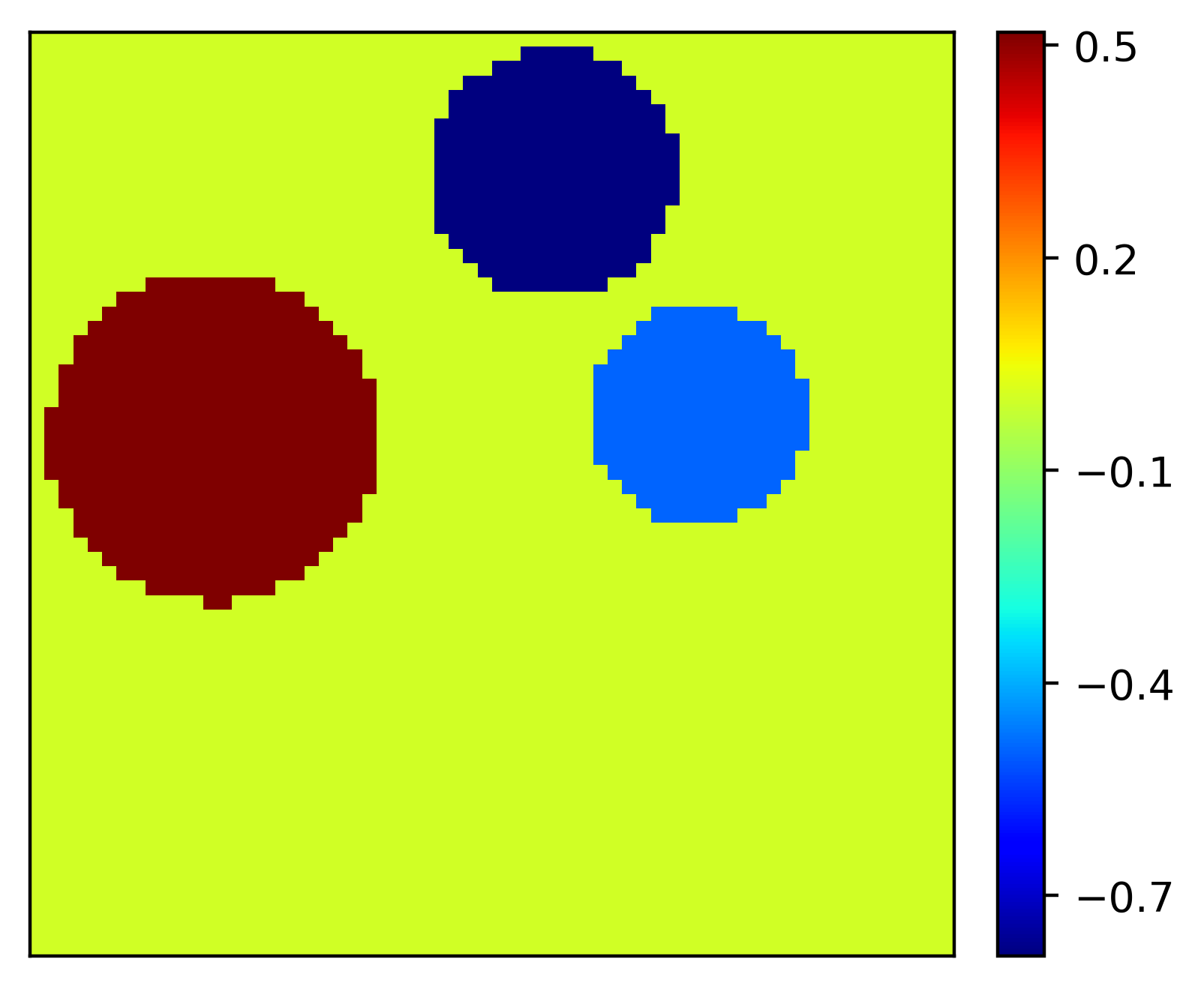}
	}
	\subcaptionbox{}{
		\includegraphics[width=0.25\linewidth]{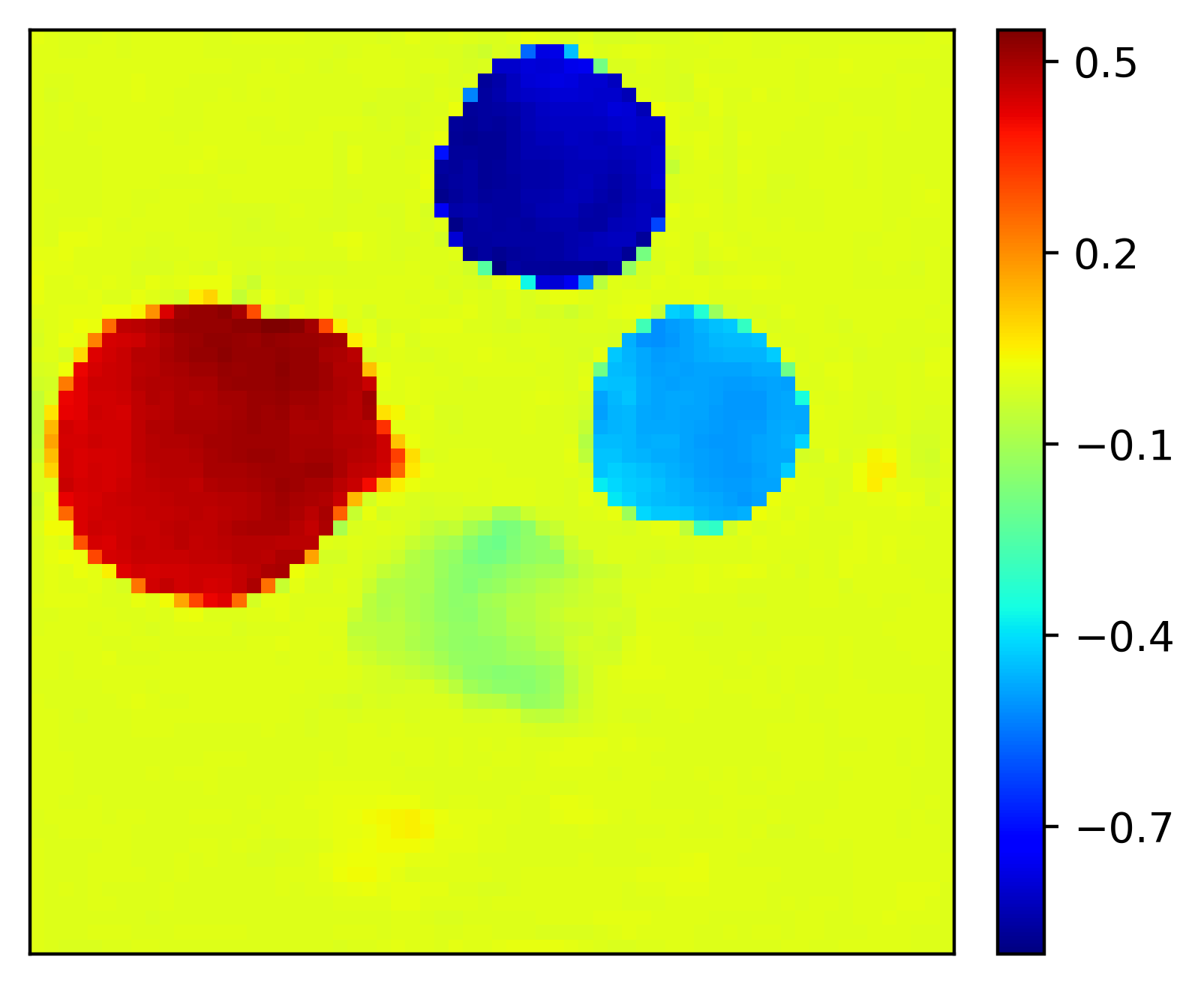}
	}
	\subcaptionbox{}{
		\includegraphics[width=0.25\linewidth]{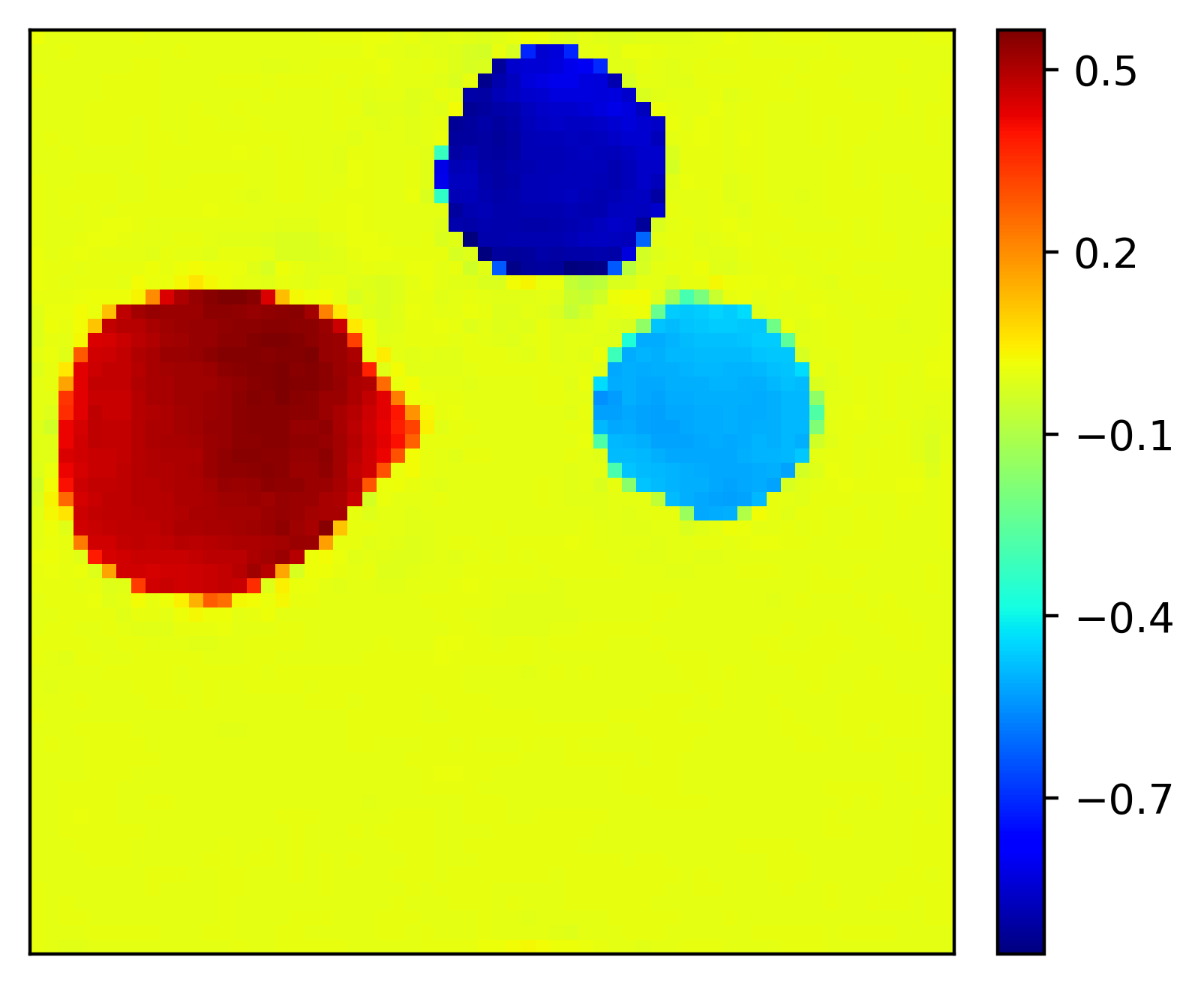}
	}\\
	\caption{Generalization performance of the network trained on $100\%$ noise when applied to lower-noise data. From left to right: ground-truth source, reconstruction for $5\%$ noise level, and reconstruction for $50\%$ noise level.}
	\label{fig:reconst_100}
\end{figure}
As illustrated in \Cref{fig:reconst_100}, the U-Net model $\mathcal{G}_\Theta$ trained exclusively on $100\%$ noise data exhibits strong generalization when directly applied to $5\%$ and $50\%$ noise levels. This suggests that training under extreme noise conditions encourages the network to learn structural priors that remain valid across varying data fidelities. Building on this observation, we investigate the high-to-low transfer learning strategy: we use the $100\%$ noise model as a pre-trained initialization for lower-noise regimes. For fine-tuning, the initial learning rate is set to $0.0005$ with a decay factor of $0.5$ every $5$ epochs, and training proceeds for only $30$ epochs while keeping other hyperparameters constant.
\begin{figure}[htbp]
	\centering
	\subcaptionbox{}{
		\adjustbox{valign=t}{\includegraphics[width=0.18\linewidth]{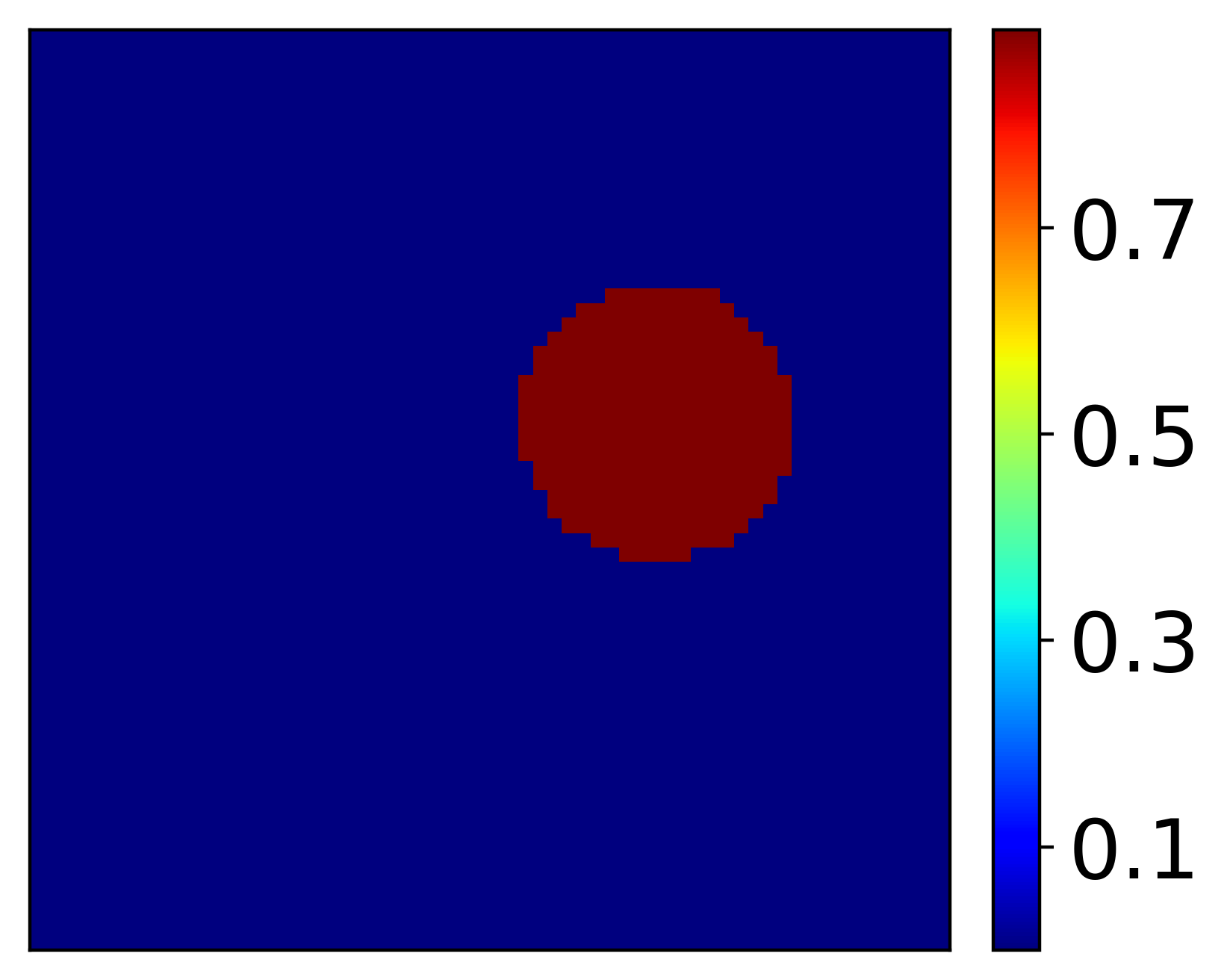}}
	}
	\hspace{-0.4cm}
	\subcaptionbox{}{
		\adjustbox{valign=t}{\includegraphics[width=0.18\linewidth]{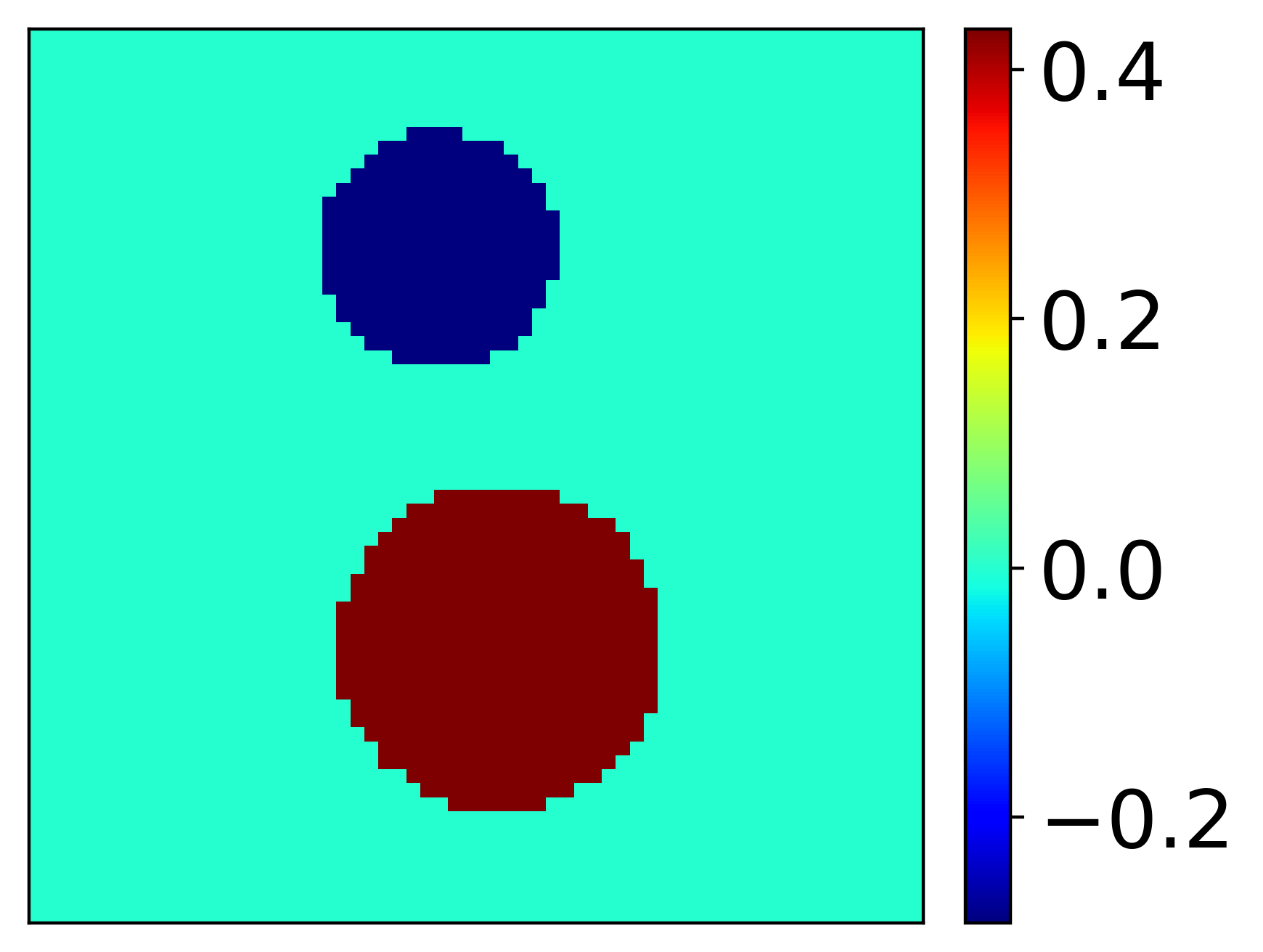}}
	}
	\hspace{-0.4cm}
	\subcaptionbox{}{
		\adjustbox{valign=t}{\includegraphics[width=0.18\linewidth]{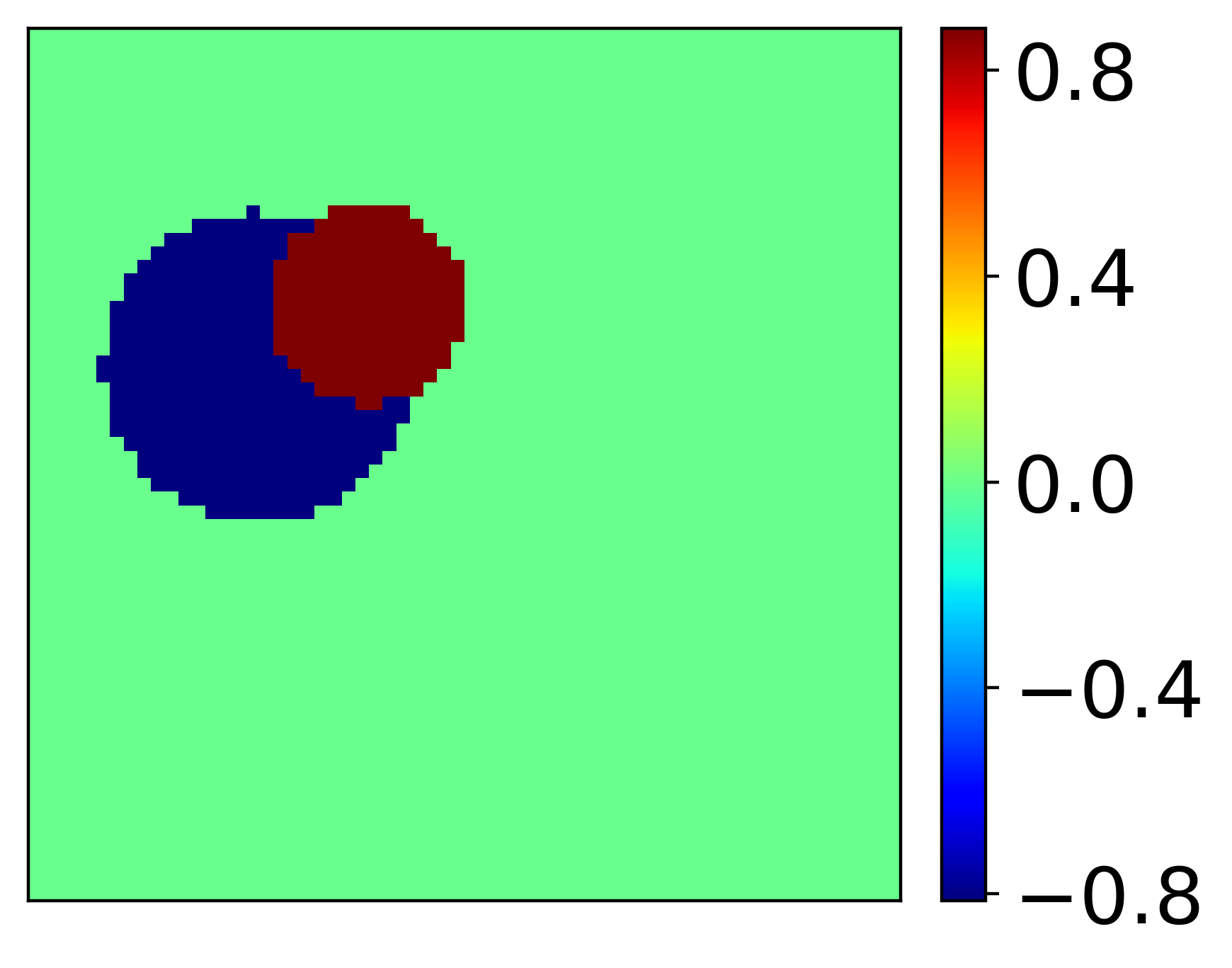}}
	}
	\hspace{-0.4cm}
	\subcaptionbox{}{
		\adjustbox{valign=t}{\includegraphics[width=0.18\linewidth]{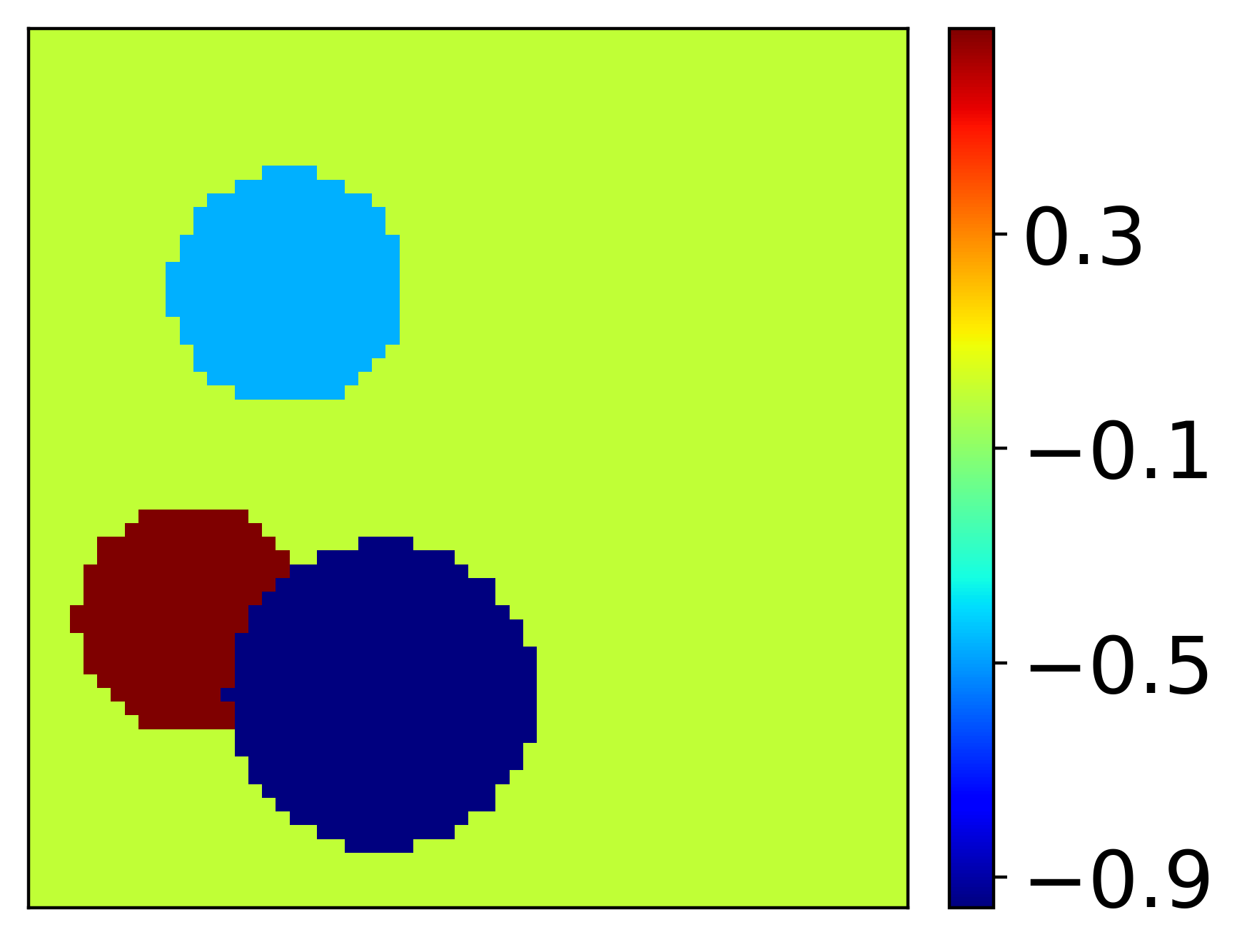}}
	}
	\hspace{-0.4cm}
	\subcaptionbox{}{
		\adjustbox{valign=t}{\includegraphics[width=0.18\linewidth]{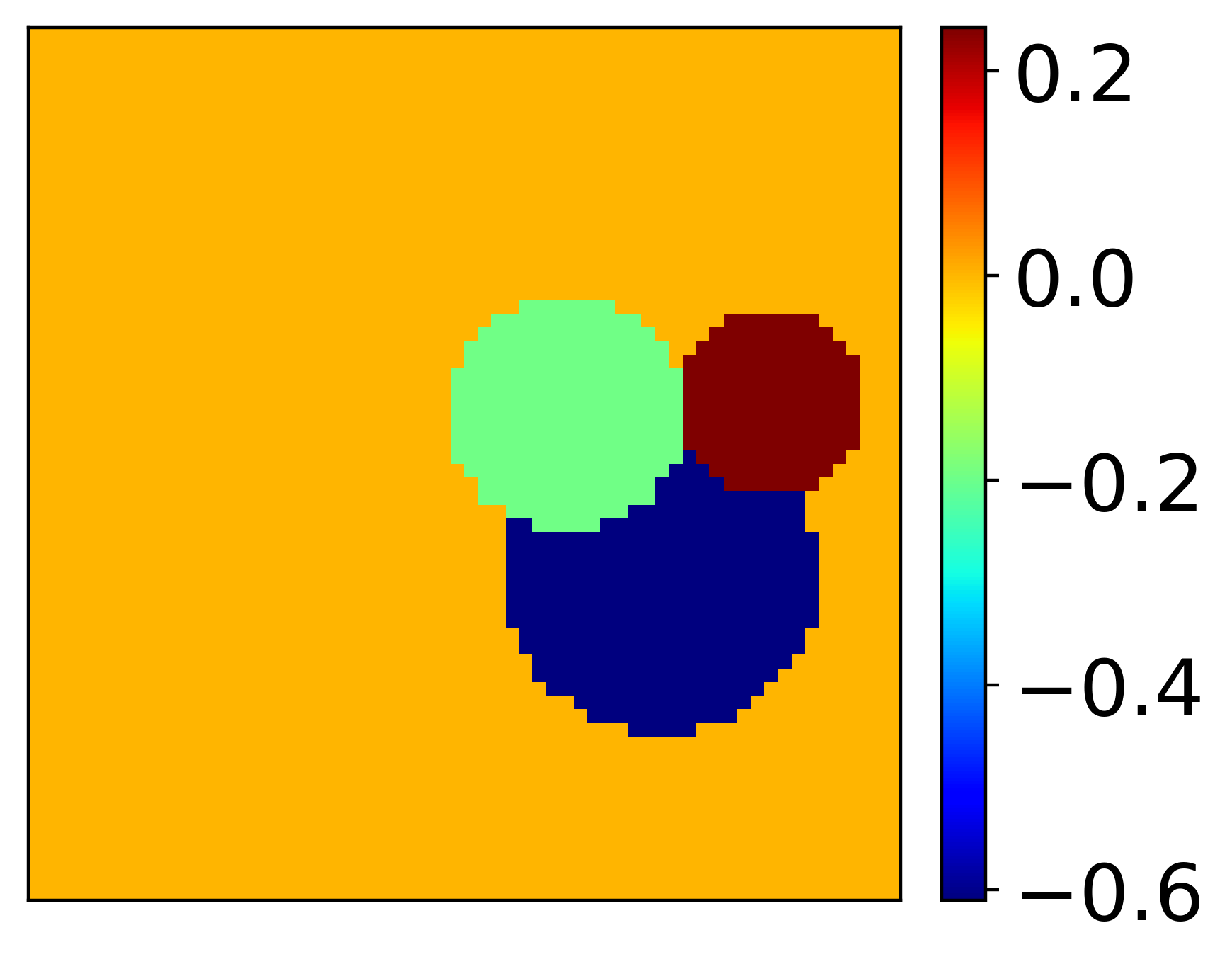}}
	}\\
	\subcaptionbox{}{
		\adjustbox{valign=t}{\includegraphics[width=0.18\linewidth]{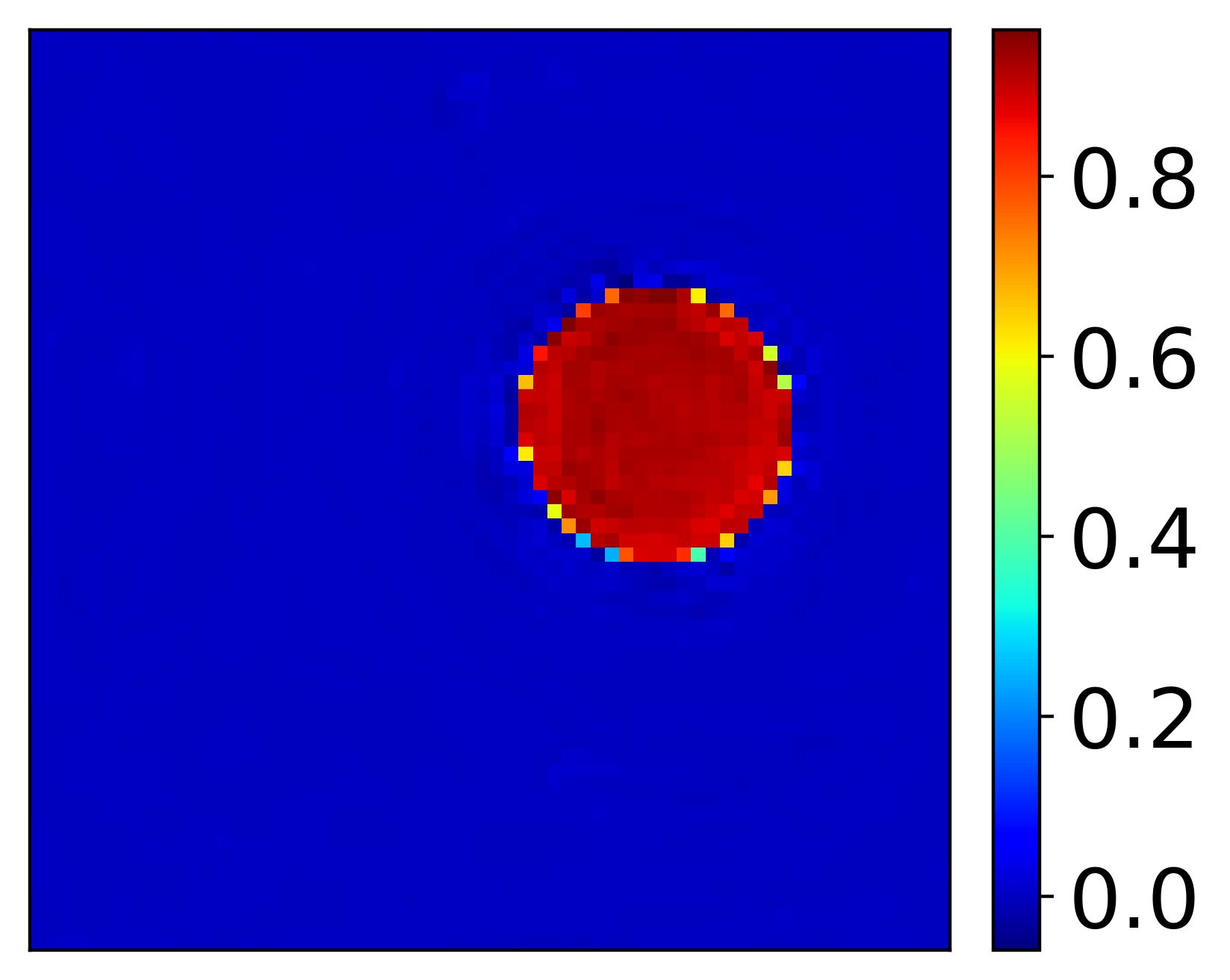}}
	}
	\hspace{-0.4cm}
	\subcaptionbox{}{
		\adjustbox{valign=t}{\includegraphics[width=0.18\linewidth]{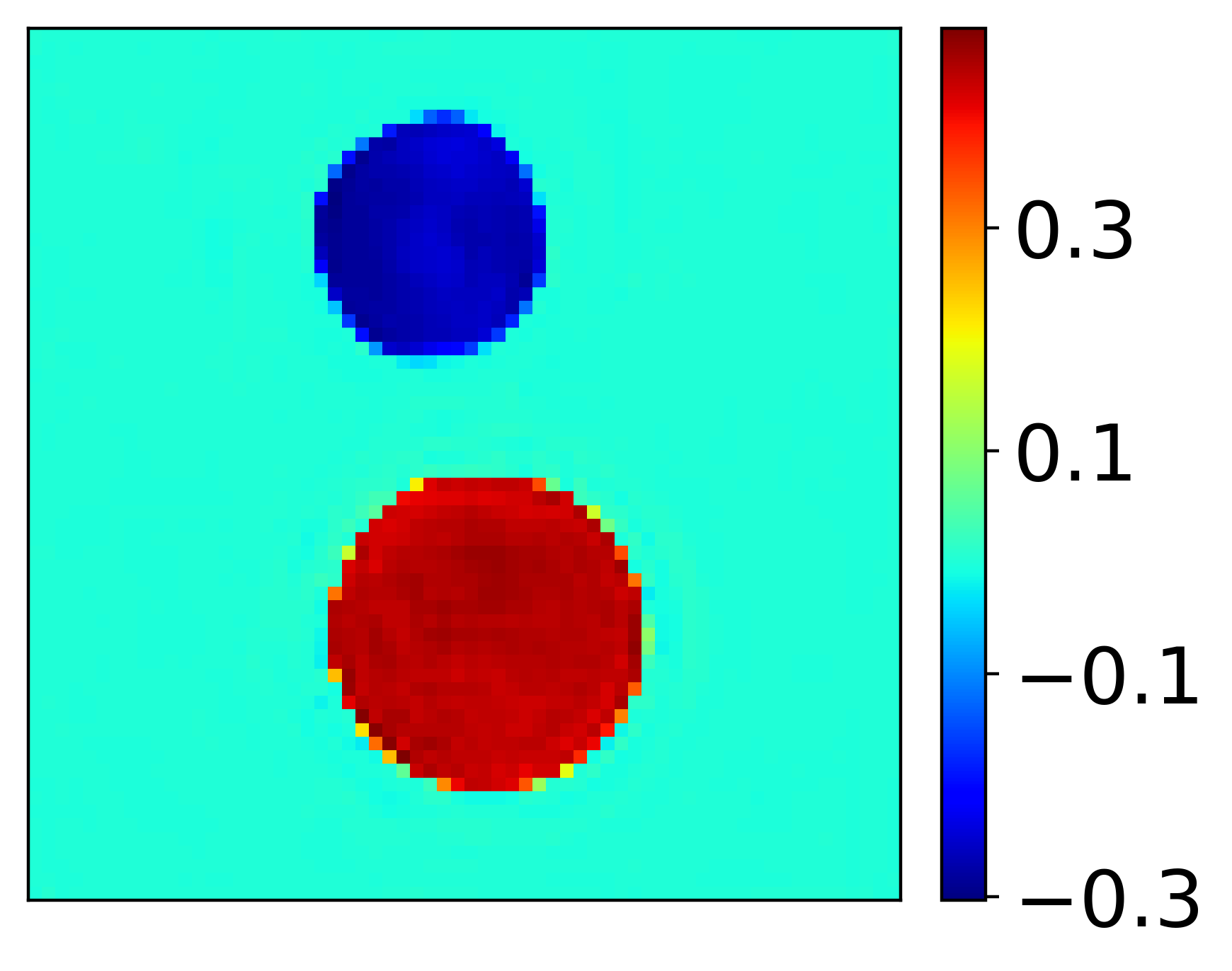}}
	}
	\hspace{-0.4cm}
	\subcaptionbox{}{
		\adjustbox{valign=t}{\includegraphics[width=0.18\linewidth]{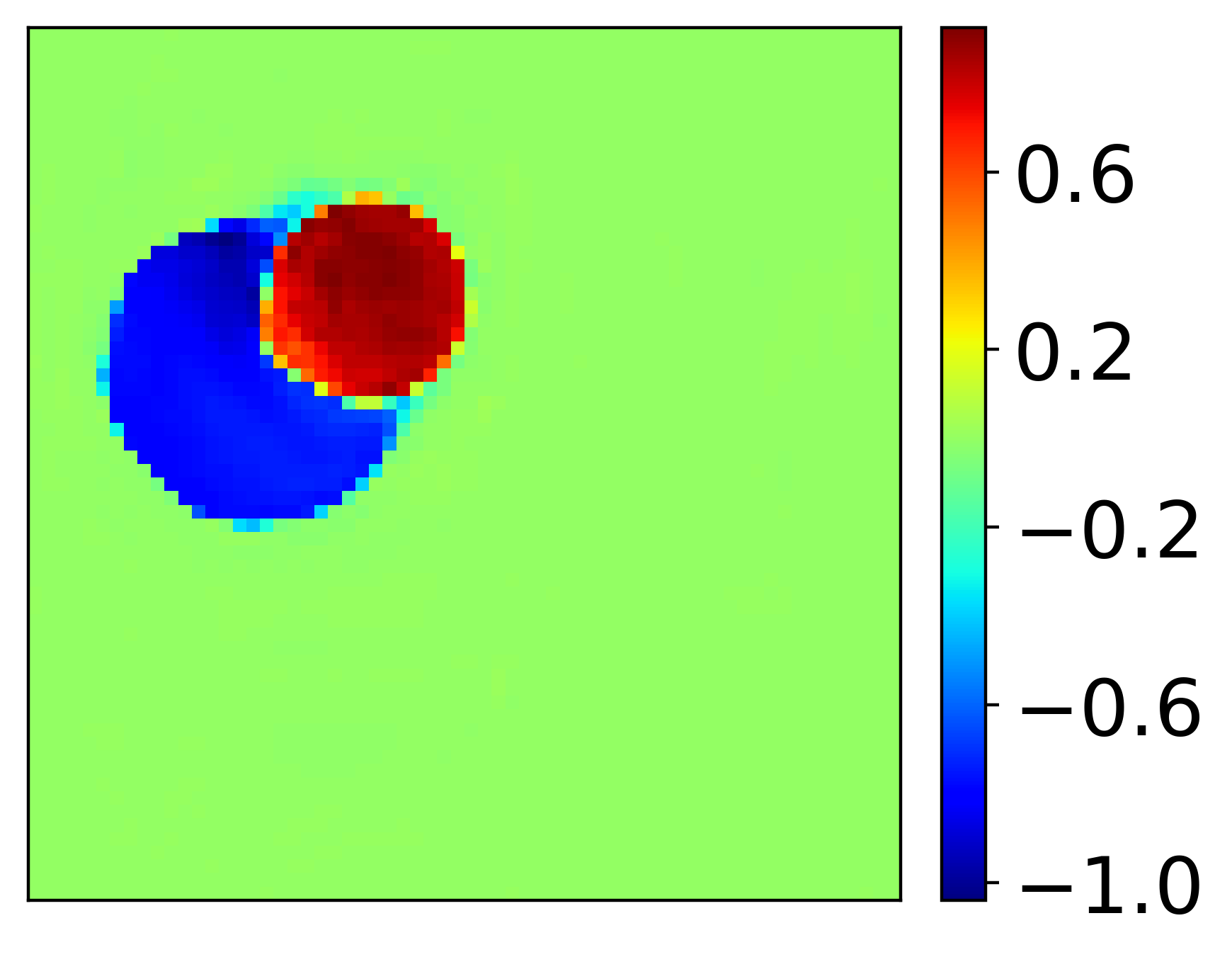}}
	}
	\hspace{-0.4cm}
	\subcaptionbox{}{
		\adjustbox{valign=t}{\includegraphics[width=0.18\linewidth]{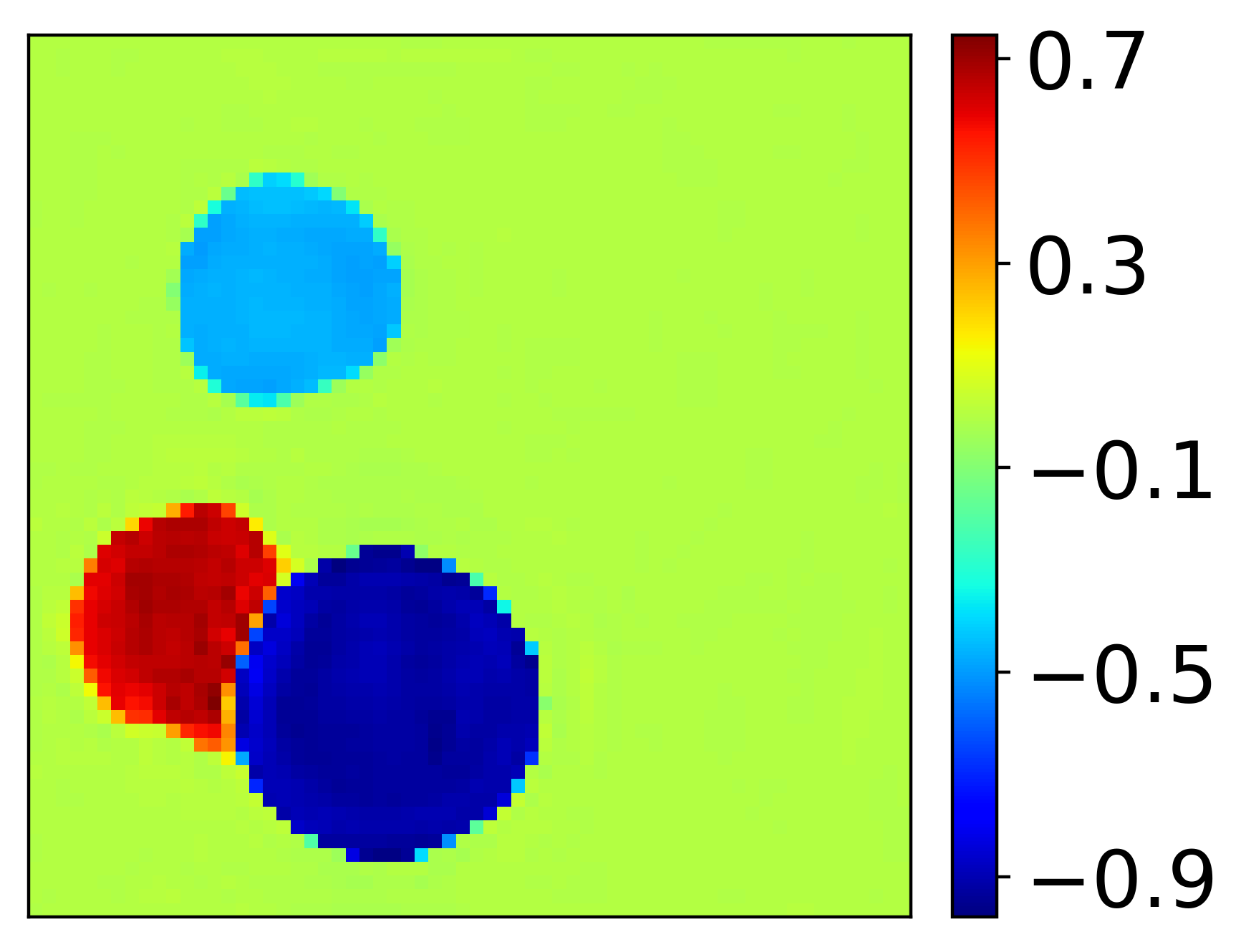}}
	}
	\hspace{-0.4cm}
	\subcaptionbox{}{
		\adjustbox{valign=t}{\includegraphics[width=0.18\linewidth]{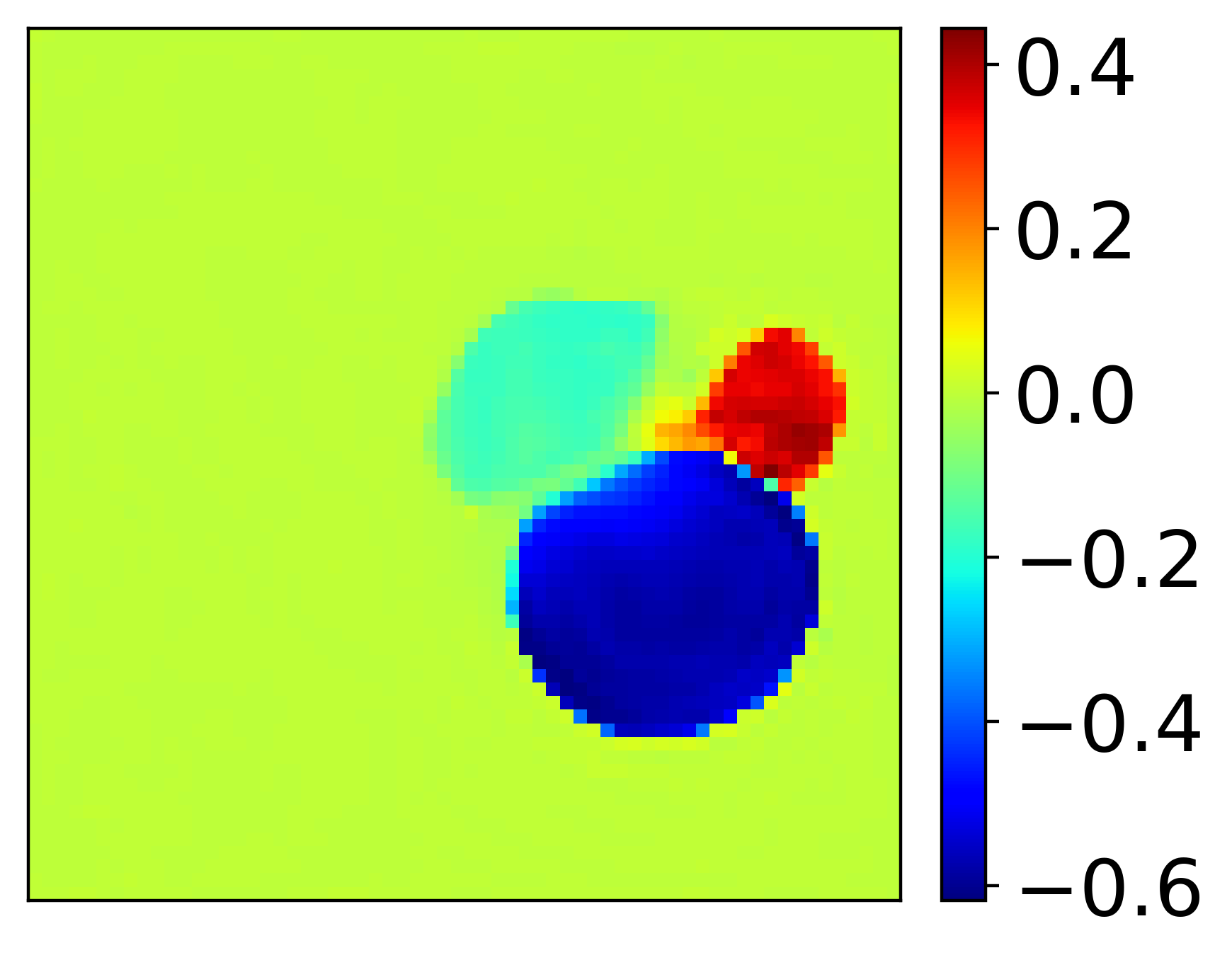}}
	}\\
	\subcaptionbox{}{
		\adjustbox{valign=t}{\includegraphics[width=0.18\linewidth]{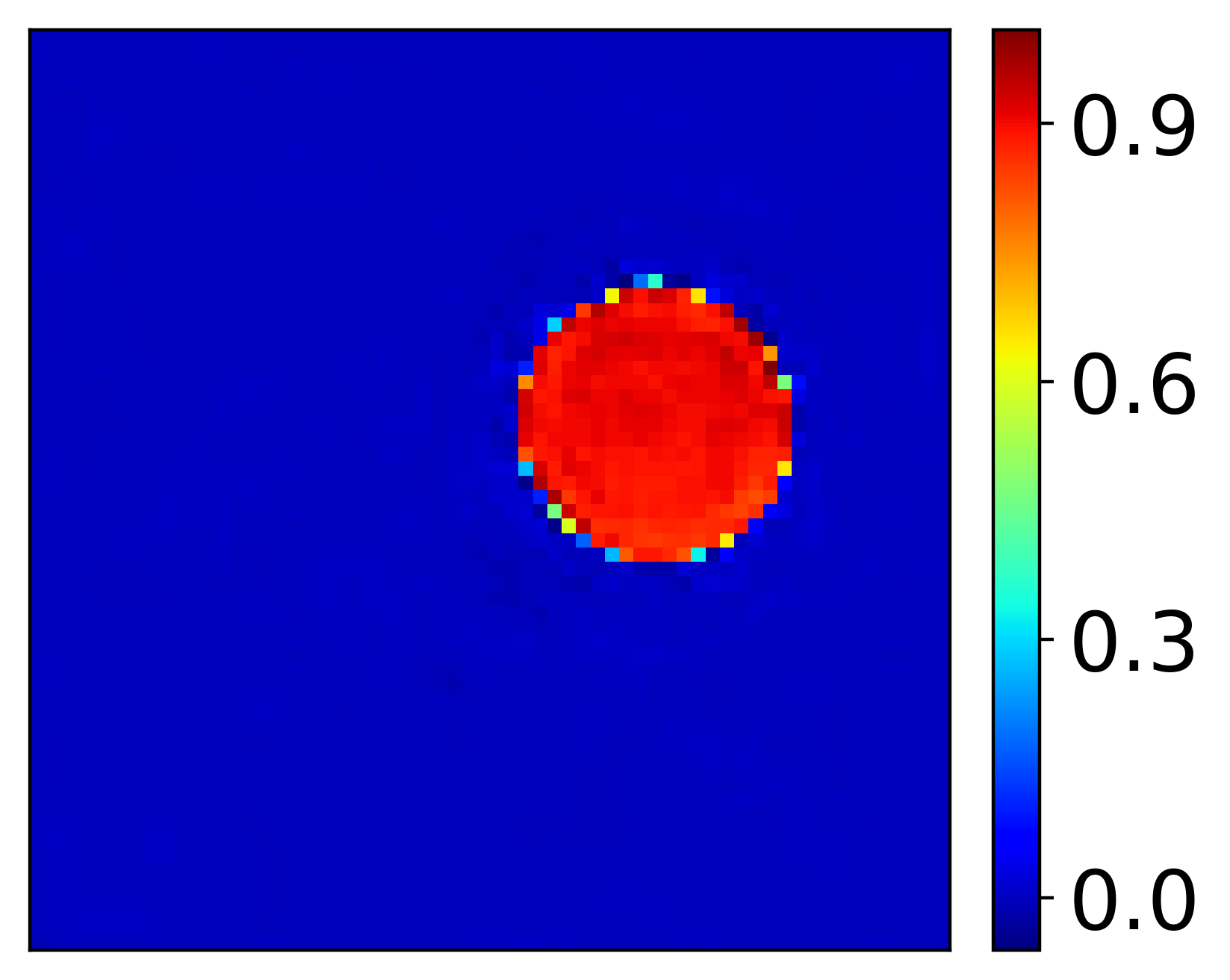}}
	}
	\hspace{-0.4cm}
	\subcaptionbox{}{
		\adjustbox{valign=t}{\includegraphics[width=0.18\linewidth]{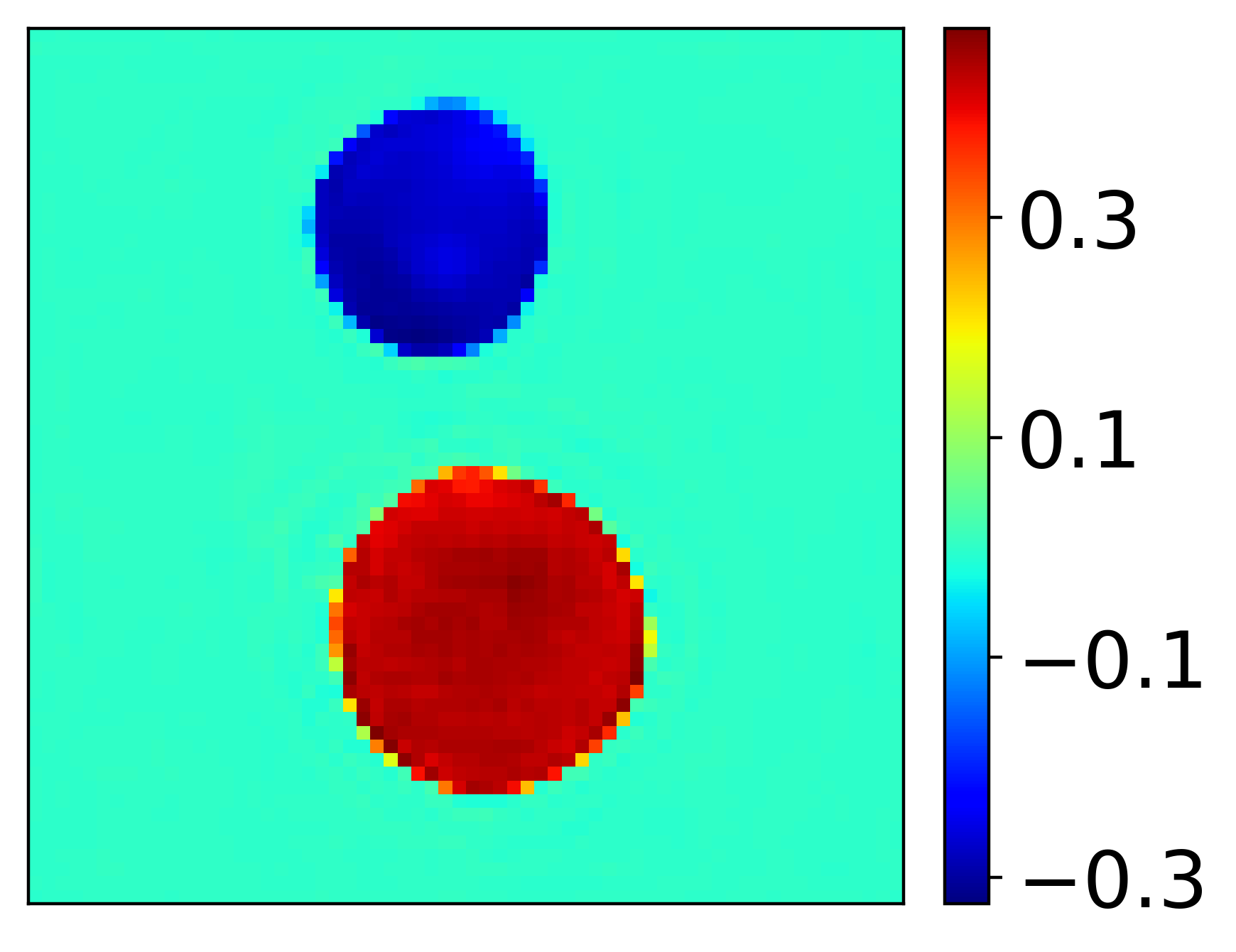}}
	}
	\hspace{-0.4cm}
	\subcaptionbox{}{
		\adjustbox{valign=t}{\includegraphics[width=0.18\linewidth]{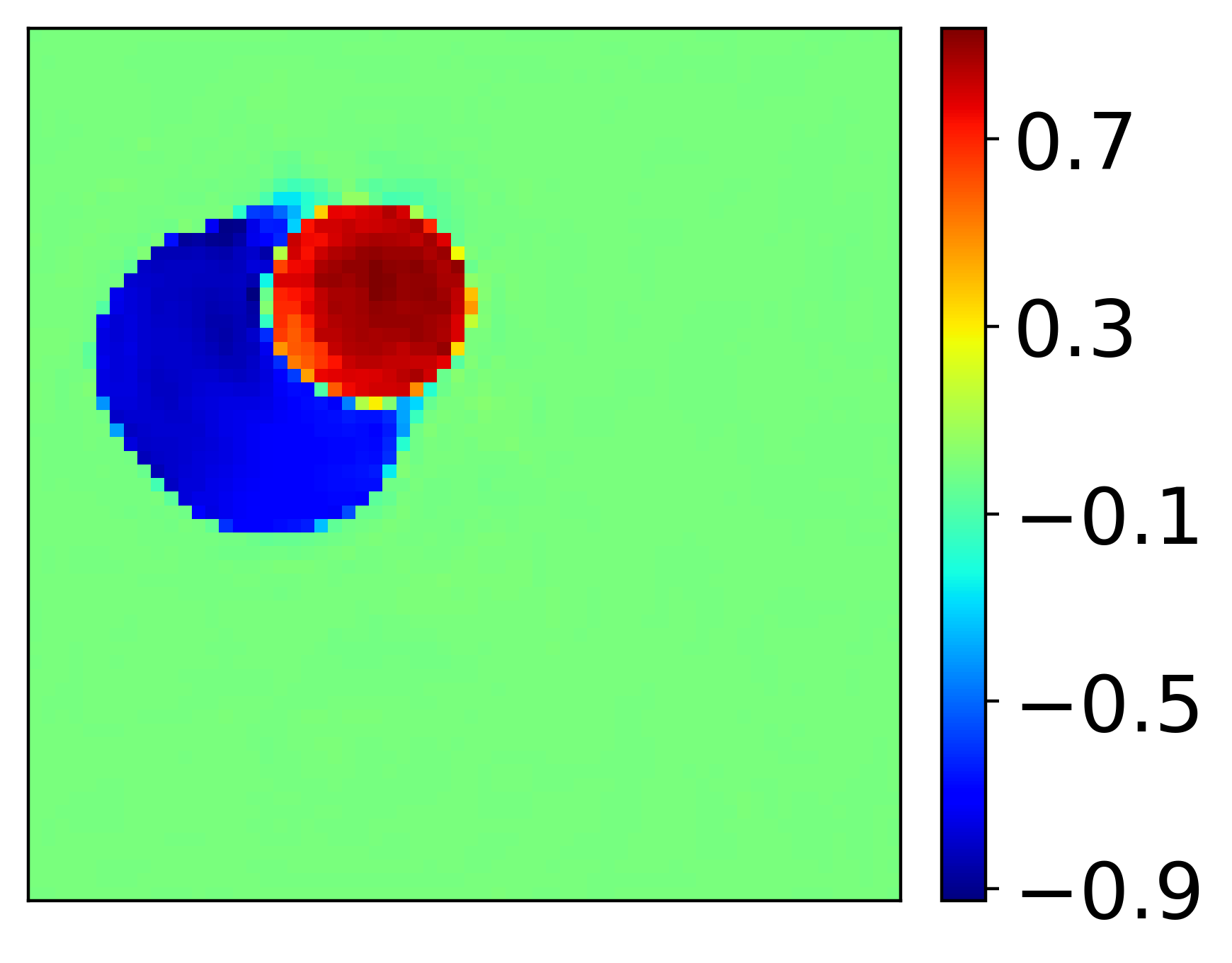}}
	}
	\hspace{-0.4cm}
	\subcaptionbox{}{
		\adjustbox{valign=t}{\includegraphics[width=0.18\linewidth]{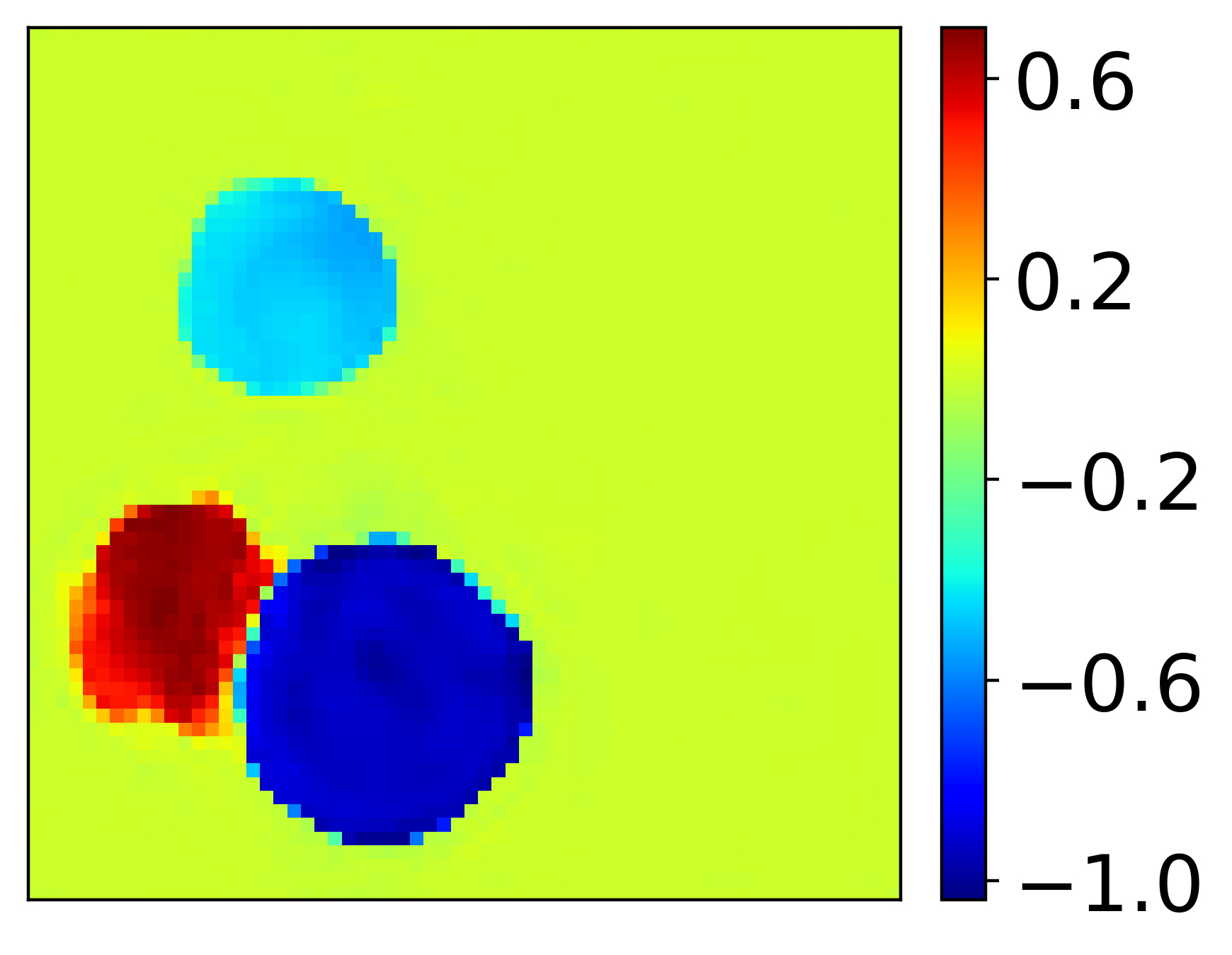}}
	}
	\hspace{-0.4cm}
	\subcaptionbox{}{
		\adjustbox{valign=t}{\includegraphics[width=0.18\linewidth]{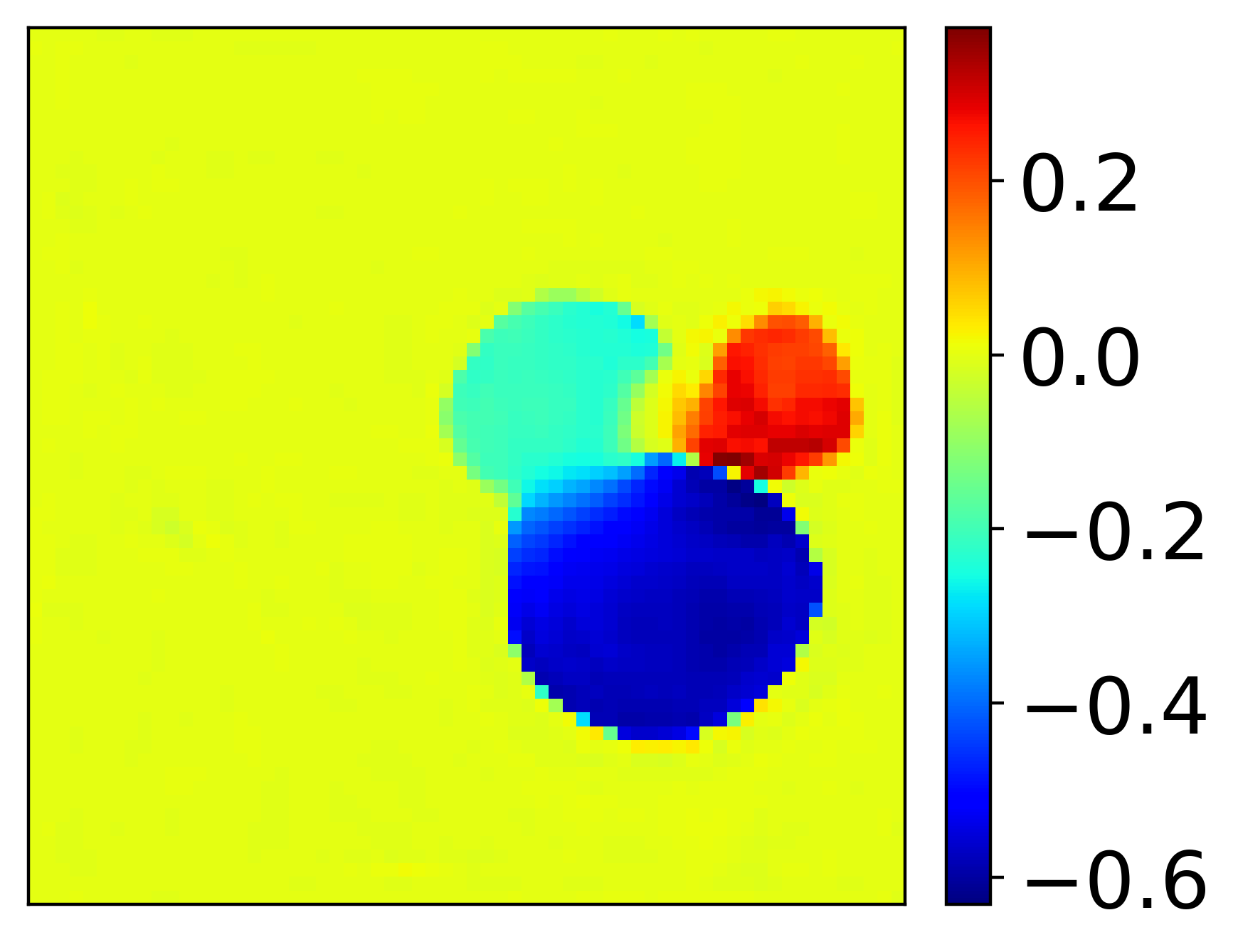}}
	}
	\caption{Fine-tuned reconstruction results using the high-to-low noise transfer learning strategy. Rows from top to bottom: ground-truth sources, $5\%$ noise reconstructions, and $50\%$ noise reconstructions.}
	\label{fig:tl_net}
\end{figure}
\begin{table}[htbp]
	\centering
	\begin{tabular}{l|ccc}
		\hline
		Method            & Noise level & NMSE     & Epochs \\
		\hline
		From scratch      & $5\%$       & $5.35\%$ & $200$  \\
		Transfer learning & $5\%$       & $5.62\%$ & $30$   \\
		\hline
		From scratch      & $50\%$      & $8.48\%$ & $200$  \\
		Transfer learning & $50\%$      & $7.44\%$ & $30$   \\
		\hline
	\end{tabular}
	\caption{Comparison of reconstruction performance and training efficiency between \enquote{from scratch} training and transfer learning initialized with the $100\%$ noise model.}
	\label{tab:tl}
\end{table}
The resulting reconstructions and associated errors are summarized in \Cref{fig:tl_net} and \Cref{tab:tl}. The transfer-learning-enhanced models achieve fidelity comparable to those trained from scratch. Importantly, the pre-trained initialization enables the network to converge in significantly fewer iterations, reducing both training time and computational overhead. These results confirm that training on high-noise data captures the underlying geometric distribution effectively, providing a superior starting point for refinement in cleaner data regimes.

\subsection{Complex Geometric Structures: MNIST Benchmark Dataset}
This experiment evaluates the method's robustness using the MNIST dataset of handwritten digits as a geometric benchmark, assessing performance under conditions of extreme sparsity and noise. Original $28 \times 28$ images are resized to $64 \times 64$ and normalized to $[0, 1]$. To suppress background noise, values below $0.1$ are thresholded to zero. We set a minimal truncation frequency of $N = 2$ and generate $5000$ samples ($4500$ for training, $500$ for testing). The model is trained for $50$ epochs with an initial learning rate of $0.001$, decaying by a factor of $0.5$ every $5$ epochs. To quantify structural fidelity beyond pixel-wise error, we employ the structural similarity index measure (SSIM) \cite{wang2003multiscale}:
\begin{equation*}
	\text{SSIM}(X, Y) = \frac{(2\mu_X\mu_Y + C_1)(2\sigma_{XY} + C_2)}{(\mu_X^2 + \mu_Y^2 +C_1)(\sigma_X^2 + \sigma_Y^2 + C_2)},
\end{equation*}
where $\mu$ and $\sigma$ denote the mean and standard deviation, respectively, and $\sigma_{XY}$ represents the covariance. Constants $C_1 = (K_1L_X)^2$ and $C_2 = (K_2L_X)^2$ ensure numerical stability near zero, with parameters fixed at $L_X = 1$, $K_1 = 0.01$, and $K_2 = 0.03$.

% ... (Figures omitted for brevity in diff) ...
\begin{figure}[htbp]
	\centering
	\subcaptionbox{}{
		\includegraphics[width=0.19\linewidth]{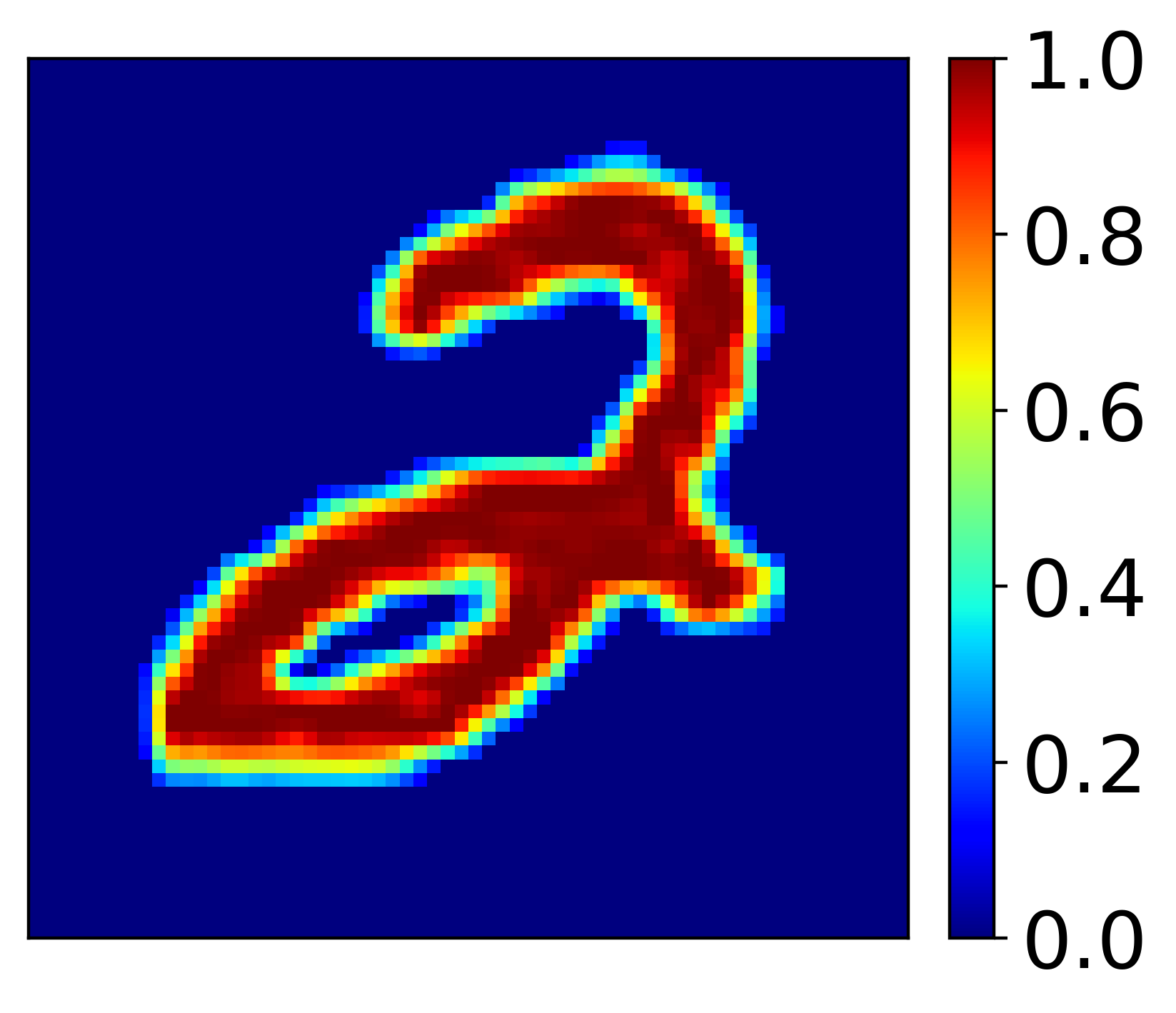}
	}
	\hspace{-0.4cm}
	\subcaptionbox{}{
		\includegraphics[width=0.19\linewidth]{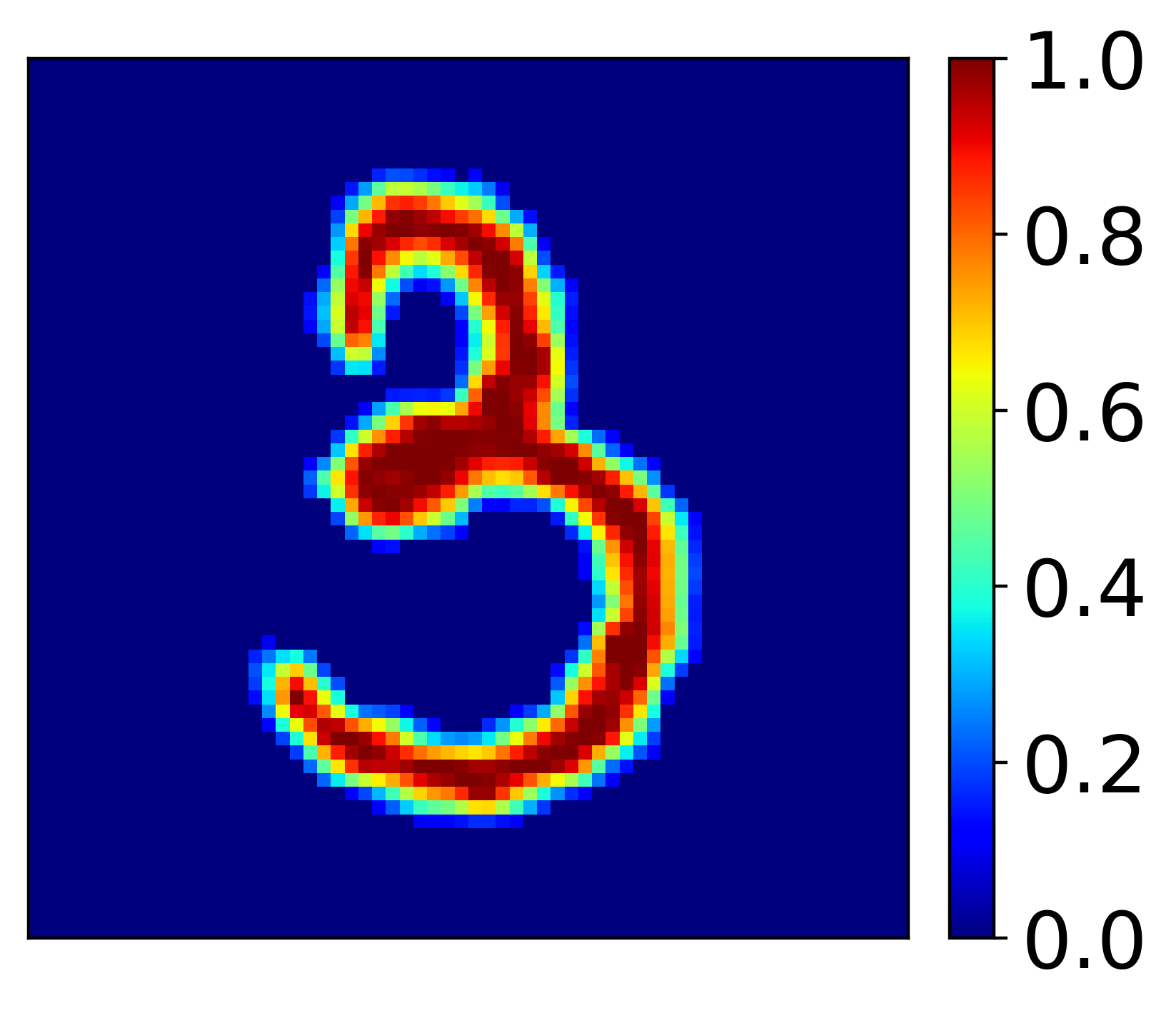}
	}
	\hspace{-0.4cm}
	\subcaptionbox{}{
		\includegraphics[width=0.19\linewidth]{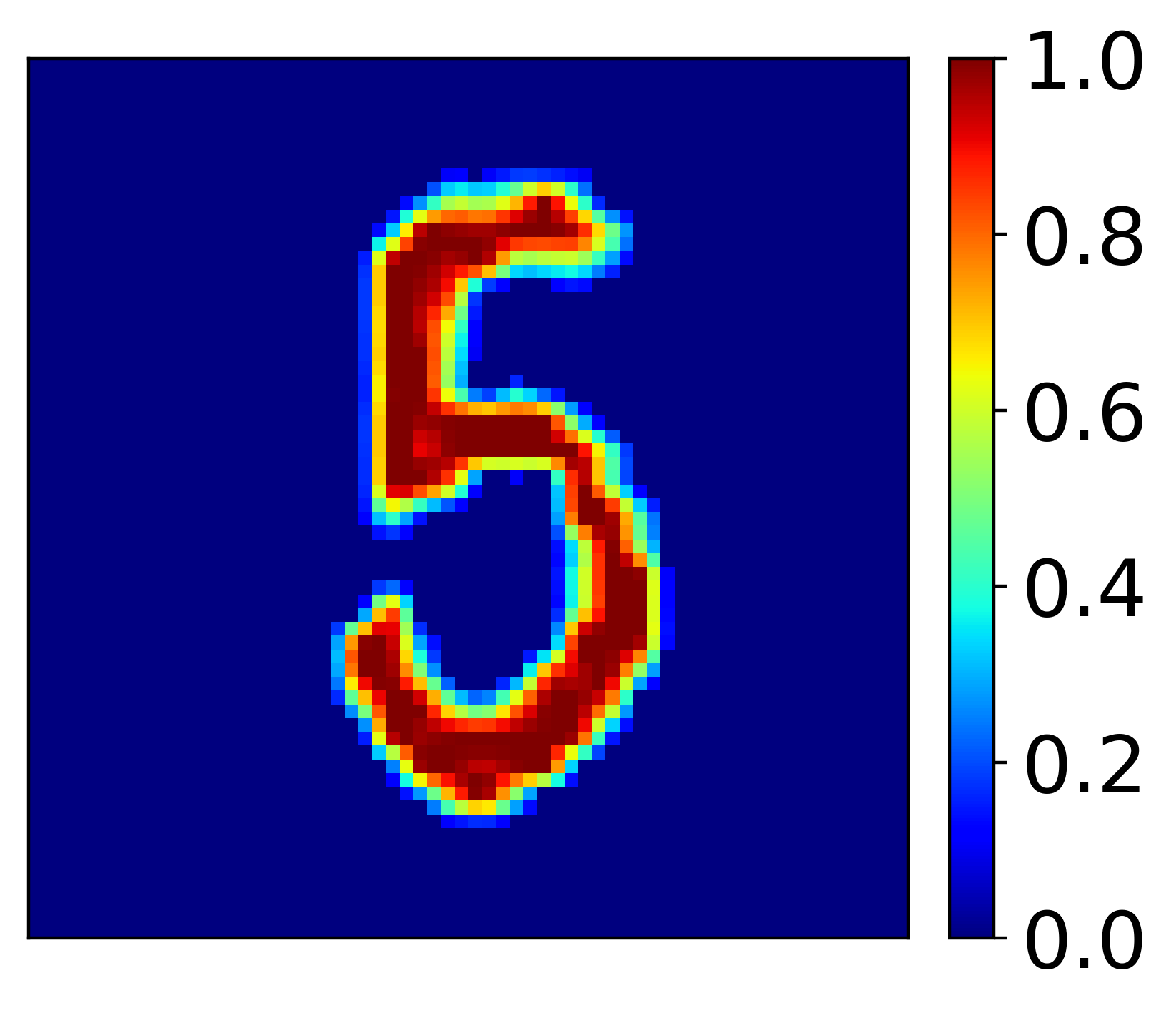}
	}
	\hspace{-0.4cm}
	\subcaptionbox{}{
		\includegraphics[width=0.19\linewidth]{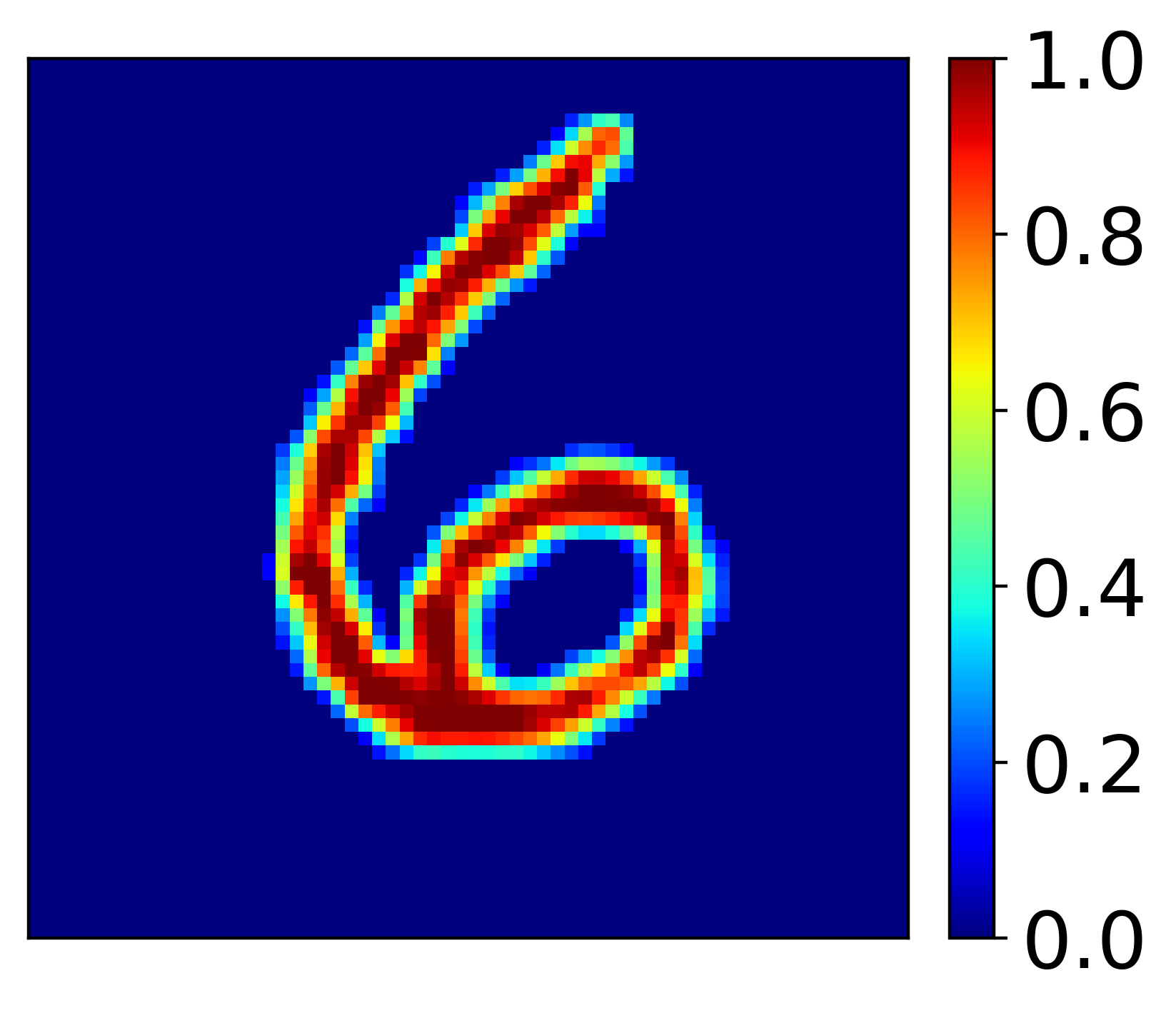}
	}
	\hspace{-0.4cm}
	\subcaptionbox{}{
		\includegraphics[width=0.19\linewidth]{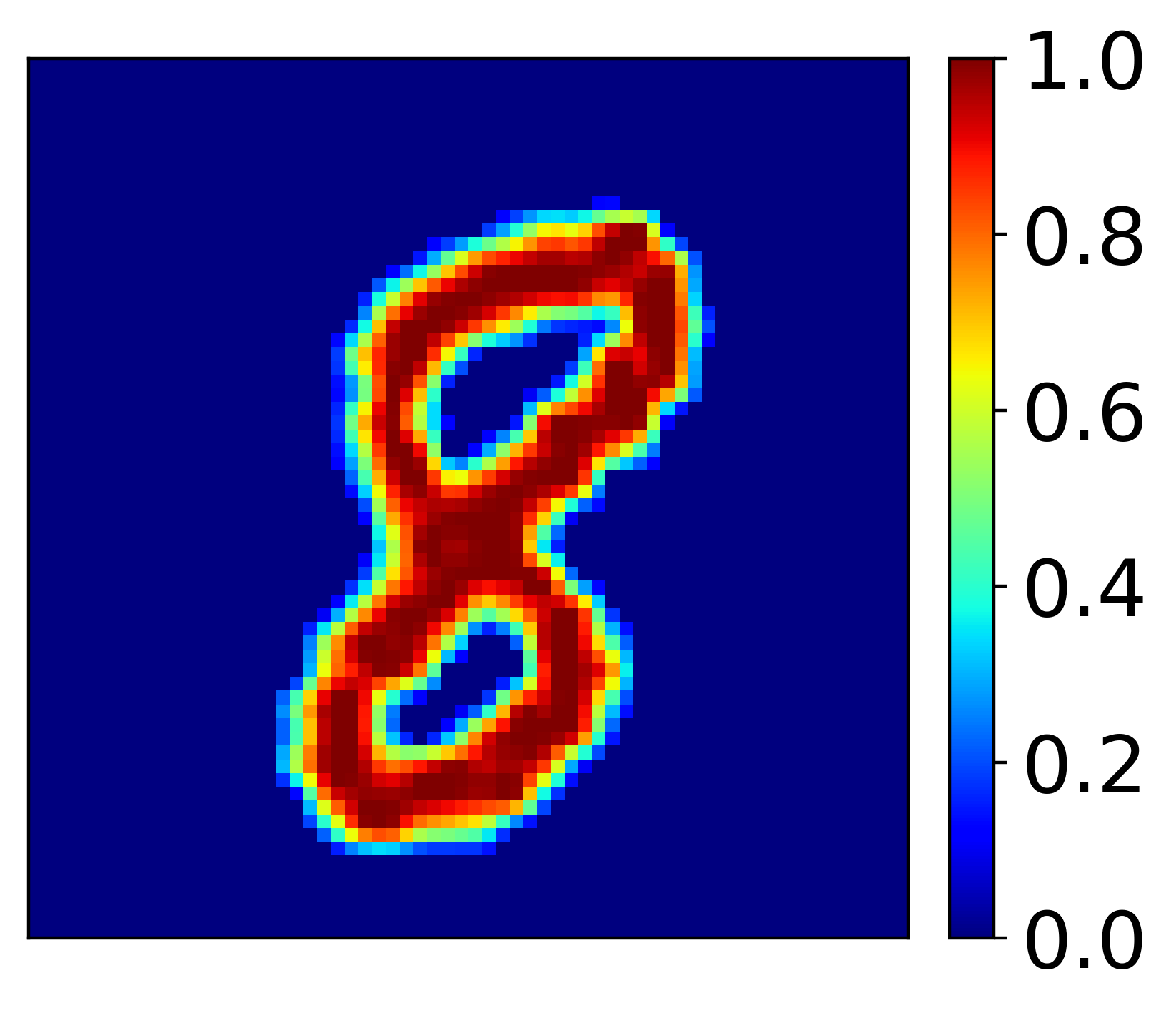}
	}\\
	\subcaptionbox{}{
		\includegraphics[width=0.19\linewidth]{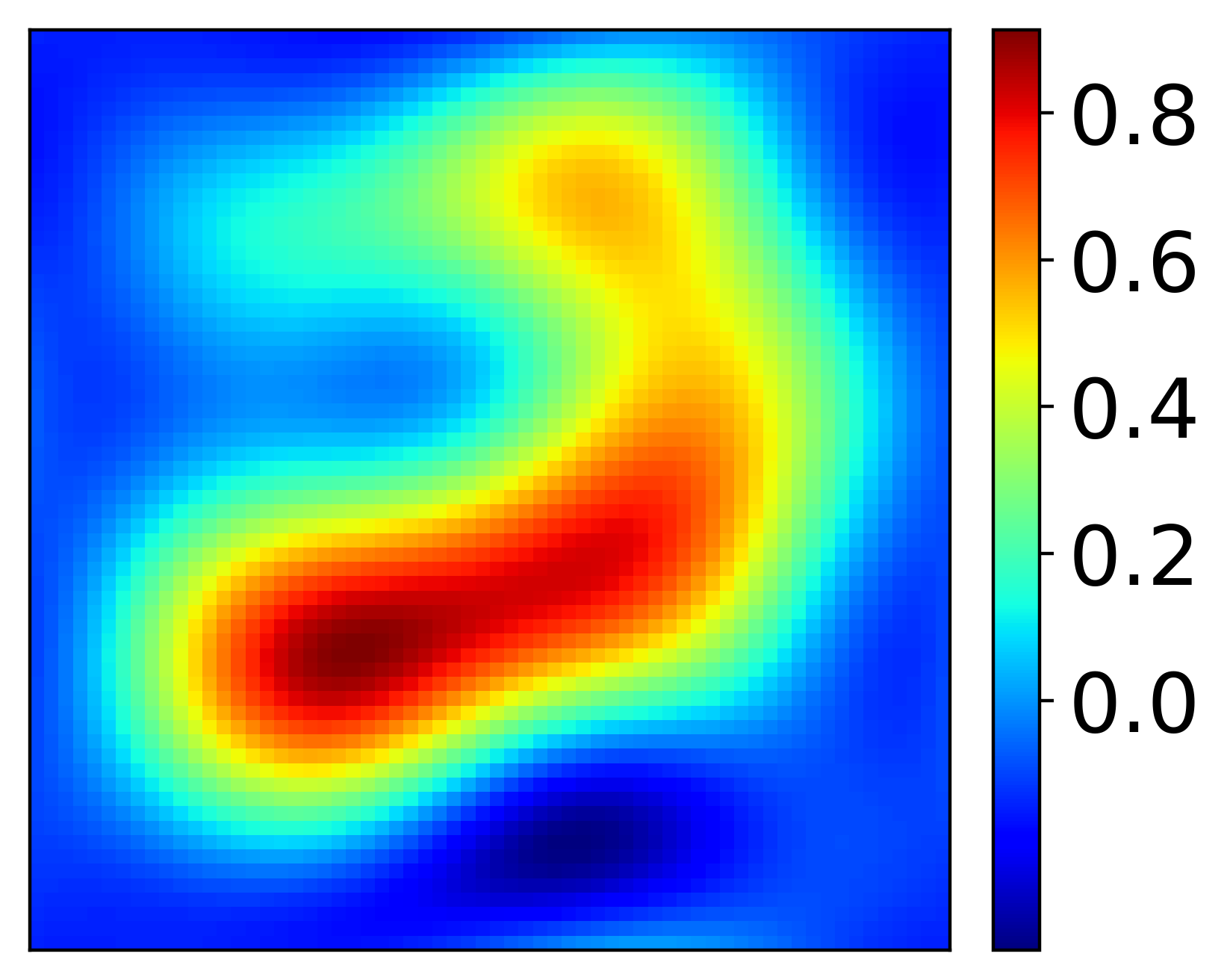}
	}
	\hspace{-0.4cm}
	\subcaptionbox{}{
		\includegraphics[width=0.19\linewidth]{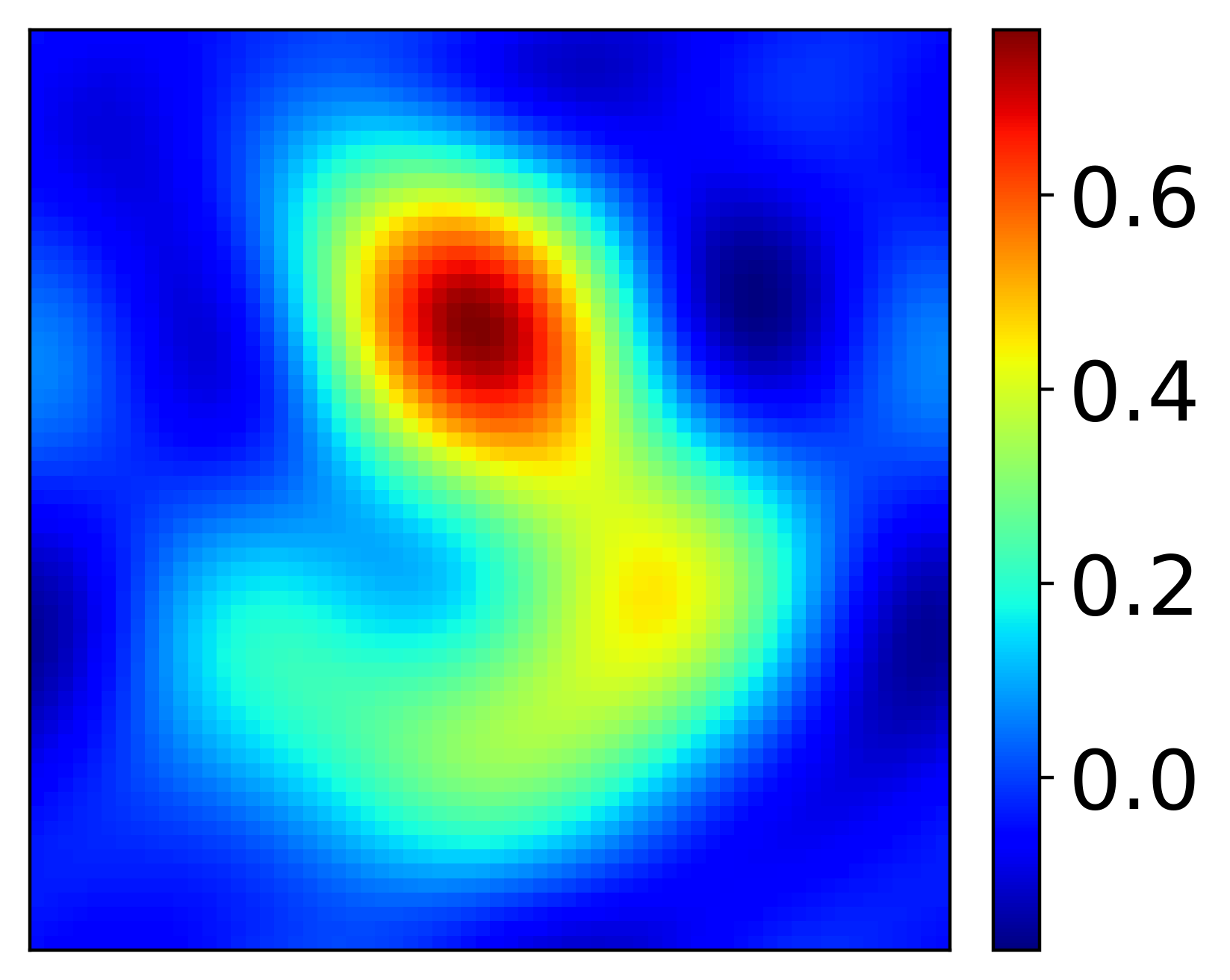}
	}
	\hspace{-0.4cm}
	\subcaptionbox{}{
		\includegraphics[width=0.19\linewidth]{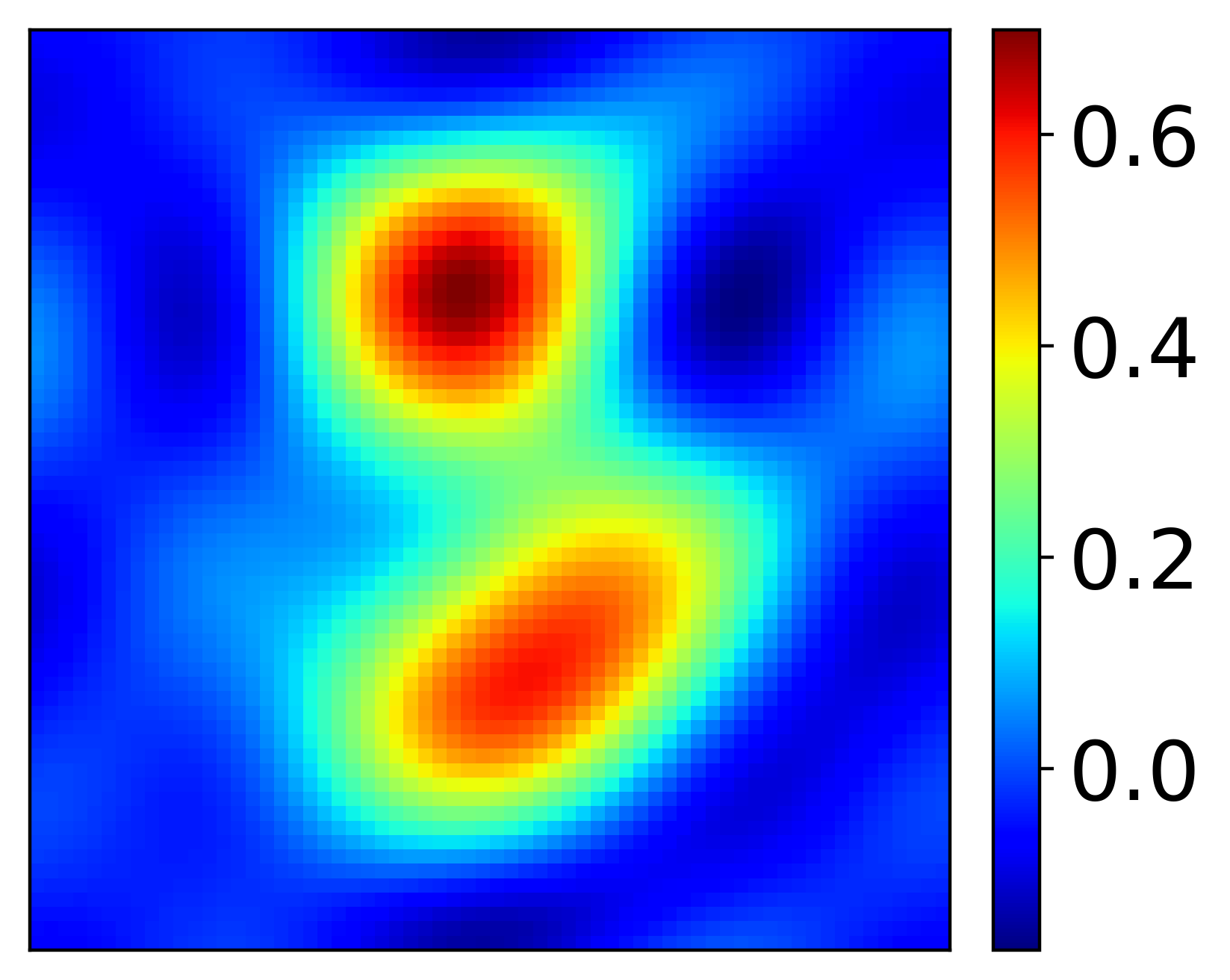}
	}
	\hspace{-0.4cm}
	\subcaptionbox{}{
		\includegraphics[width=0.19\linewidth]{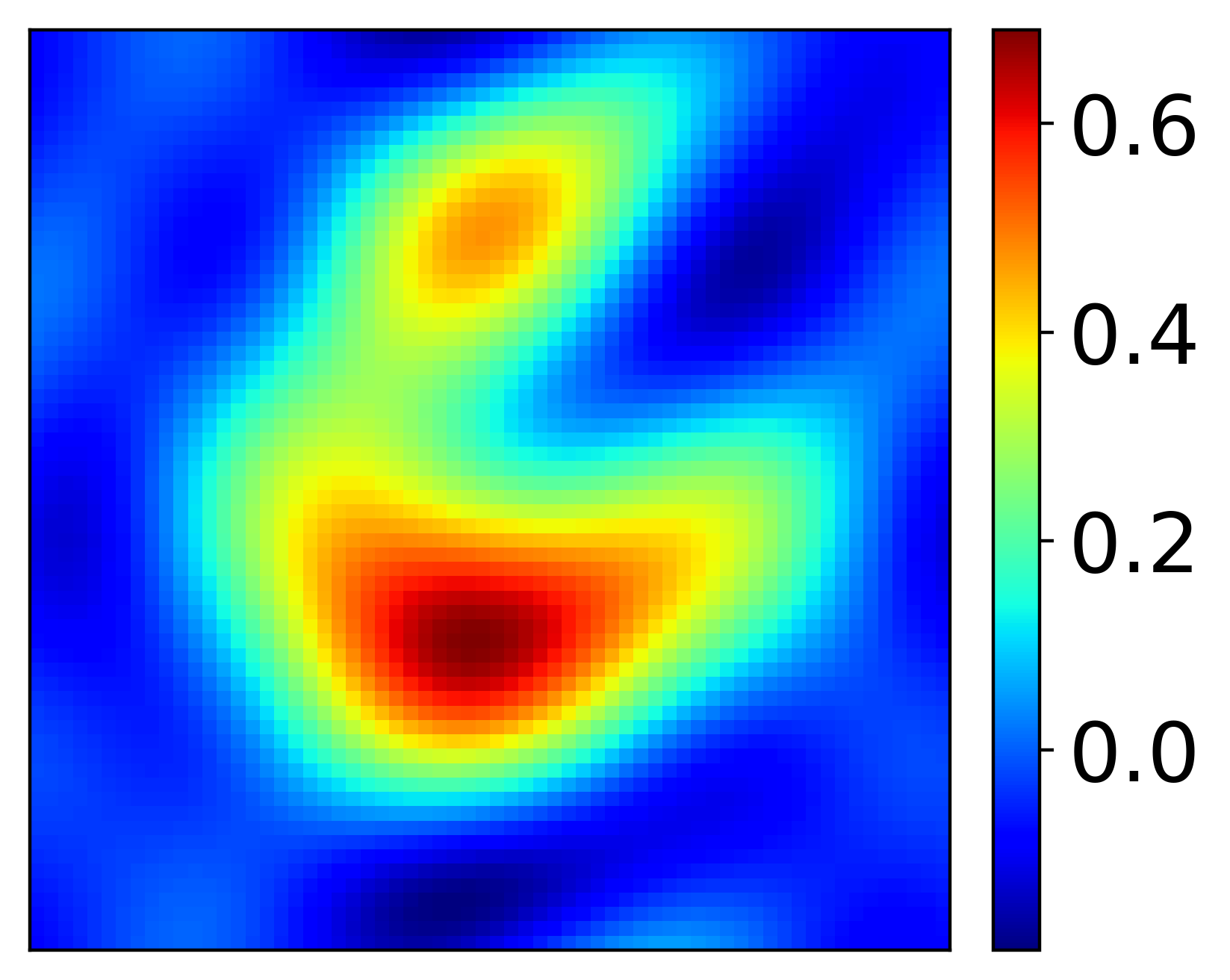}
	}
	\hspace{-0.4cm}
	\subcaptionbox{}{
		\includegraphics[width=0.19\linewidth]{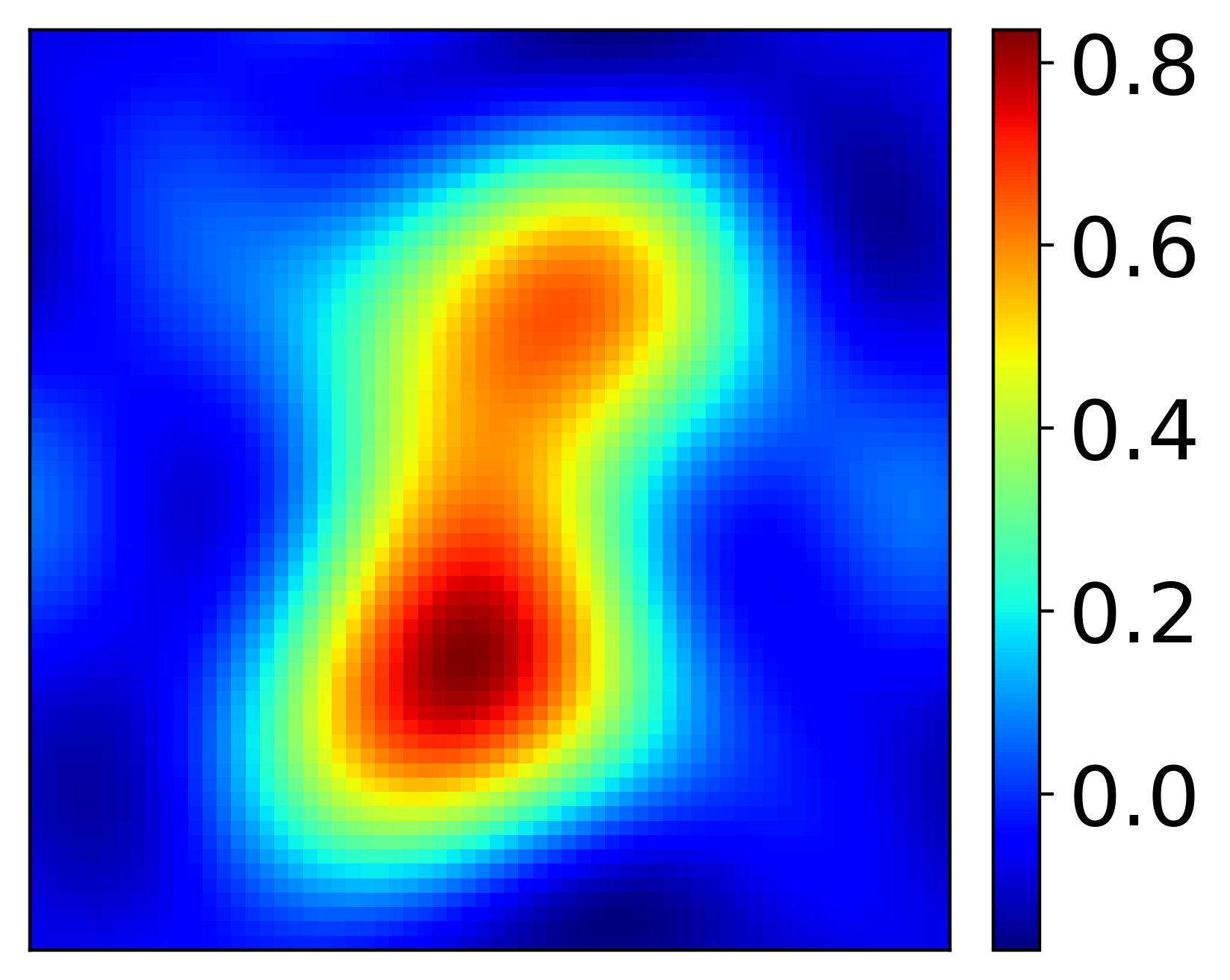}
	}\\
	\subcaptionbox{}{
		\includegraphics[width=0.19\linewidth]{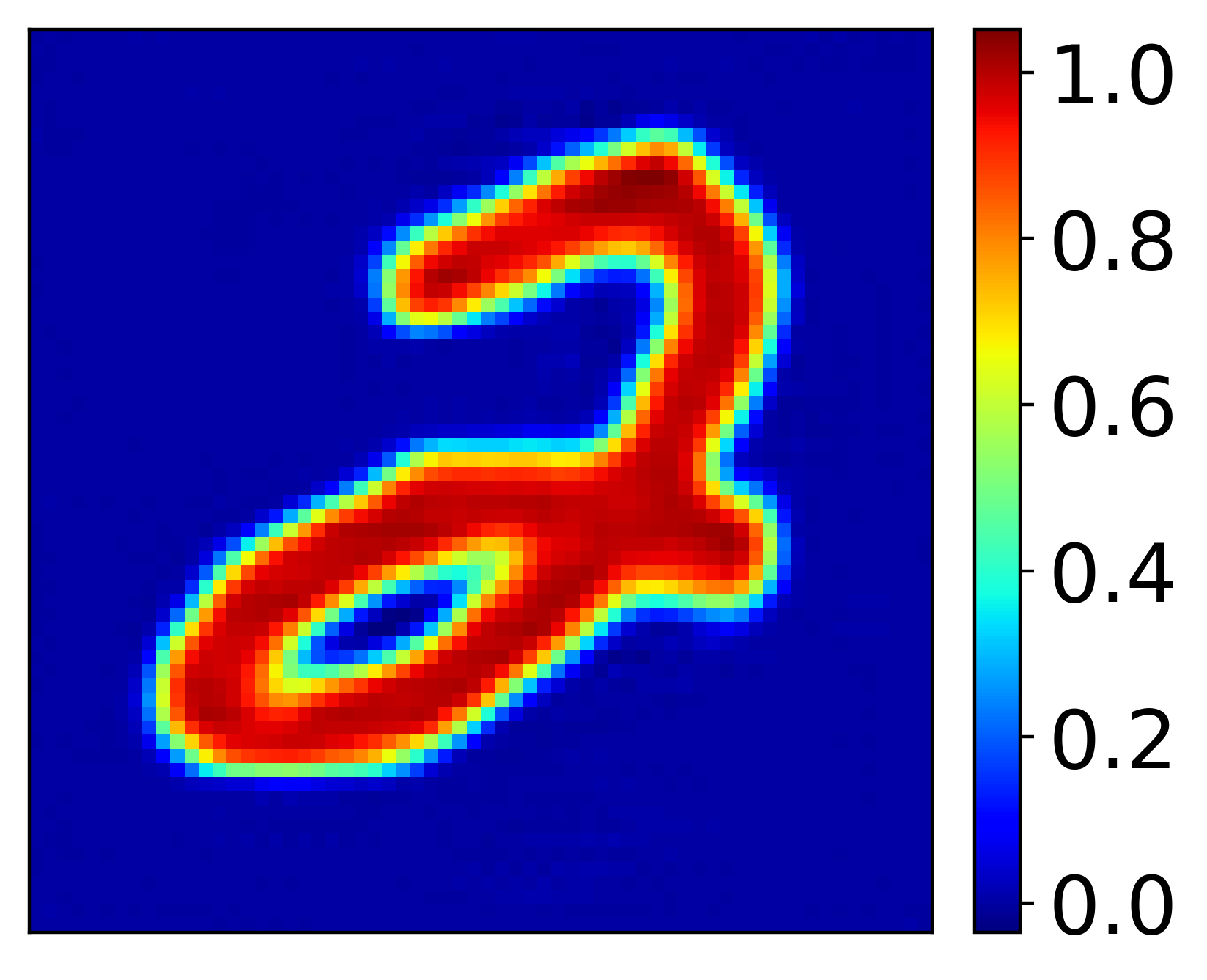}
	}
	\hspace{-0.4cm}
	\subcaptionbox{}{
		\includegraphics[width=0.19\linewidth]{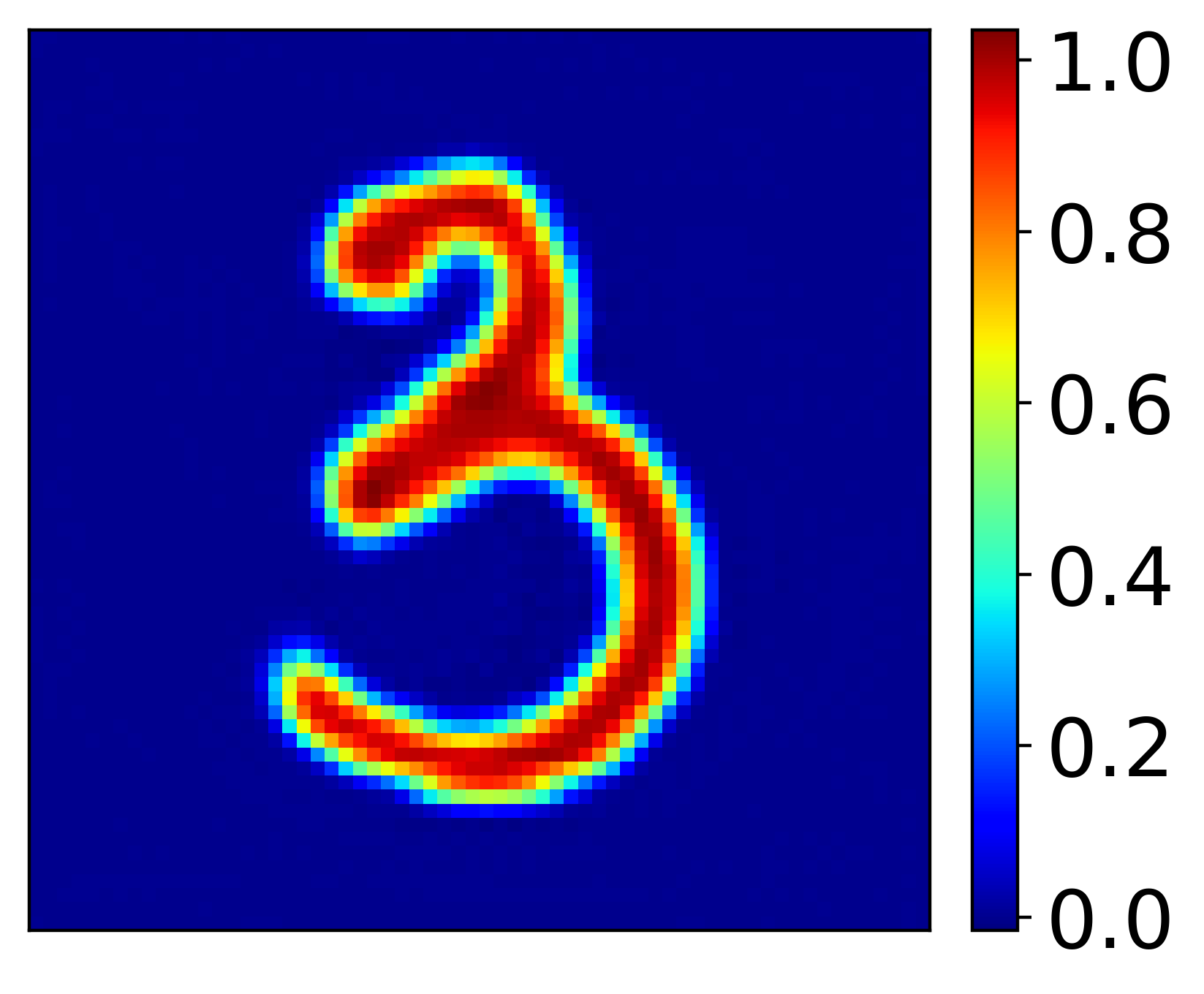}
	}
	\hspace{-0.4cm}
	\subcaptionbox{}{
		\includegraphics[width=0.19\linewidth]{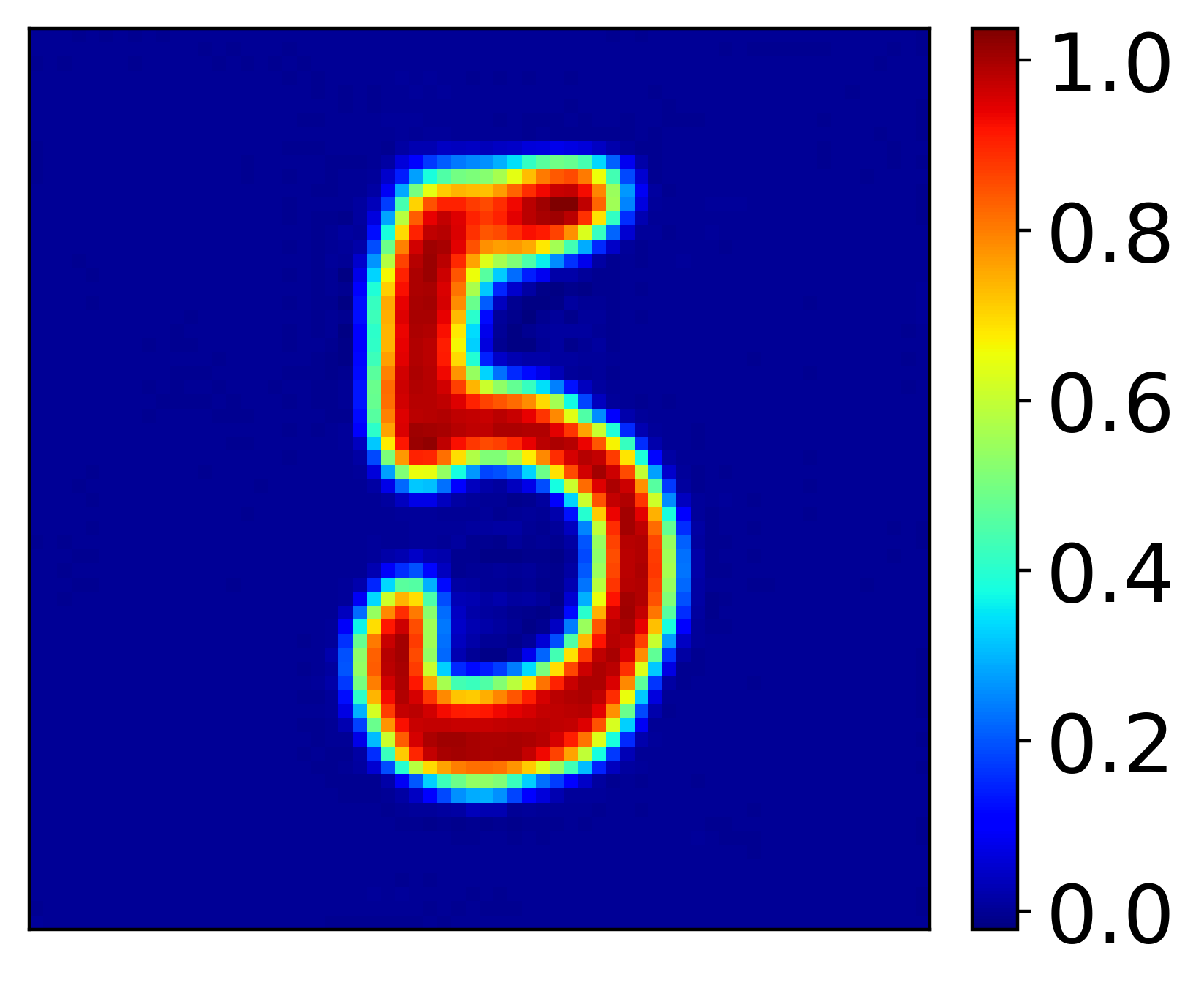}
	}
	\hspace{-0.4cm}
	\subcaptionbox{}{
		\includegraphics[width=0.19\linewidth]{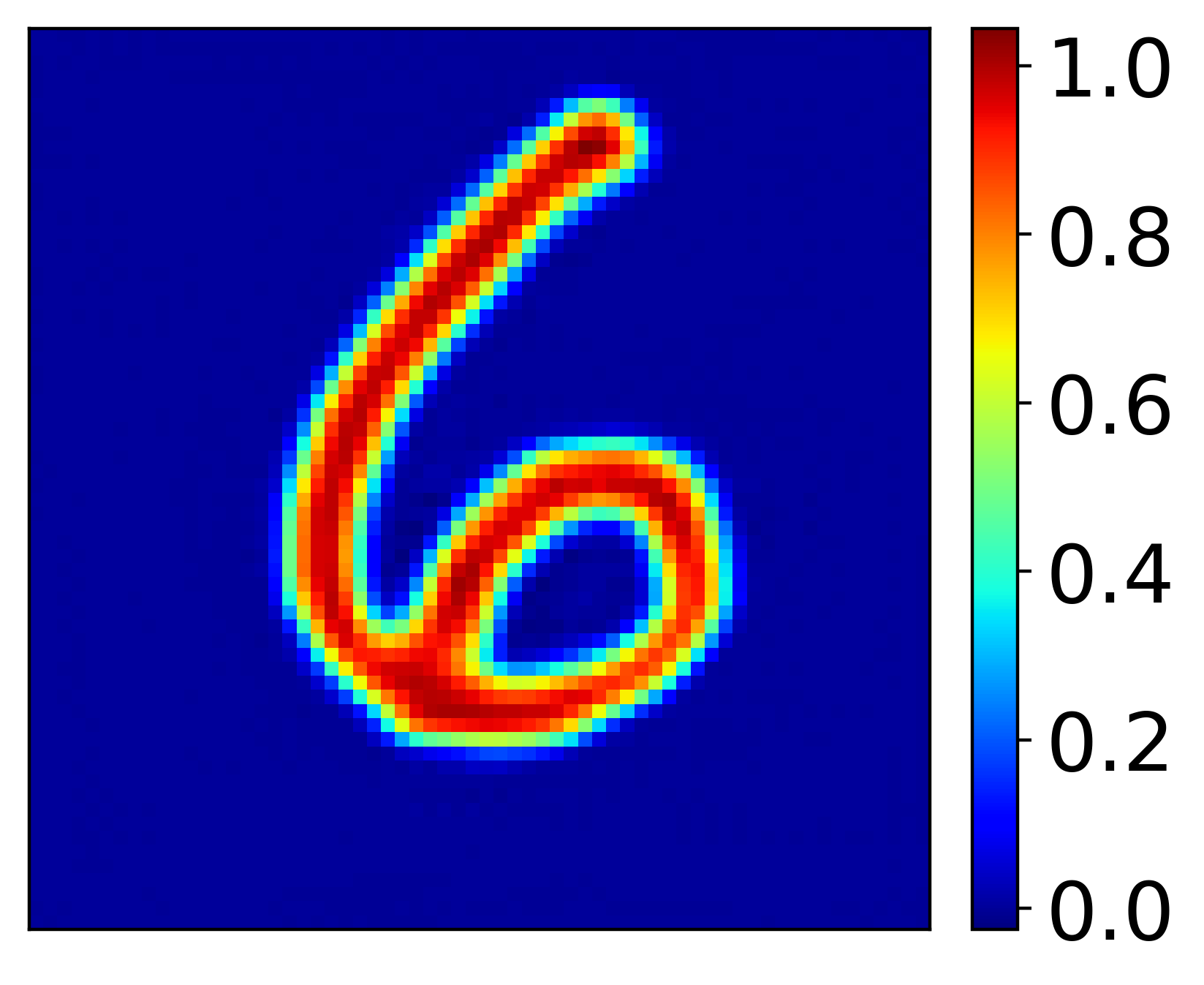}
	}
	\hspace{-0.4cm}
	\subcaptionbox{}{
		\includegraphics[width=0.19\linewidth]{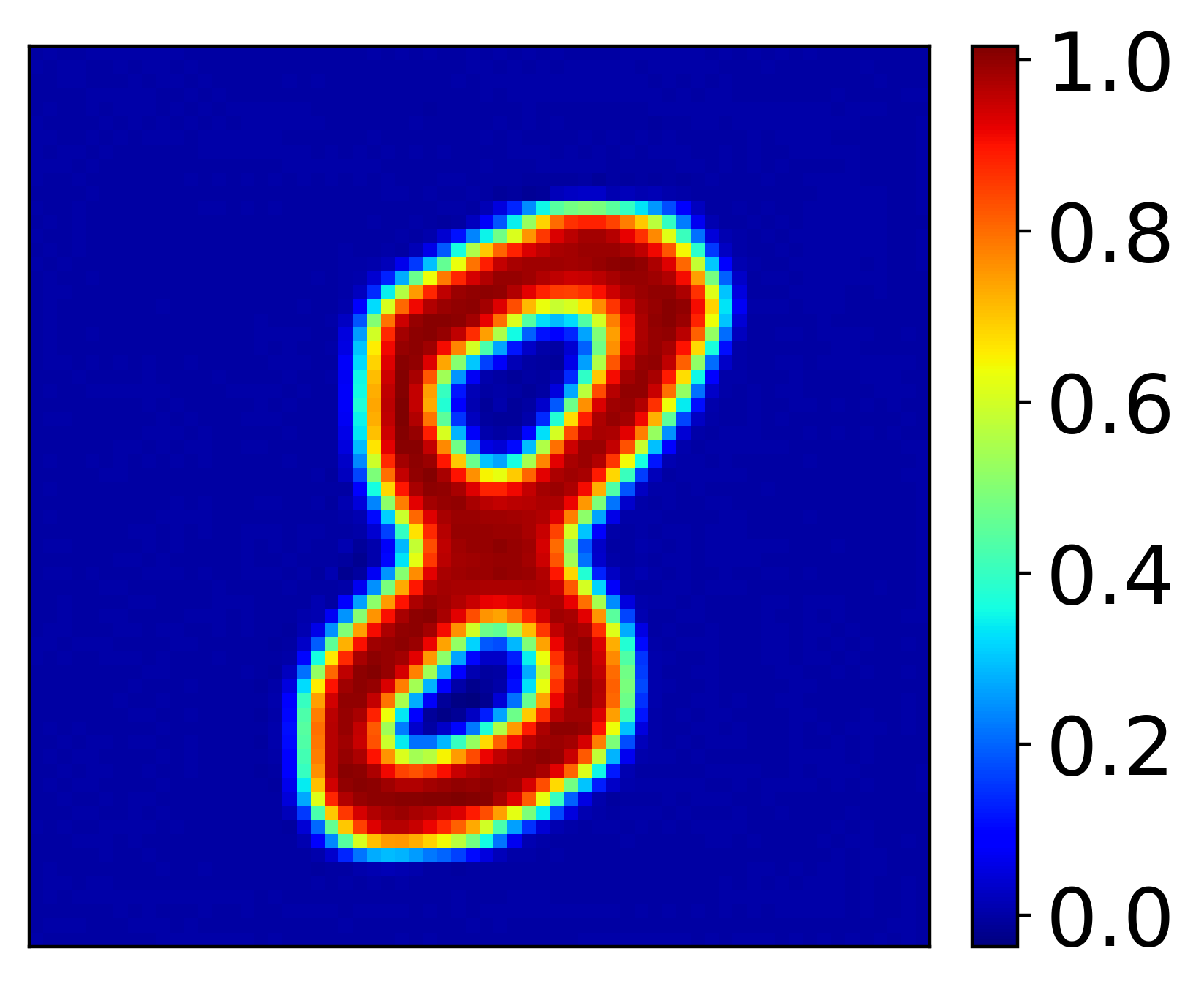}
	}\\
	\caption{Reconstruction results for MNIST handwritten digits under $50\%$ noise. Rows from top to bottom: ground-truth sources, Fourier reconstructions ($N=2$, which serve as the input for the U-Net in the third row), and U-Net-enhanced reconstructions.} 
	\label{fig:mnist_5}
\end{figure}

\begin{figure}[htbp]
	\centering
	\subcaptionbox{}{
		\includegraphics[width=0.24\linewidth]{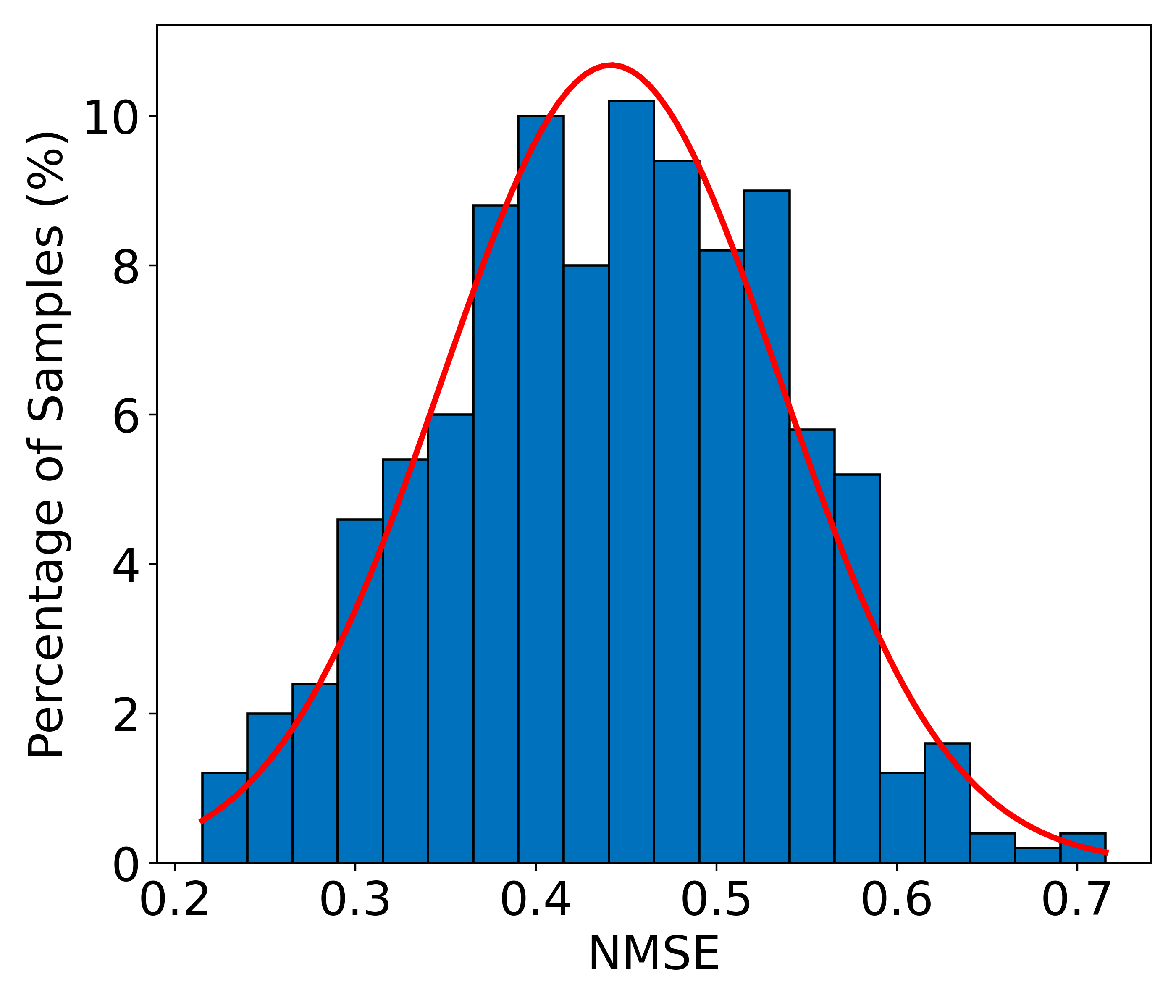}
	}
	\hspace{-0.4cm}
	\subcaptionbox{}{
		\includegraphics[width=0.24\linewidth]{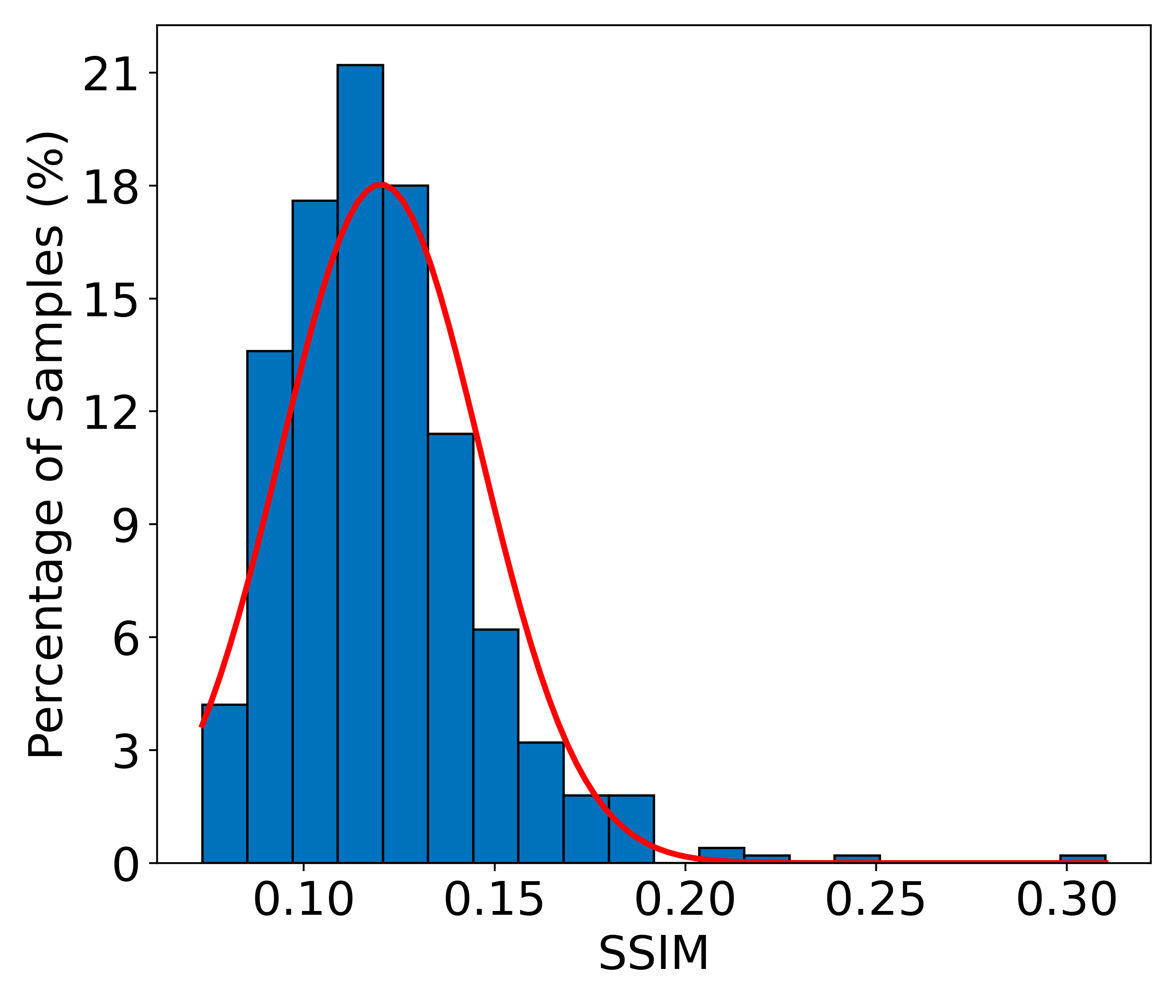}
	}
	\hspace{-0.4cm}
	\subcaptionbox{}{
		\includegraphics[width=0.24\linewidth]{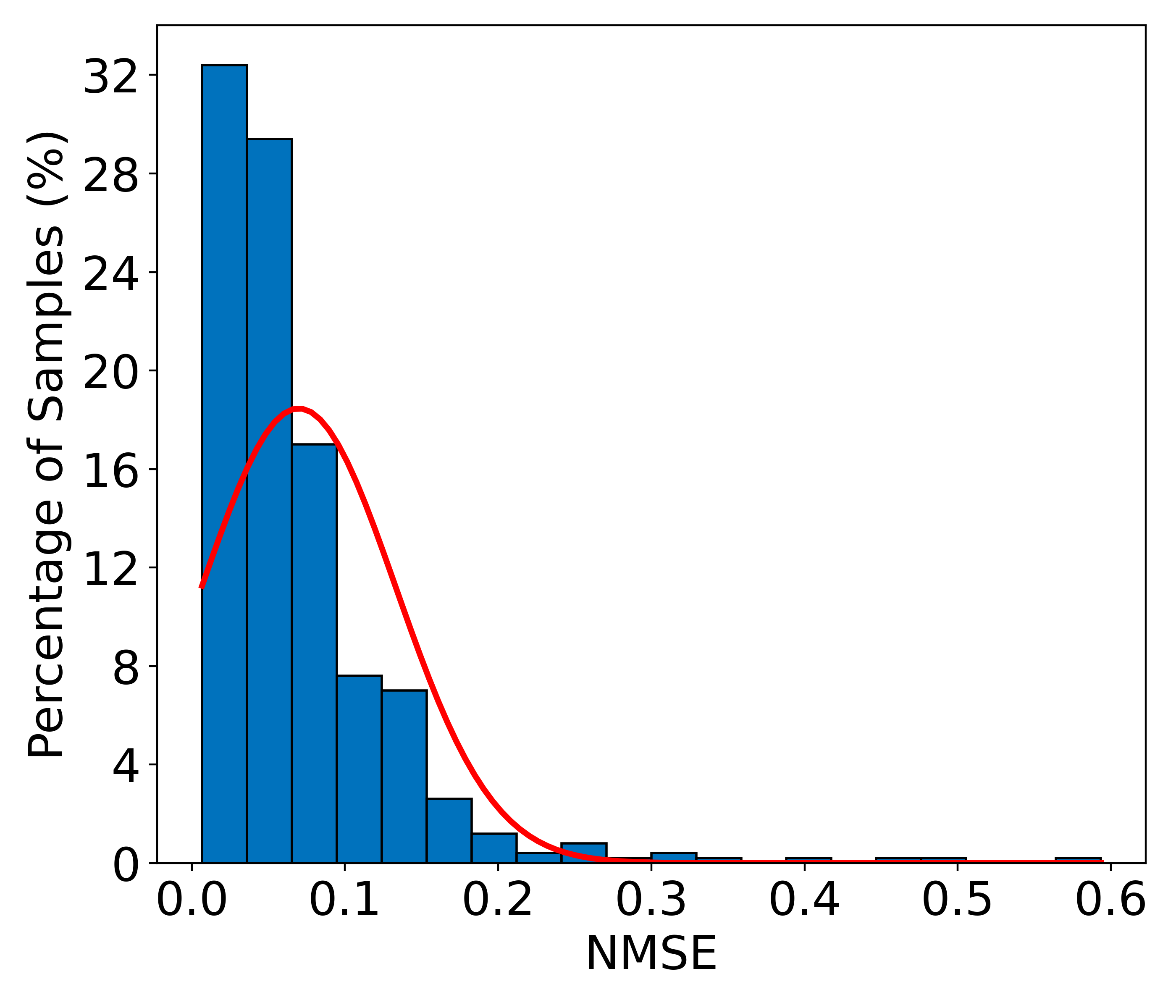}
	}
	\hspace{-0.4cm}
	\subcaptionbox{}{
		\includegraphics[width=0.24\linewidth]{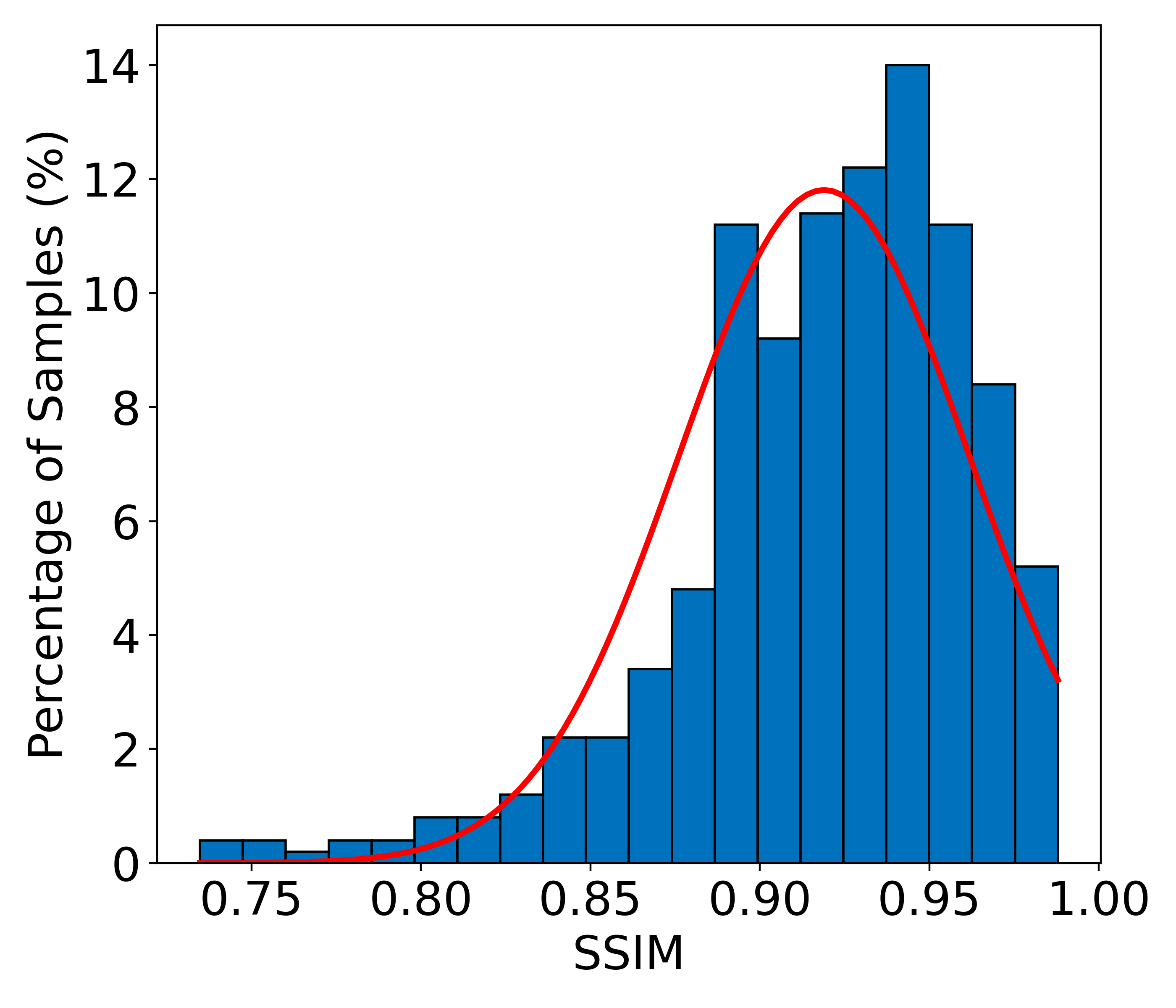}
	}
	\caption{Statistical distribution of NMSE and SSIM for the MNIST test set under $50\%$ noise. The first two subfigures represent the Fourier reconstruction ($N=2$), and the last two represent the U-Net-enhanced results. The histograms are fitted with normal density functions.}
	\label{fig:nmse_mnist}
\end{figure}

As shown in \Cref{fig:mnist_5} and \Cref{fig:nmse_mnist}, the proposed method successfully recovers complex handwritten structures even under high noise and minimal measurements ($N=2$). The statistical results in \Cref{fig:nmse_mnist} provide a direct comparison between the initial Fourier reconstruction and the U-Net-enhanced output. Although the initial Fourier reconstruction—serving as the network's input—is highly blurred and dominated by artifacts, the U-Net-enhanced model dramatically improves fidelity under $50\%$ noise, achieving an average NMSE of $0.07$ and an SSIM of $0.92$. These results demonstrate stable performance under strong interference and poor initial physics-based estimates.

\subsection{Out-of-Distribution Generalization: Letter Dataset}

% ... (Figures omitted for brevity in diff) ...
\begin{figure}[htbp]
	\centering
	\subcaptionbox{}{
		\includegraphics[width=0.19\linewidth]{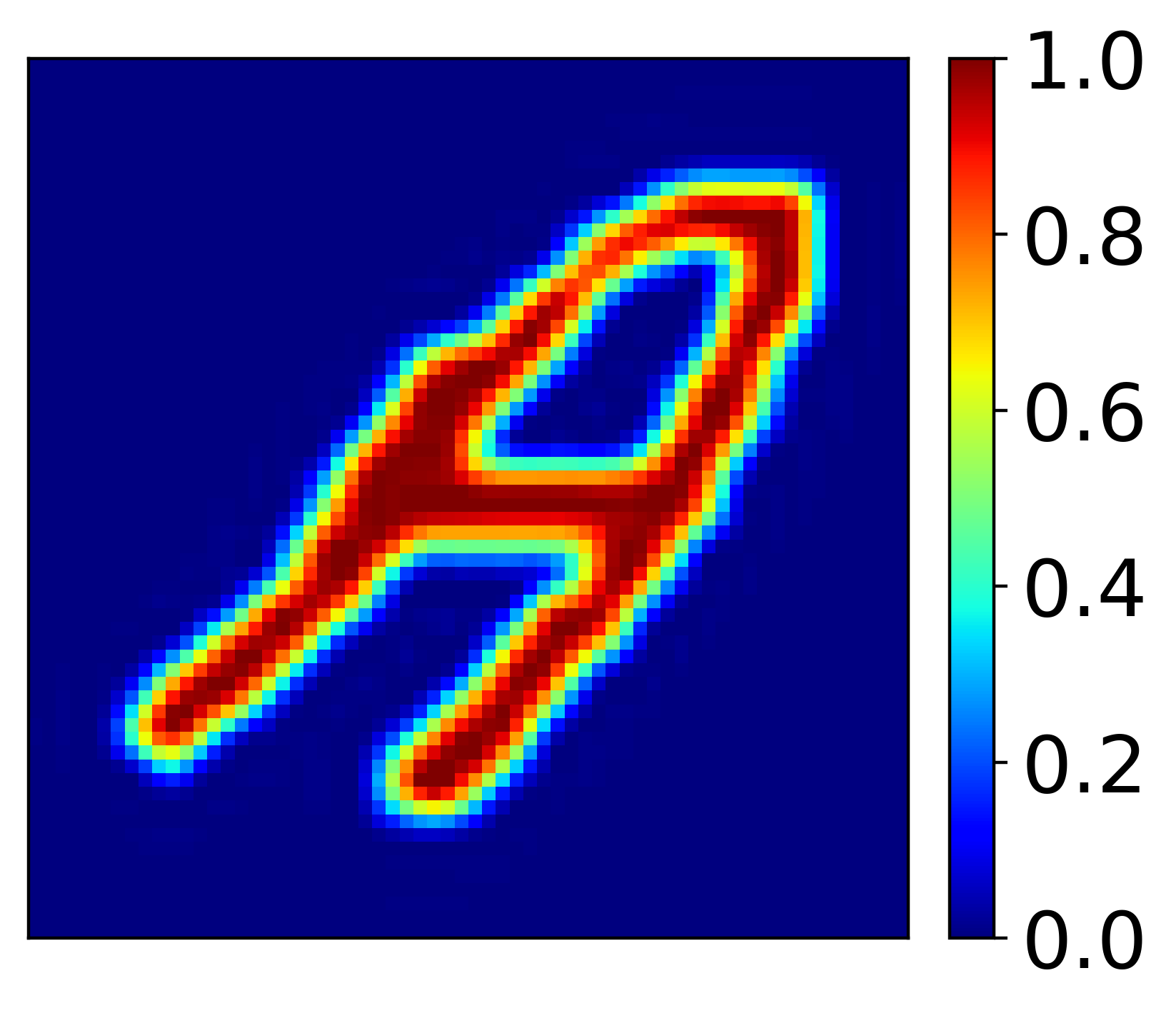}
	}
	\hspace{-0.4cm}
	\subcaptionbox{}{
		\includegraphics[width=0.19\linewidth]{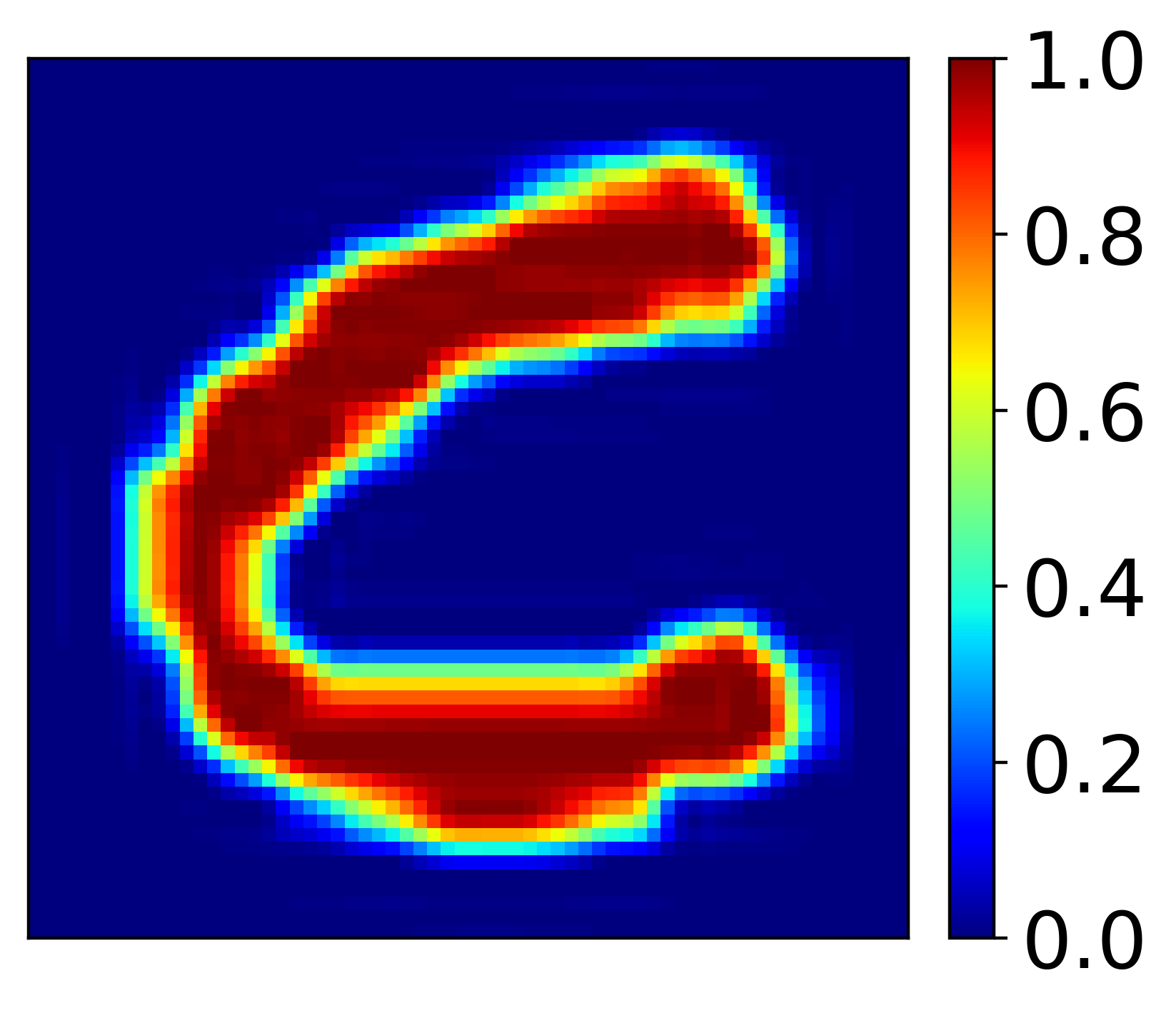}
	}
	\hspace{-0.4cm}
	\subcaptionbox{}{
		\includegraphics[width=0.19\linewidth]{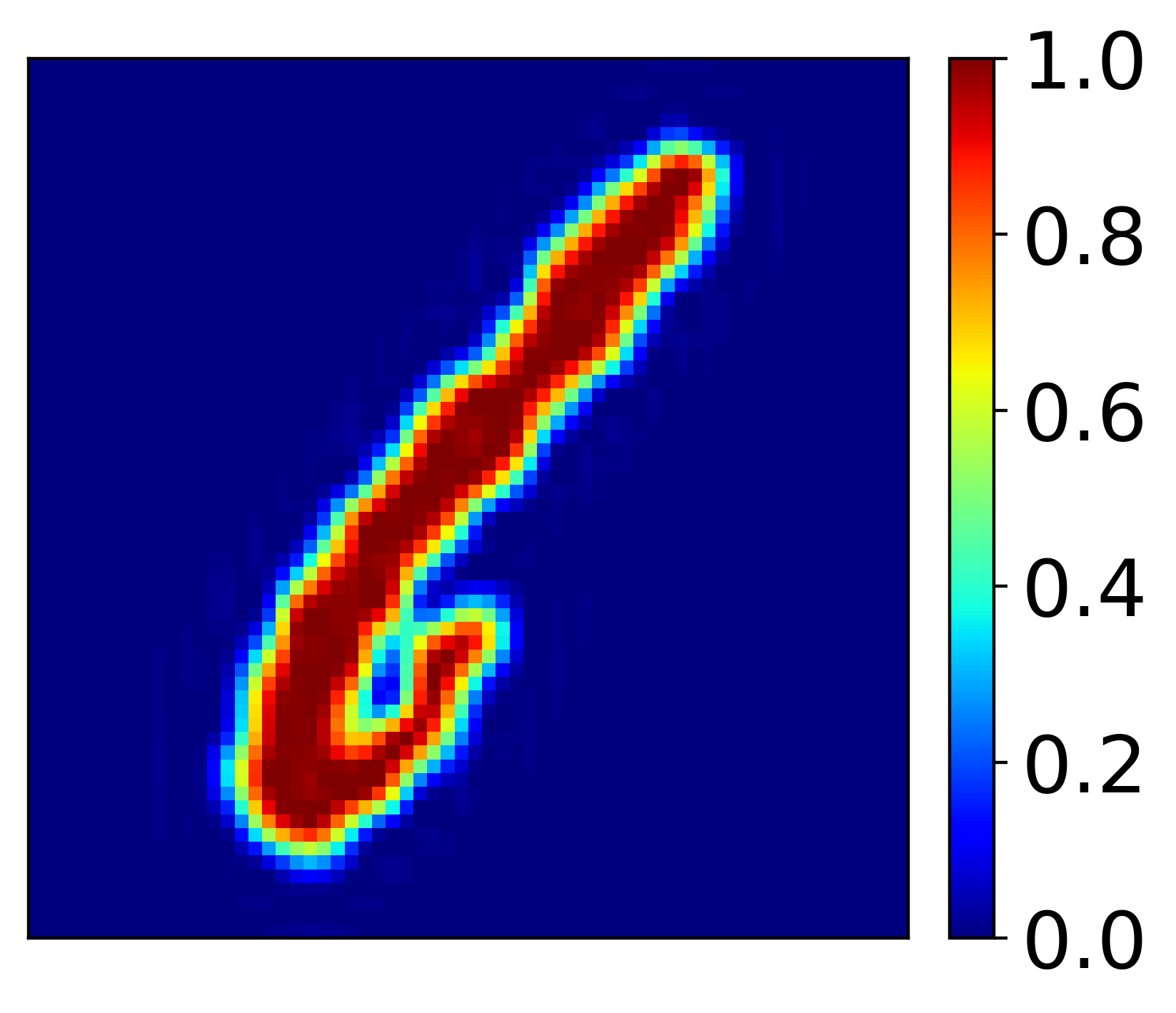}
	}
	\hspace{-0.4cm}
	\subcaptionbox{}{
		\includegraphics[width=0.19\linewidth]{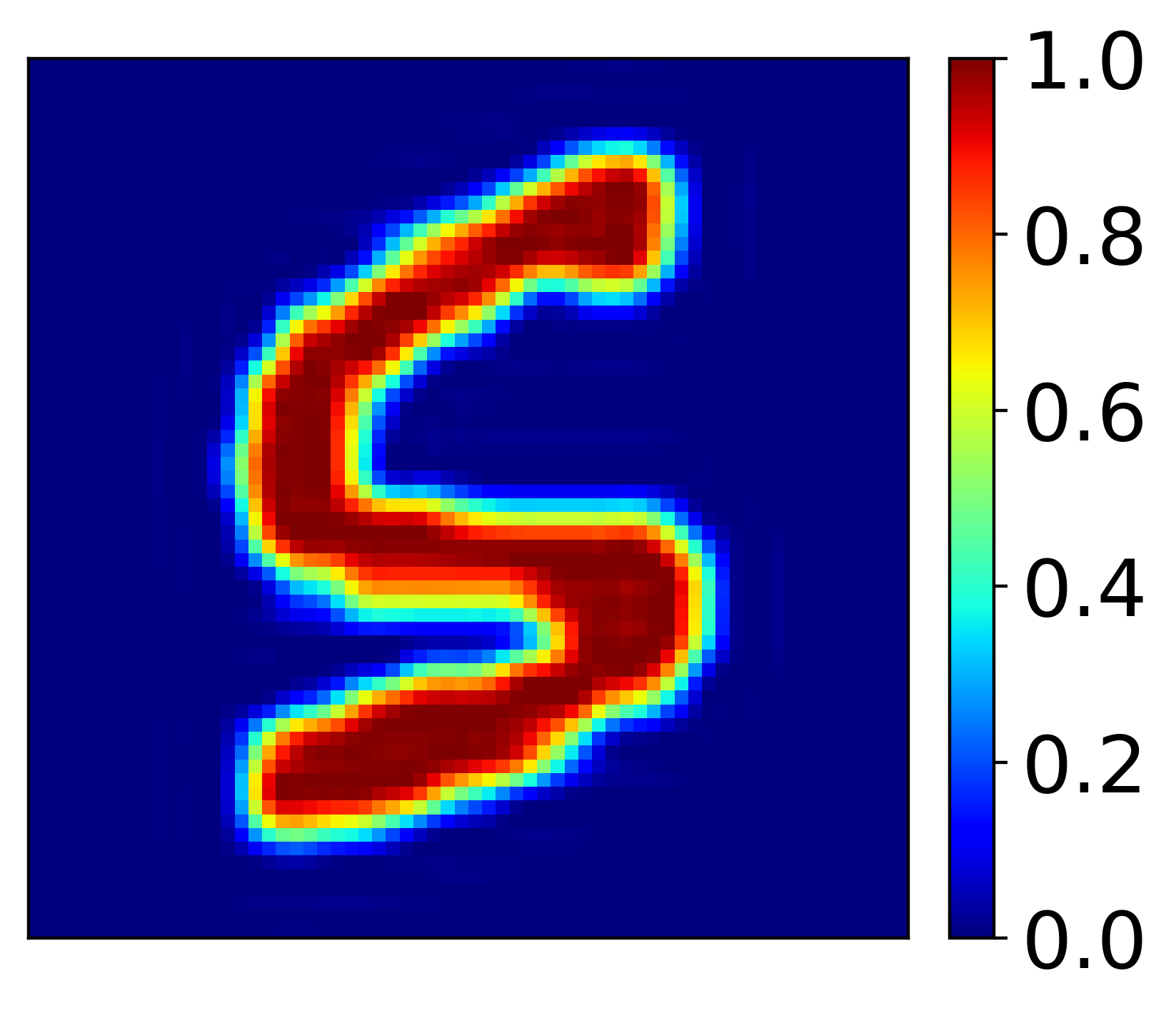}
	}
	\hspace{-0.4cm}
	\subcaptionbox{}{
		\includegraphics[width=0.19\linewidth]{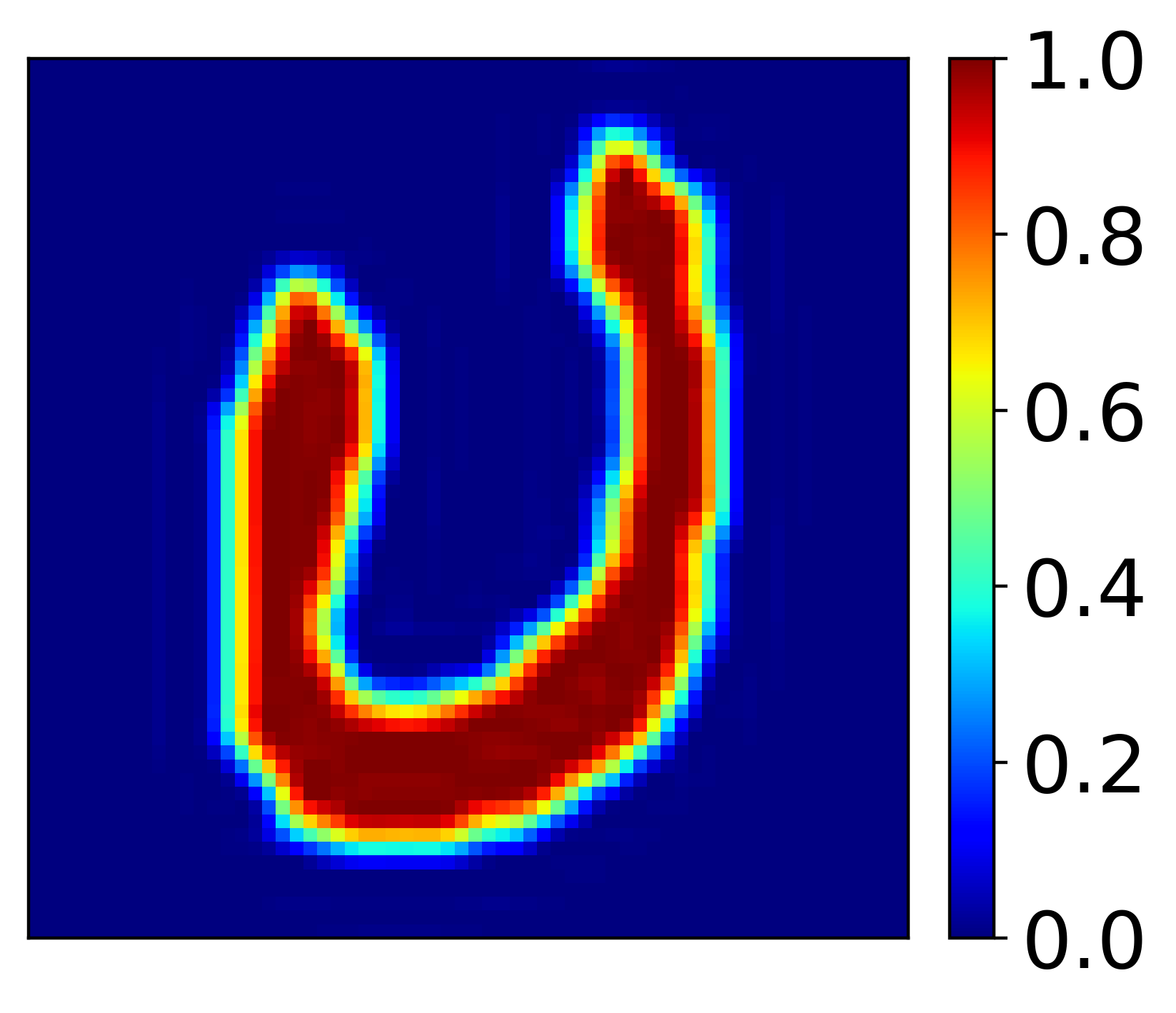}
	}\\%
	\subcaptionbox{}{
		\includegraphics[width=0.19\linewidth]{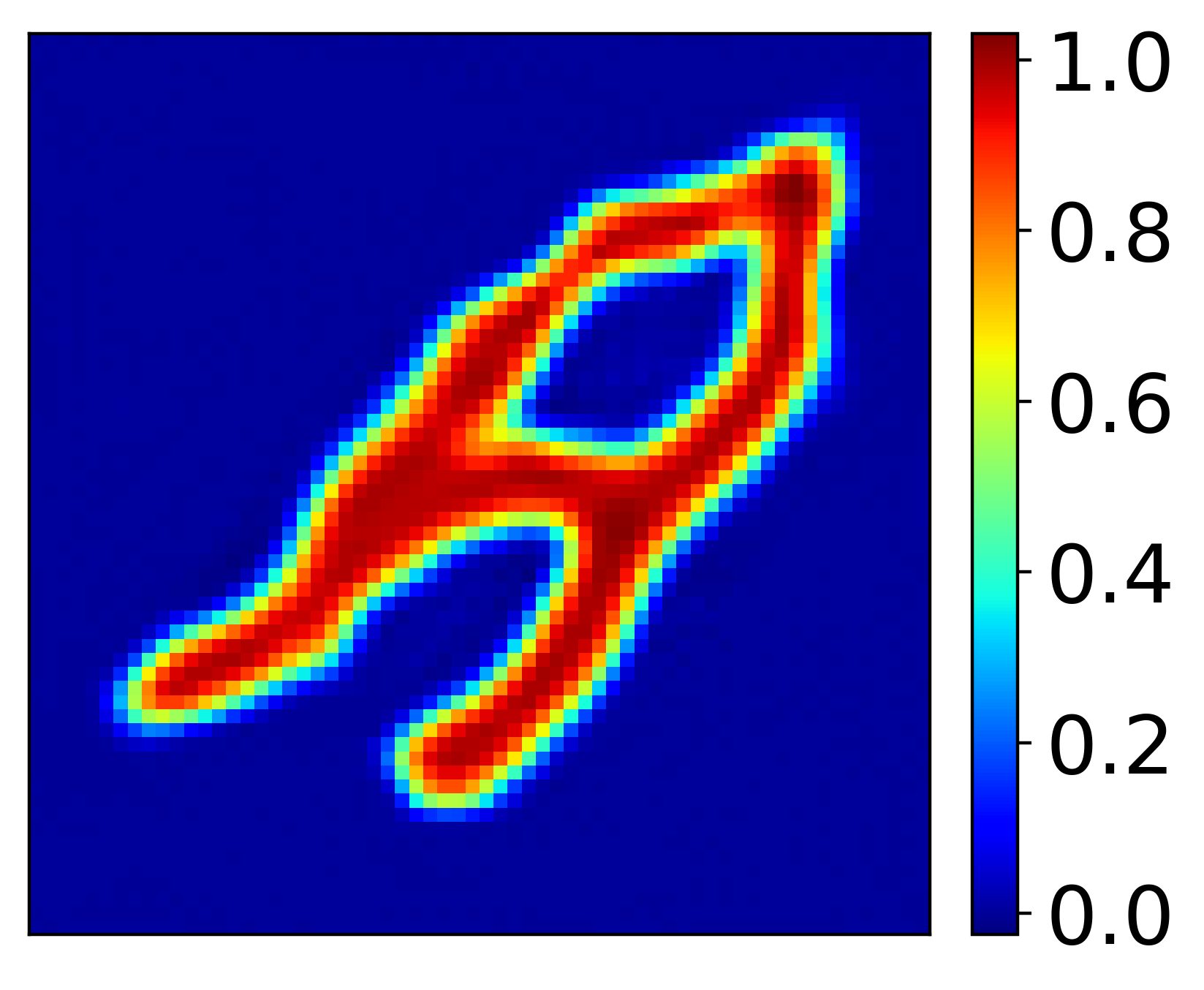}
	}
	\hspace{-0.4cm}
	\subcaptionbox{}{
		\includegraphics[width=0.19\linewidth]{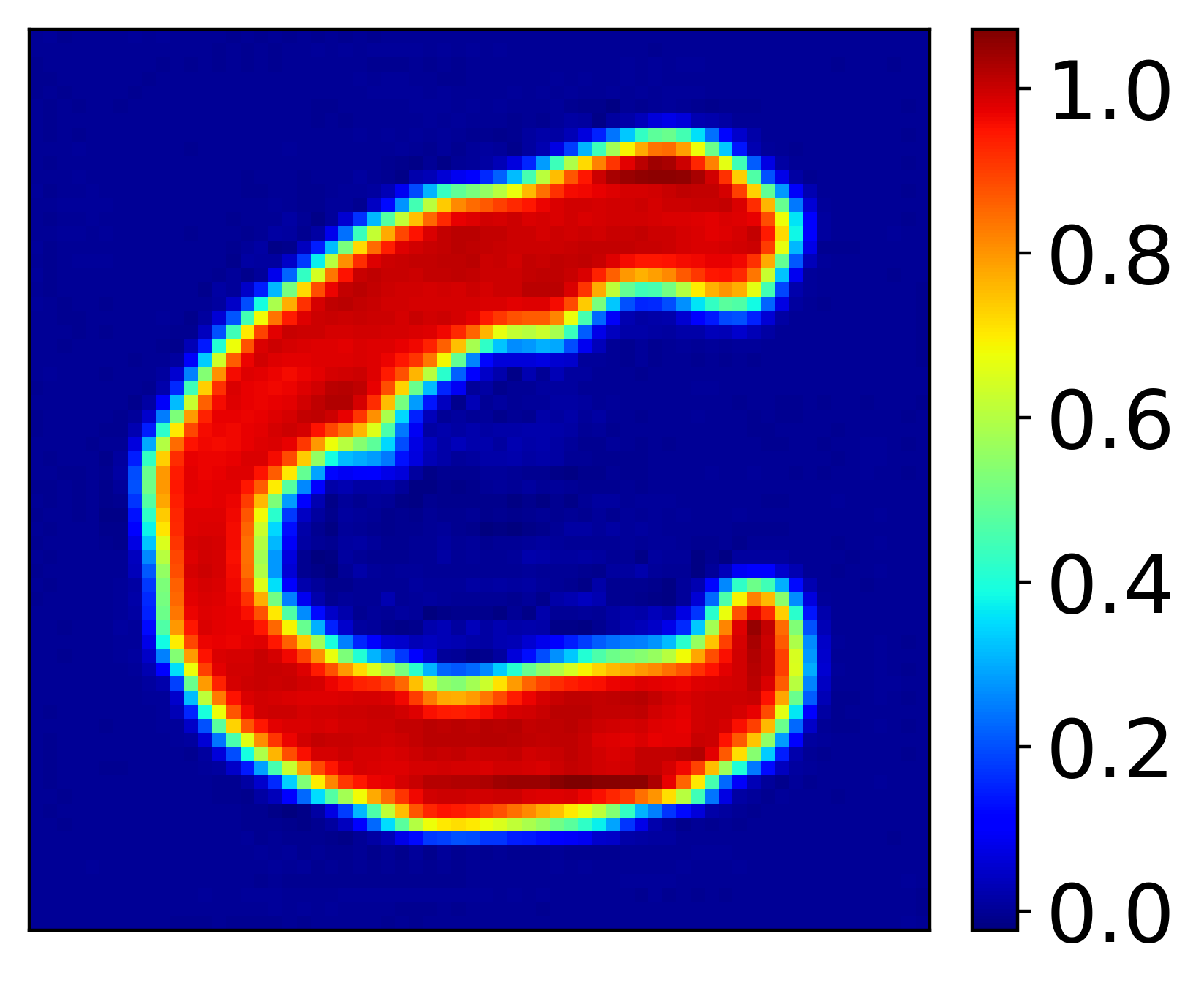}
	}
	\hspace{-0.4cm}
	\subcaptionbox{}{
		\includegraphics[width=0.19\linewidth]{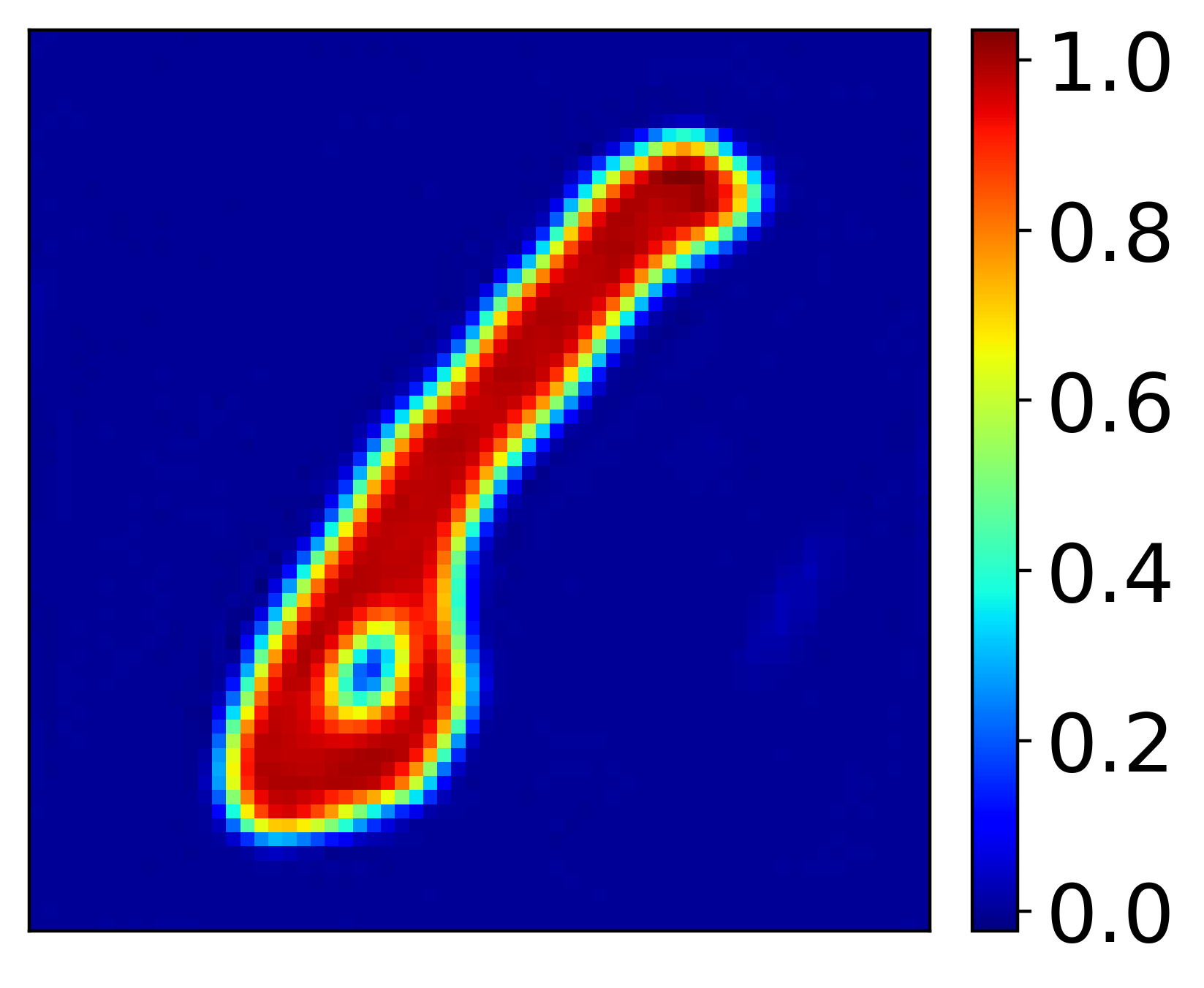}
	}
	\hspace{-0.4cm}
	\subcaptionbox{}{
		\includegraphics[width=0.19\linewidth]{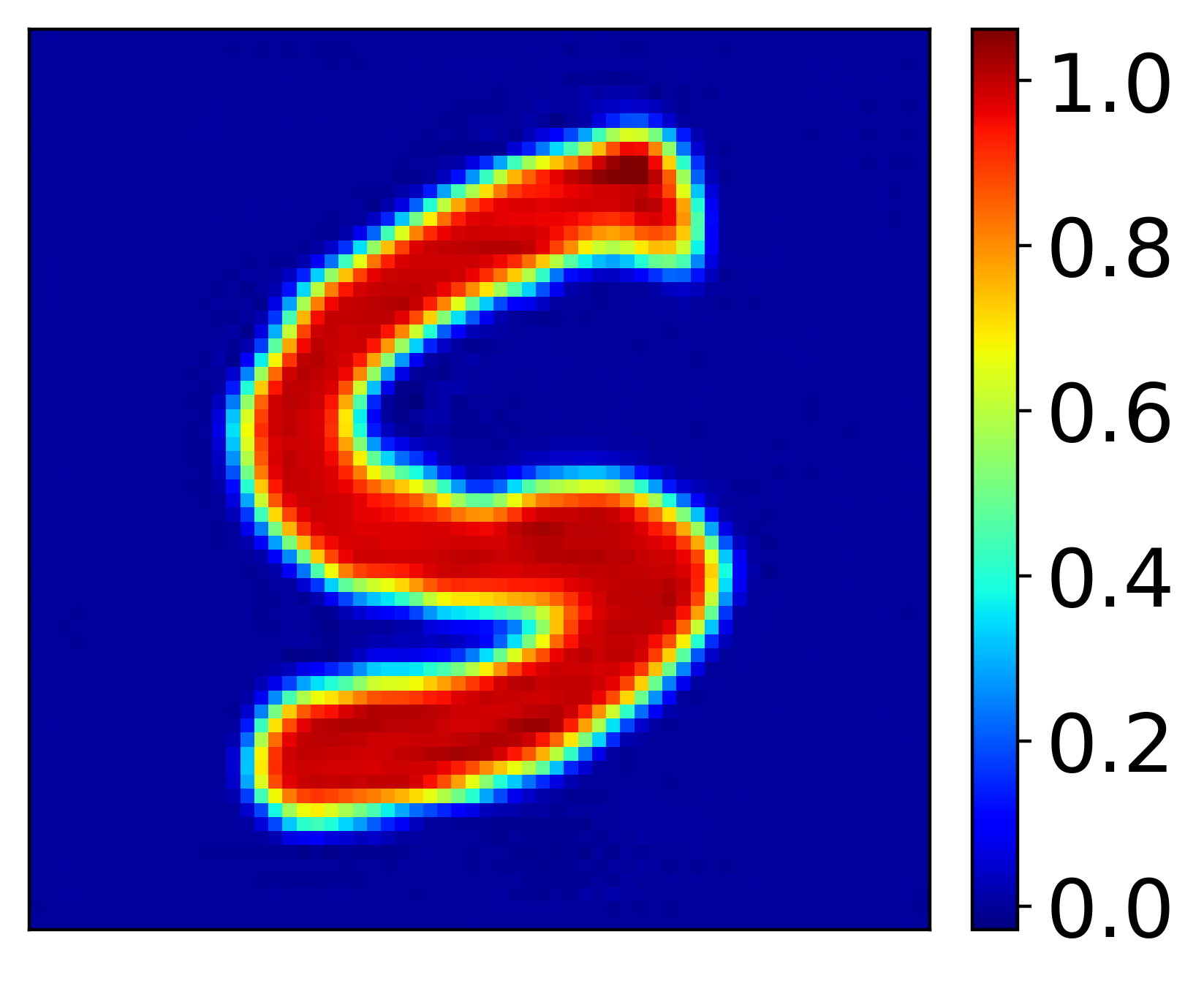}
	}
	\hspace{-0.4cm}
	\subcaptionbox{}{
		\includegraphics[width=0.19\linewidth]{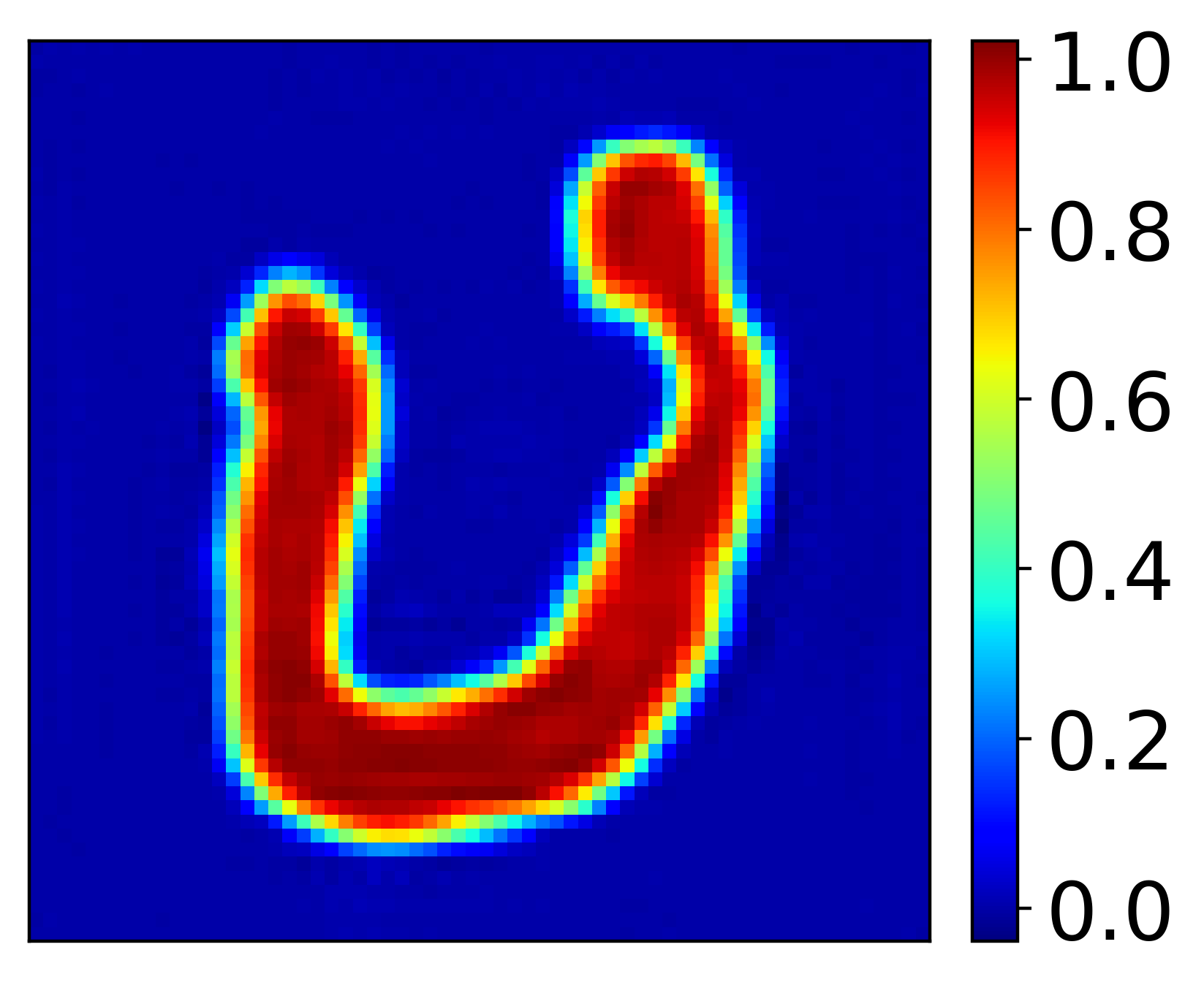}
	}\\
	\caption{Reconstruction results for the `Letter' dataset under $50\%$ noise. The model was trained on MNIST digits, highlighting its cross-domain generalization. Rows from top to bottom: ground-truth sources and U-Net-enhanced reconstructions.}
	\label{fig:letter_5}
\end{figure}

\begin{figure}[htbp]
	\centering
	\subcaptionbox{}{
		\includegraphics[width=0.3\linewidth]{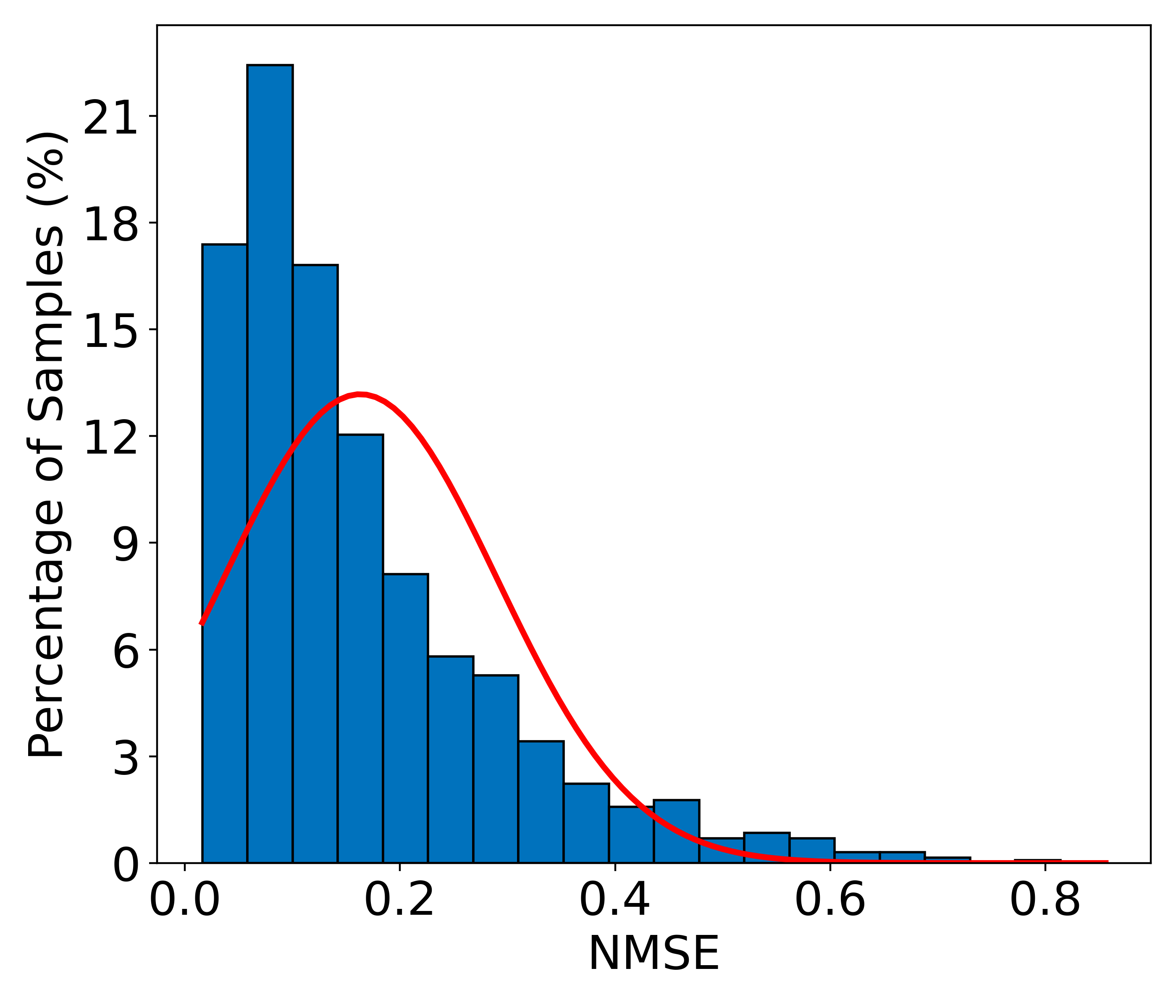}
	}
	\hspace{0.6cm}
	\subcaptionbox{}{
		\includegraphics[width=0.3\linewidth]{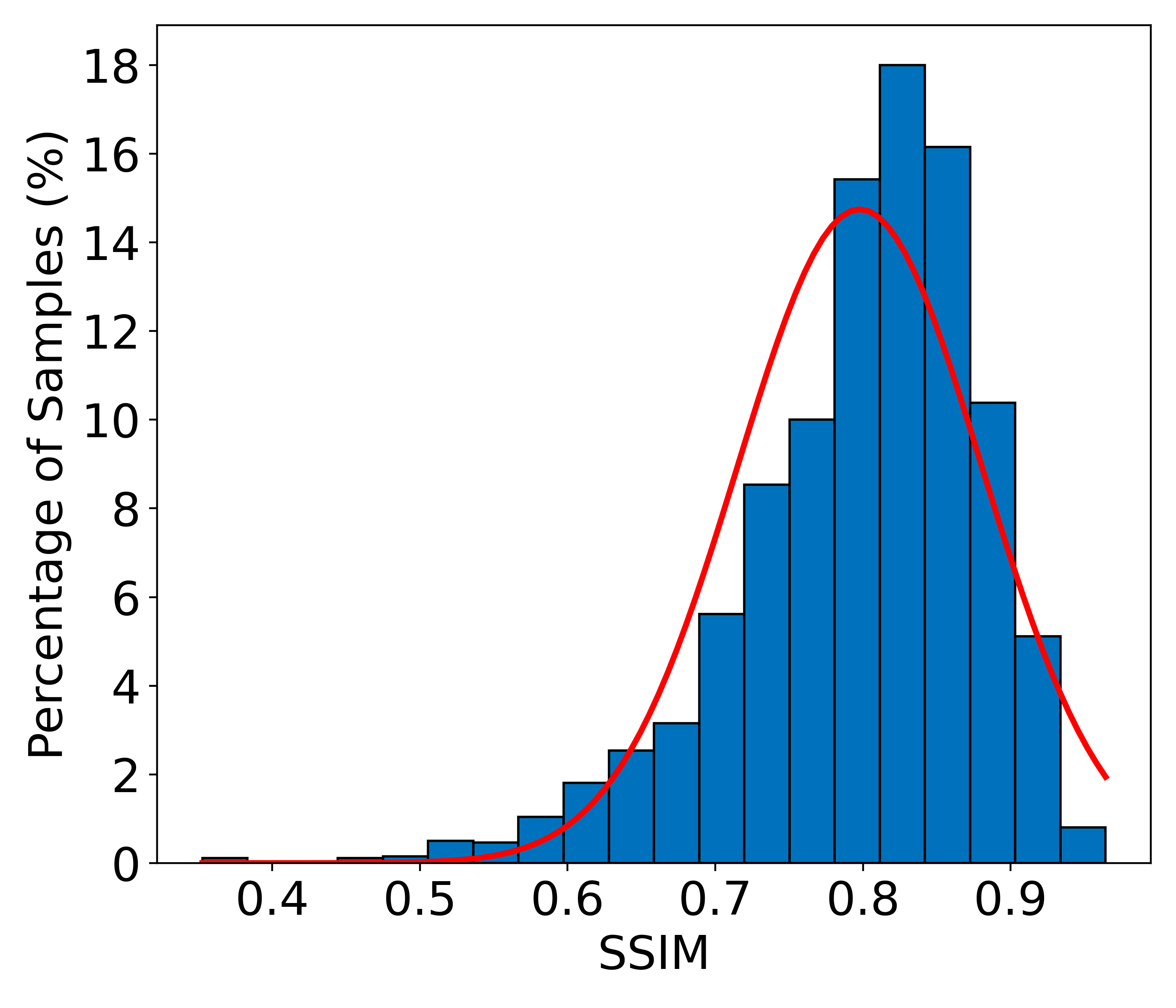}
	}
	\caption{Statistical distribution of NMSE and SSIM for the \enquote{Letter} dataset under $50\%$ noise. The fitting indicates consistent performance on out-of-distribution samples.}
	\label{fig:nmse_letter}
\end{figure}
To assess the out-of-distribution generalization capabilities of the proposed method, we evaluate the model on the \enquote{Letter} dataset, which comprises images of the English alphabet.
These characters exhibit different geometric topologies and edge distributions compared to the numerical digits used during training.
Qualitative results for representative samples under $50\%$ noise levels are presented in \Cref{fig:letter_5}.
Despite the distribution shift and the presence of significant noise, the network preserves the intricate topological features of the letters based on the highly limited $N=2$ initial inputs.
This demonstrates that the model has learned structural priors that generalize well even under challenging conditions where the initial physics-based estimates are poor.
Quantitative performance, summarized in \Cref{fig:nmse_letter}, corroborates this observation, showing consistently low NMSE and high SSIM scores across the dataset.

The numerical experiments confirm the accuracy and efficiency of the method under various conditions. Once trained, the proposed approach requires only a single low-cost Fourier reconstruction, followed by a forward pass through the network. This reduces the computational cost compared to traditional iterative solvers, enabling rapid reconstructions even under sparse data regimes. Thus, our framework offers a robust and computationally efficient solution for inverse source scattering problems with limited data.

\section{Conclusions}

This paper presents a deep-learning-enhanced Fourier method designed to mitigate the ill-posedness and high data requirements of the inverse source scattering problem. By integrating a U-Net architecture with the classical Fourier method, our approach learns an image-to-image mapping that transforms artifact-prone initial reconstructions into high-fidelity source representations. Extensive numerical experiments demonstrate that the proposed hybrid framework achieves accurate and stable reconstructions using sparse, low-frequency data and performs well under extreme noise levels. The model also generalizes across out-of-distribution source geometries and benefits from a high-to-low transfer learning strategy that accelerates convergence for varied noise regimes.

By bridging the physical interpretability of spectral methods with the nonlinear approximation capabilities of deep neural networks, this framework provides a viable approach for inverse problems involving incomplete and noisy measurements. Future work will explore quantitative numerical error analysis, extensions to three-dimensional configurations, and broader classes of wave propagation models.

\section*{Acknowledgments}
Yan Chang was supported by the National Natural Science Foundation of China (NSFC; Grant Nos. 12501584 and 12571452) and the China Postdoctoral Science Foundation (Grant No. GZC20252762). Yukun Guo was supported by the NSFC (Grant Nos. 11971133 and 12571452). Yuliang Wang was supported by the Guangdong and Hong Kong Universities \enquote*{1+1+1} Joint Research Collaboration Scheme (Grant No. 2025A0505000007) and the General Program of the Guangdong Basic and Applied Basic Research Foundation (Grant No. 2025A1515010201).

\bibliographystyle{unsrt}
\bibliography{reference.bib}

\end{document}